\documentclass[fleqn,a4paper,12pt]{article}
\usepackage{amssymb,amsfonts,amsmath,amsthm,bbm,graphicx}

\usepackage[dvips]{epsfig}

\numberwithin{equation}{section}
\numberwithin{figure}{section}

\newcommand{\bs}{\boldsymbol}

\begin{document}


\begin{center}



{\textbf{ASPHERICITY OF A LENGTH FOUR RELATIVE GROUP PRESENTATION}}

\medskip Abd Ghafur Bin Ahmad, Muna A Al-Mulla and Martin Edjvet

\end{center}

\medskip

\begin{abstract}
We consider the relative group presentation $\mathcal{P} = \langle G, \bs{x} | \bs{r} \rangle$ where $\bs{x} = \{ x \}$ and $\bs{r} = \{ xg_1 xg_2 xg_3 x^{-1} g_4 \}$.  We show modulo a small number of exceptional cases exactly when $\mathcal{P}$ is aspherical.  If $H = \langle g_1^{-1} g_2, g_1^{-1} g_3 g_1 , g_4 \rangle \leq G$ then the exceptional cases occur when $H$ is isomorphic to one of $C_5,C_6,C_8$ or $C_2 \times C_4$.\newline
\end{abstract}

\section{Introduction}

A \textit{relative group presentation} is a presentation of the form $\mathcal{P} = \langle G, \bs{x} | \bs{r} \rangle$ where $G$ is a group,
$\bs{x}$ a set disjoint from $G$ and $\bs{r}$ is a set of cyclically reduced words in the free product $G \ast \langle \bs{x} \rangle$ where
$\langle \bs{x} \rangle$ denotes the free group on $\bs{x}$.  If $G(\mathcal{P})$ denotes the group defined by $\mathcal{P}$ then $G(\mathcal{P})$ is the quotient group $G \ast \langle \bs{x} \rangle / N$ where $N$ denotes the normal closure in $G \ast \langle \bs{x} \rangle$ of $\bs{r}$. A relative presentation is defined in [2] to be \textit{aspherical} if every spherical picture over it contains a \textit{dipole}.  If $\mathcal{P}$ is aspherical then statements about
$G(\mathcal{P})$ can be deduced and the reader is referred to [2] for a discussion of these; in particular torsion in $G(\mathcal{P})$ can easily be described.

We will consider the case when both $\bs{x}$ and $\bs{r}$ consists of a single element.  Thus $\bs{r}=\{ r \}$ where
$r=x^{\varepsilon_1} g_1 \ldots x^{\varepsilon_k} g_k$ where $g_i \in G$, $\varepsilon _i = \pm 1$ and $g_i = 1$ implies
$\varepsilon_i + \varepsilon_{i+1} \neq 0$ for $1 \leq i \leq k$, subscripts \textrm{mod}$\, k$.  If $k \leq 3$ or if
$r \in \{ xg_1 xg_2 xg_3 xg_4 , xg_1 xg_2 xg_3 xg_4 xg_5 \}$ then, modulo some open cases, a complete classification of when $\mathcal{P}$ is
aspherical has been obtained in [1], [2], [7] and [11].
The case $r=(xg_1)^p (xg_2)^q (xg_3)^r$ for $p,q,r > 1$ was studied in [12] and $r=x^n g_1 x^{-1} g_2$ ($n \geq 4$) was studied in [5].
The authors of [10] used results from [7] in which $r=xg_1 xg_2 x^{-1} g_3$ to prove asphericity for certain LOG groups.
In this paper we continue the study of asphericity and consider $r=xg_1 xg_2 xg_3 x^{-1} g_4$.
Observe that $r=1$ if and only if $x^{-1} g_2^{-1} x^{-1} g_1^{-1} x^{-1} g_4^{-1} xg_3^{-1} =1$ so replacing $x^{-1}$ by $x$ it follows that we can work modulo $g_1 \leftrightarrow g_2^{-1}$ and $g_3 \leftrightarrow g_4^{-1}$.
A standard approach is to make the transformation $t=xg_1$ and then consider the subgroup $H$ of $G$ generated by the resulting coefficients.
In our case $r$ becomes $t^2 g_1^{-1} g_2 t g_1^{-1} g_3 g_1 t^{-1} g_4$ and so $H = \langle g_1^{-1} g_2, g_1^{-1} g_3 g_1, g_4 \rangle$.
One then usually proceeds according to either when $H$ is non-cyclic or when $H$ is cyclic.
(Note that $\langle g_1^{-1} g_2, g_1^{-1} g_3 g_1, g_4 \rangle$ is cyclic if and only if $\langle g_2 g_1^{-1},g_2 g_4^{-1} g_2^{-1}, g_3^{-1} \rangle$ is cyclic.)
The latter case seems to be the more complicated -- indeed the open cases referred to in the above paragraph almost all involve $H$ being cyclic.  Our results reflect this difference in difficulty.  When $H$ is non-cyclic we obtain a complete answer except for the following case (modulo $g_1 \leftrightarrow g_2^{-1}$, $g_3 \leftrightarrow g_4^{-1}$) in which $H \cong C_2 \times C_4$:
\begin{enumerate}
\item[(\textbf{E})]
$|g_3|=2$; $|g_4|=4$; $g_1^{-1} g_2=1$; $g_1^{-1} g_3 g_1 g_4 = g_4 g_1^{-1} g_3 g_1$.
\end{enumerate}
\textbf{Theorem 1.1}\quad \textit{Let $\mathcal{P}$ be the relative presentation}
\[
\mathcal{P} = \langle G, x | xg_1 xg_2 xg_3 x^{-1} g_4 \rangle ,
\]
\textit{where $g_i \in G$ ($1 \leq i \leq 4$), $g_3 \neq 1$, $g_4 \neq 1$ and $x \notin G$.  Let $H=\langle g_1^{-1} g_2, g_1^{-1} g_3 g_1 , g_4 \rangle$ and assume that $H$ is non-cyclic and the exceptional case (\textbf{E}) does not hold.  Then $\mathcal{P}$ is aspherical if and only if (modulo $g_1 \leftrightarrow g_2^{-1}$, $g_3 \leftrightarrow g_4^{-1}$) none of the following conditions holds:}
\begin{enumerate}
\item[(i)]
$| g_4 | < \infty$ \textit{and $g_3 g_1 g_2^{-1} = 1$;}
\item[(ii)]
$| g_1^{-1} g_2 | < \infty$, $| g_3 | = | g_4 | = 2$ \textit{and $g_1^{-1} g_3 g_1 g_4 = g_2 g_4^{-1} g_2^{-1} g_3^{-1}=1$;}
\item[(iii)]
$\dfrac{1}{|g_1^{-1} g_2|} + \dfrac{1}{| g_3 |} + \dfrac{1}{| g_4 |} + \dfrac{1}{| g_2 g_4 g_1^{-1} g_3^{-1} |} > 2$.
\end{enumerate}

Now let $H$ be cyclic.  Before stating the theorem we make a list of exceptions (modulo $g_1 \leftrightarrow g_2^{-1}$, $g_3 \leftrightarrow g_4^{-1}$).
\begin{enumerate}
\item[(\textbf{E1})]
$| g_4 | = 5$; $g_1^{-1} g_2 = g_4^2$; $g_1^{-1} g_3 g_1 = g_4^3$.
\item[(\textbf{E2})]
$|g_4|=6$; $g_1^{-1} g_2=1$; $g_1^{-1} g_3 g_1 = g_4^2$.
\item[(\textbf{E3})]
$|g_4|=6$; $g_1^{-1} g_2=1$; $g_1^{-1} g_3 g_1 = g_4^4$.
\item[(\textbf{E4})]
$| g_4 | = 8$; $g_1^{-1} g_2 = 1$; $g_1^{-1} g_3 g_1 = g_4^4$.
\end{enumerate}
Observe that (\textbf{E1}) implies $H \cong C_5$; (\textbf{E2}) and (\textbf{E3}) imply $H \cong C_6$; and (\textbf{E4}) implies $H \cong C_8$.

\textbf{Theorem 1.2}\quad \textit{Let $\mathcal{P}$ be the relative presentation}
\[
\mathcal{P}= \langle G,x | xg_1 xg_2 xg_3 x^{-1} g_4 \rangle
\]
\textit{where $g_i \in G$ ($1 \leq i \leq 4$), $g_3 \neq 1$, $g_4 \neq 1$ and $x \notin G$.  Let
$H=\langle g_1^{-1} g_2, g_1^{-1} g_3 g_1, g_4 \rangle$ be a cyclic group.  Suppose that none of the exceptional conditions (\textbf{E1})--(\textbf{E4}) holds. Then $\mathcal{P}$ is aspherical if and only if either $H$ is infinite or $H$ is finite and (modulo $g_1 \leftrightarrow g_2^{-1}$, $g_3 \leftrightarrow g_4^{-1}$) none of the following conditions holds:}

\begin{tabular}{rlrl}
(i)&$g_3 g_1 g_2^{-1}=1$\textit{;}&(vi)&$|g_3|=2$\textit{;} $g_1^{-1} g_3 g_2 g_4 = (g_1^{-1} g_2)^2 g_4^{-1} =1$\textit{;}\\
(ii)&$g_3^{-1} g_1 g_2^{-1} = g_2 g_4^{-1} g_1^{-1} g_3^{-1}=1$\textit{;}&(vii)&$g_1^{-1} g_2 = 1$\textit{;} $g_1^{-1} g_3 g_1 g_4^{\pm 1} =1$\textit{;}\\
(iii)&$g_3^{-1} g_1 g_2^{-1} = g_4 g_2^{-1} g_1 =1$\textit{;}&(viii)&$g_1^{-1} g_2=1$\textit{;} $|g_3|=2$\textit{;} $|g_4|=3$\textit{;}\\
(iv)&$|g_3|=2$\textit{;} $|g_4|=2$\textit{;}&(ix)&$g_1^{-1} g_2=1$\textit{;} $4 \leq | g_3 | \leq 5$\textit{;} $g_1^{-1} g_3^2 g_1 g_4$\textit{;}\\
(v)&$|g_3|=2$\textit{;} $g_1^{-1} g_3 g_2 g_4 = g_1^{-1} g_2 g_4^{-2}=1$\textit{;}&(x)&
$g_1^{-1} g_2=1$\textit{;} $| g_3 | = 6$\textit{;} $g_1^{-1} g_3^3 g_1 g_4$.
\end{tabular}

In Section 2 we discuss pictures and curvature; in Section 3 there are some preliminary results; Theorem 1.1 and Theorem 1.2 are proved in Sections 4 and 5.

Acknowledgement: the authors thank R M Thomas for his helpful comments concerning this paper.


\section{Pictures} 

The definitions of this section are taken from [2].  The reader should consult [2] and [1] for more details.

A \textit{picture} $\bs{P}$ is a finite collection of pairwise disjoint discs $\{ \Delta_1, \ldots, \Delta_m \}$ in the interior of a disc $D^2$, together with a finite collection of pairwise disjoint simple arcs $\{ \alpha_1 , \ldots, \alpha _n \}$ embedded in the closure of
$D^2 - \bigcup_{i=1}^m \Delta_i$ in such a way that each arc meets $\partial D^2 \cup \bigcup_{i=1}^m \Delta_i$ transversely in its end points.
The \textit{boundary} of $\bs{P}$ is the circle $\partial D^2$, denoted $\partial \bs{P}$.  For $1 \leq i \leq m$, the \textit{corners} of $\Delta_i$ are the closures of the connected components of $\partial \Delta_i - \bigcup_{j=1}^n \alpha_j$, where $\partial \Delta_i$ is the boundary of $\Delta_i$.
The \textit{regions} of $\bs{P}$ are the closures of the connected components of
$D^2 - \left( \bigcup_{i=1}^m \Delta_i \cup \bigcup_{j=1}^n \alpha _j \right)$.

An \textit{inner region} of $\bs{P}$ is a simply connected region of $\bs{P}$ that does not meet $\partial \bs{P}$.
The picture $\bs{P}$ is \textit{non-trivial} if $m \geq 1$, is \textit{connected} if
$\bigcup_{i=1}^m \Delta_i \cup \bigcup_{j=1}^n \alpha_j$ is connected, and is \textit{spherical} if it is non-trivial and if none of the arcs meets the boundary of $D^2$.  The number of edges in a region $\Delta$ is called the \textit{degree} of $\Delta$ and is denoted by $d(\Delta)$.  If $\bs{P}$ is a spherical picture, the number of different discs to which a disc $\Delta_i$ is connected is called the \textit{degree} of $\Delta_i$, denoted by $\deg (\Delta_i)$.

With $\mathcal{P} = \langle G, \bs{x} | \bs{r} \rangle$ define the following labelling: each arc $\alpha_j$ is equipped with a normal orientation, indicated by a short arrow meeting the arc transversely, and labelled by an element of $\bs{x} \cup \bs{x}^{-1}$.  Each corner of $\bs{P}$ is oriented \textit{anticlockwise} (with respect to $D^2$) and labelled by an element of $G$.  If $\kappa$ is a corner of a disc $\Delta_i$ of $\bs{P}$, then $W(\kappa)$ is the word obtained by reading in an anticlockwise order the labels on the arcs and corners meeting $\partial \Delta_i$ beginning with the label on the first arc we meet as we read the anticlockwise corner $\kappa$.  If we cross an arc labelled $x$ in the direction of its normal orientation, we read $x$, otherwise we read $x^{-1}$.

A picture $\bs{P}$ is called a \textit{picture over the relative presentation} $\mathcal{P}$ if the above labelling satisfies the following conditions.
\begin{enumerate}
\item[(1)]
For each corner $\kappa$ of $\bs{P}$, $W(\kappa ) \in \bs{r}^{\ast}$, the set of all cyclic permutations of the members of $\bs{r} \cup \bs{r}^{-1}$ which begin with a member of $\bs{x}$.
\item[(2)]
If $g_1, \ldots, g_l$ is the sequence of corner labels encountered in a \textit{clockwise} traversal of the boundary of an inner region $\Delta$ of $\bs{P}$, then the product $g_1 \ldots g_l = 1$ in $G$.  We say that $g_1 \ldots g_l$ is the \textit{label} of $\Delta$.
\end{enumerate}
A \textit{dipole} in a labelled picture $\bs{P}$ over $\mathcal{P}$ consists of corners $\kappa$ and $\kappa '$ of $\bs{P}$ together with an arc joining the two corners such that $\kappa$ and $\kappa '$ belong to the same region and such that if $W(\kappa)=Sg$ where $g \in G$ and $S$ begins and ends with a member of $\bs{x} \cup \bs{x}^{-1}$, then $W(\kappa ')=S^{-1} h^{-1}$.  The picture $\bs{P}$ is \textit{reduced} if it does not contain a dipole.
A relative presentation $\mathcal{P}$ is called \textit{aspherical} if every connected spherical picture over $\mathcal{P}$ contains a dipole.

The \textit{star graph} $\mathcal{P}^{\textrm{st}}$ of a relative presentation $\mathcal{P}$ is a graph whose vertex set is
$\bs{x} \cup \bs{x}^{-1}$ and edge set is $\bs{r}^{\ast}$.  For $R \in \bs{r}^{\ast}$, write $R=Sg$ where $g \in G$ and $S$ begins and ends with a member of $\bs{x} \cup \bs{x}^{-1}$.  The initial and terminal functions are given as follows: $\iota (R)$ is the first symbol of $S$, and $\tau (R)$ is the inverse of the last symbol of $S$.  The labelling function on the edges is defined by $\lambda (R)=g^{-1}$ and is extended to paths in the usual way.
A non-empty cyclically reduced cycle (closed path) in $\mathcal{P}^{\textrm{st}}$ will be called \textit{admissible} if it has trivial label in $G$.
In our case $G$ embeds in $\langle G,t|r \rangle$ [8] so \textit{every} region of a reduced picture over $\mathcal{P}$ supports an admissible cycle in $\mathcal{P}^{\textrm{st}}$.

As described in the introduction we will consider spherical pictures over $\mathcal{P}=\langle G,t|r \rangle$ where
$r=t^2 g_1^{-1} g_2 tg_1^{-1} g_3 g_1 t^{-1} g_4$.
For ease of presentation we introduce the following notation: $a=1$, $b=g_1^{-1} g_2$, $c=g_1^{-1} g_3 g_1$ and $d=g_4$ and consider $tatbtct^{-1} d$. Exception (\textbf{E}) and conditions (i)--(iii) of Theorem 1.1 can then be re-written as
\begin{enumerate}
\item[(\textbf{E})]
$|c|=2$; $|d|=4$; $b=1$; $cd=dc$.
\item[(i)]
$|d| < \infty$ and $cab^{-1} =1$;
\item[(ii)]
$|a^{-1} b| < \infty$, $|c| = |d| = 2$ and $a^{-1} cad=bd^{-1} b^{-1} c^{-1}=1$;
\item[(iii)]
$\dfrac{1}{|a^{-1} b|} + \dfrac{1}{|c|} + \dfrac{1}{|d|} + \dfrac{1}{| bda^{-1} c^{-1}|} > 2$.
\end{enumerate}
The exceptions (\textbf{E1})--(\textbf{E4}) and conditions (i)--(x) of Theorem 1.2 can be rewritten as
\begin{enumerate}
\item[(\textbf{E1})]
$|d|=5$; $b=d^2$; $c=d^3$;
\item[(\textbf{E2})]
$|d|=6$; $b=1$; $c=d^2$;
\item[(\textbf{E3})]
$|d|=6$; $b=1$; $c=d^4$;
\item[(\textbf{E4})]
$|d|=8$; $b=1$; $c=d^4$;
\end{enumerate}

\begin{tabular}{rlrl}
(i)&$cab^{-1}=1$;&(vi)&$|c|=2$; $cbda^{-1}=(a^{-1}b)^2 d^{-1}=1$;\\
(ii)&$c^{-1}ab^{-1}=cadb^{-1}=1$;&(vii)&$a^{-1}b=1$; $cad^{\pm 1} a^{-1} =1$;\\
(iii)&$c^{-1} ab^{-1}=db^{-1} a=1$;&(viii)&$a^{-1}b=1$; $|c|=2$; $|d|=3$;\\
(iv)&$|c|=|d|=2$;&(ix)&$a^{-1}b=1$; $4 \leq |c| \leq 5$; $c^2 ada^{-1}=1$;\\
(v)&$|c|=2$; $cbda^{-1}=a^{-1}bd^{-2}=1$;&(x)&$a^{-1}b=1$; $|c|=6$; $c^3 ada^{-1} =1$.
\end{tabular}

\newpage
\begin{figure}
\begin{center}
\psfig{file=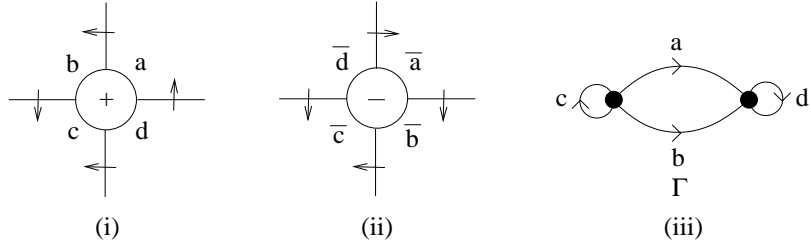}
\end{center}
\caption{vertices and star graph}
\end{figure}

Let $\bs{P}$ be a reduced connected spherical picture over $\mathcal{P}$.  Then the vertices (discs) of $\bs{P}$ are given by Figure 2.1(i) and (ii)
where $\bar{x}$ denotes $x^{-1}$ for $x \in \{a,b,c,d \}$; and the star graph $\Gamma$ is 
given by Figure 2.1(iii).

We make the following assumptions.
\begin{enumerate}
\item[(\textbf{A1})]
$\bs{P}$ has a minimum number of vertices.
\item[(\textbf{A2})]
If $|c|=2$ then, subject to (\textbf{A1}), $\bs{P}$ has a maximum number of regions of degree 2 with label $c^{\pm 2}$.
\end{enumerate}
Observe that (\textbf{A1}) implies that $c^{\varepsilon_1} wc^{\varepsilon_2}$, $d^{\varepsilon_1} w d^{\varepsilon_2}$ where
$\varepsilon_1 =-\varepsilon_2 = \pm 1$ and $w=1$ in $G$ cannot occur as sublabels of a region.  For otherwise a sequence of bridge moves [4] can be applied to produce a dipole which can then be deleted to obtain a picture with fewer vertices.  Moreover if $|c| =2$ then (\textbf{A2}) implies that
$c^{\pm 2}$ cannot be a proper sublabel and $c^{\varepsilon} wc^{\varepsilon}$ where $\varepsilon = \pm 1$, $w=1$ in $G$ cannot be a sublabel of a region in $\bs{P}$.  For otherwise bridge moves can be applied to increase the number of regions labelled $c^{\pm 2}$ while leaving the number of vertices unchanged.

To prove asphericity we adopt the approach of [6].  Let each corner in every region of $\bs{P}$ be given an angle.  The \textit{curvature} of a vertex is defined to be $2 \pi$ less the sum of the angles at that vertex.  The curvature $c(\Delta)$ of a $k$-gonal region $\Delta$ of $\bs{P}$ is the sum of all the angles of the corners of this region less $(k-2)\pi$.  Our method of associating angles is to give each corner at a vertex of degree $d$ an angle $2 \pi /d$.
This way the vertices have zero curvature and we need consider only the regions.  Thus if $\Delta$ is a $k$-gonal region of $\bs{P}$ (a $k$-gon), denoted by $d(\Delta)=k$, and the degree of the vertices of $\Delta$ are $d_i$ ($1 \leq i \leq k$) then
\[
c(\Delta)=c(d_1,d_2,\ldots,d_k)=(2-k)\pi + 2 \pi \sum_{i=1}^k (1/d_i).
\]
In fact since each $d_i=4$ ($1 \leq i \leq k$) we obtain
\[
c(\Delta)=\pi (2- k/2)
\]
so if $d(\Delta) \geq 4$ then $c(\Delta) \leq 0$.

It follows from the fundamental curvature formula that $\sum c(\Delta)=4 \pi$ is where the sum is taken over all the regions $\Delta$ of $\bs{P}$.
Our strategy to show asphericity will be to show that the positive curvature that exists in $\bs{P}$ can be sufficiently compensated by the negative curvature.
To this end, as a first step, we located the regions $\Delta$ of $\bs{P}$ satisfying $c(\Delta) >0$, that is, of positive curvature.  For each such $\Delta$ we distribute all of $c(\Delta)$ to regions $\hat{\Delta}$ near $\Delta$.  
For such regions $\hat{\Delta}$ of $\bs{P}$ define $c^{\ast}(\hat{\Delta})$ to equal $c(\hat{\Delta})$ plus all the positive curvature $\hat{\Delta}$ receives in the distribution procedure mentioned above with the understanding that if
$\hat{\Delta}$ receives no positive curvature then $c^{\ast}(\hat{\Delta})=c(\hat{\Delta})$.  Observe then that the total curvature of $\bs{P}$ is at most
$\sum (c^{\ast}(\hat{\Delta}))$ where the sum is taken over all regions $\hat{\Delta}$ of $D$ that are not positively curved regions.  Therefore to prove $\mathcal{P}$ is aspherical it suffices to show that $c^{\ast}(\hat{\Delta}) \leq 0$ for each $\hat{\Delta}$.

Using the star graph $\Gamma$ of Figure 2.1(iii) we can list the possible labels of regions of small degree (up to cyclic permutation and inversion).
\[
d(\Delta)=2 \Rightarrow l(\Delta) \in S_2 = \{ c^2,d^2,a^{-1} b \}
\]
\[
d(\Delta)=3 \Rightarrow l(\Delta) \in S_3 = \{ c^3, cab^{-1}, c^{-1} ab^{-1}, d^3, db^{-1} a, d^{-1} b^{-1} a \}
\]
Allowing each element in $S_2 \cup S_3$ to be either trivial or non-trivial yields 512 possibilities.
This number can be reduced without any loss as follows.
\begin{enumerate}
\item[1.]
Work modulo \textit{T-equivalence}, that is, $a \leftrightarrow b^{-1}$, $c \leftrightarrow d^{-1}$.
\item[2.]
Delete any combination that implies $c=1$ or $d=1$.
\item[3.]
Delete any combination that yields a contradiction (for example $c^2=1$, $cab^{-1}=1$, $c^{-1}ab^{-1} \neq 1$.
\item[4.]
Delete any combination that yields $|d| < \infty$ and $cab^{-1}=1$ or $|c| < \infty$ and $d^{-1} b^{-1} a=1$ (see Lemma 3.1(i)).
\item[5.]
When $H = \langle b,c,d \rangle$ is cyclic it can be assumed that $H$ is finite (see Lemma 3.4(i)).
\end{enumerate}
It can be readily verified that there remain 23 cases partitioned according
to the existence in $\bs{P}$ of regions of degree 2 and are listed below.

$\textbf{Case A}$: There are no regions of degree two.\\
\\
(\textbf{A0}) $|c|>3$, $|d|>3$, $a^{-1}b\ne1$, $c^{\pm1}ab^{-1}\ne1$, $d^{\pm1}b^{-1}a\ne1$.\\
(\textbf{A1}) $|c|=3$, $|d|>3$, $a^{-1}b\ne1$, $c^{\pm1}ab^{-1}\ne1$, $d^{\pm1}b^{-1}a\ne1$.\\
(\textbf{A2}) $|c|=3$, $|d|=3$, $a^{-1}b\ne1$, $c^{\pm1}ab^{-1}\ne1$, $d^{\pm1}b^{-1}a\ne1$.\\
(\textbf{A3}) $|c|>3$, $|d|>3$, $a^{-1}b\ne1$, $cab^{-1}=1$, $c^{-1}ab^{-1}\ne1$, $d^{\pm1}b^{-1}a\ne1$.\\
(\textbf{A4}) $|c|>3$, $|d|>3$, $a^{-1}b\ne1$, $cab^{-1}\ne1$, $c^{-1}ab^{-1}=1$, $d^{\pm1}b^{-1}a\ne1$.\\
(\textbf{A5}) $|c|=3$, $|d|>3$, $a^{-1}b\ne1$, $cab^{-1}=1$, $c^{-1}ab^{-1}\ne1$, $d^{\pm1}b^{-1}a\ne1$.\\
(\textbf{A6}) $|c|=3$, $|d|>3$, $a^{-1}b\ne1$, $cab^{-1}\ne1$, $c^{-1}ab^{-1}=1$, $d^{\pm1}b^{-1}a\ne1$.\\
(\textbf{A7}) $|c|=3$, $|d|>3$, $a^{-1}b\ne1$, $c^{\pm1}ab^{-1}\ne1$, $db^{-1}a=1$, $d^{-1}b^{-1}a\ne1$.\\
(\textbf{A8}) $|c|=3$, $|d|=3$, $a^{-1}b\ne1$, $cab^{-1}\ne1$, $c^{-1}ab^{-1}=1$, $d^{\pm1}b^{-1}a\ne1$.\\
(\textbf{A9}) $|c|=3$, $|d|=3$, $a^{-1}b\ne1$, $cab^{-1}\ne1$, $c^{-1}ab^{-1}=1$, $db^{-1}a=1$, $d^{-1}b^{-1}a\ne 1$.\\
(\textbf{A10}) $|c|>3$, $|d|>3$, $a^{-1}b\ne1$, $cab^{-1}\ne1$, $c^{-1}ab^{-1}=1$, $db^{-1}a=1$, $d^{-1}b^{-1}a\ne1$.\\
\\
\textbf{Case B}: Regions of degree two are possible.\\
\\
(\textbf{B1}) $|c|=2$, $|d|>3$, $a^{-1}b\ne1$, $c^{\pm1}ab^{-1}\ne1$, $d^{\pm1}b^{-1}a\ne1$.\\
(\textbf{B2}) $|c|=2$, $|d|=2$, $a^{-1}b\ne1$, $c^{\pm1}ab^{-1}\ne1$, $d^{\pm1}b^{-1}a\ne1$.\\
(\textbf{B3}) $|c|=2$, $|d|=3$, $a^{-1}b\ne1$, $c^{\pm1}ab^{-1}\ne1$, $d^{\pm1}b^{-1}a\ne1$.\\
(\textbf{B4}) $|c|=2$, $|d|>3$, $a^{-1}b\ne1$, $c^{\pm1}ab^{-1}=1$, $d^{\pm1}b^{-1}a\ne1$.\\
(\textbf{B5}) $|c|=2$, $|d|>3$, $a^{-1}b\ne1$, $c^{\pm1}ab^{-1}\ne1$, $db^{-1}a=1$, $d^{-1}b^{-1}a\ne 1$.\\
(\textbf{B6}) $|c|=2$, $|d|=3$, $a^{-1}b\ne1$, $c^{\pm1}ab^{-1}\ne1$, $db^{-1}a=1$, $d^{-1}b^{-1}a\ne1$.\\
(\textbf{B7}) $|c|>3$, $|d|>3$, $a^{-1}b=1$, $c^{\pm1}ab^{-1}\ne1$, $d^{\pm1}b^{-1}a\ne1$.\\
(\textbf{B8}) $|c|=2$, $|d|>3$, $a^{-1}b=1$, $c^{\pm1}ab^{-1}\ne1$, $d^{\pm1}b^{-1}a\ne1$.\\
(\textbf{B9}) $|c|=3$, $|d|>3$, $a^{-1}b=1$, $c^{\pm1}ab^{-1}\ne1$, $d^{\pm1}b^{-1}a\ne1$.\\
(\textbf{B10}) $|c|=2$, $|d|=2$, $a^{-1}b=1$, $c^{\pm1}ab^{-1}\ne1$, $d^{\pm1}b^{-1}a\ne1$.\\
(\textbf{B11}) $|c|=2$, $|d|=3$, $a^{-1}b=1$, $c^{\pm1}ab^{-1}\ne1$, $d^{\pm1}b^{-1}a\ne1$.\\
(\textbf{B12}) $|c|=3$, $|d|=3$, $a^{-1}b=1$, $c^{\pm1}ab^{-1}\ne1$, $d^{\pm1}b^{-1}a\ne1$.\\


\section{Preliminary results} 

We first exhibit spherical diagrams corresponding to some of the conditions of Theorems 1.1 and 1.2.
\textbf{Note}: when drawing figures the discs (vertices) will often be represented by points;
the edge arrows shown in Figure 2.1 will be omitted;
and regions with label $c^{\pm 2},d^{\pm 2}$ will be labelled simply by $c^{\pm 1}, d^{\pm 1}$.

\medskip

\textbf{Lemma 3.1}
\begin{enumerate}
\item[(a)]
\textit{If any of the following conditions holds then $\mathcal{P}$ fails to be aspherical:}
\begin{enumerate}
\item[(i)]
\textit{$|d| < \infty$ and $cab^{-1}=1$;}
\item[(ii)]
$|d| < \infty$ \textit{and $c^{-1} ab^{-1} = bd^{-1} a^{-1} c^{-1} = 1$;}
\item[(iii)]
$|d| < \infty$ and $a^{-1}b = cad^{-1} b^{-1} = 1$.
\end{enumerate}

\medskip

\item[(b)]
\textit{If $bda^{-1} c^{-1} =1$ and $(f_1,f_2,f_3)$ is any of the following then $\mathcal{P}$ fails to be aspherical.}
\begin{enumerate}
\item[(i)]
$(2,2, < \infty)$\textit{;}
\item[(ii)]
$(< \infty,2,2,)$\textit{;}
\item[(iii)]
$(2,3,k)$ ($3 \leq k \leq 5$)\textit{;}
\item[(iv)]
$(3,2,l)$ ($4 \leq l \leq 5$)\textit{;}
\item[(v)]
$(k,2,3)$ ($3 \leq k \leq 5$)\textit{;}
\end{enumerate}
\textit{where} $(f_1,f_2,f_3)=(|a^{-1}b|, |c|, |d|)$.

\medskip

\item[(c)]
\textit{If $a^{-1} b=1$ and $(f_1,f_2,f_3)$ is any of the following then $\mathcal{P}$ fails to be aspherical.}
\begin{enumerate}
\item[(i)]
$(2,2, < \infty)$\textit{;}
\item[(ii)]
$(2,< \infty,2)$\textit{;}
\item[(iii)]
$(2,l,3)$ ($4 \leq l \leq 5$)\textit{;}
\item[(iv)]
$(3,k,2)$ ($3 \leq k \leq 5$)\textit{;}
\item[(v)]
$(2,3,k)$ ($3 \leq k \leq 5$)\textit{,}
\end{enumerate}
\textit{where} $(f_1,f_2,f_3)=(|c|,|d|,|bda^{-1} c^{-1}|)$.
\end{enumerate}

\textit{Proof}
(a)(i) If $|d|=n < \infty$ there is a spherical picture over $\mathcal{P}$ consisting of two $n$-gons labelled $d^n, d^{-n}$ together with $2n$ 3-gons labelled $cab^{-1}$.  The sphere for $n=5$ is shown in Figure 3.1(i).

(ii) If $|d|=n < \infty$ then there is a spherical picture over $\mathcal{P}$ consisting of two regions with labels $d^n,d^{-n}$ and between them two layers.  Each layer has $n$ 4-gons and $n$ 3-gons with labels $cadb^{-1}$, 
$cba^{-1}$ and $bd^{-1} a^{-1} c^{-1}$, $c^{-1} ab^{-1}$.  The sphere for $n=6$ is shown in Figure 3.1(ii).

(iii) If $|d| = n < \infty$ then there is a spherical picture over $\mathcal{P}$ consisting of a region labelled $d^n$, a region labelled $c^{-n}$ and a layer of $2n$ regions between them.  The layer consists of $n$ 2-gons labelled 
$a^{-1}b$ and $n$ 4-gons labelled $cad^{-1}b^{-1}$.  The sphere for $n=6$ is given by Figure 3.1(iii).

\medskip

(b)(i) If $|d|=n < \infty$ then there is a spherical picture over $\mathcal{P}$ consisting of two regions with label $d^n$ and between them a layer consisting of $n$ 4-gons with label $(a^{-1} b)^2$, $n$ 2-gons with label $c^{-2}$ 
and $2n$ 4-gons with label $ad^{-1} b^{-1} c$.  The case $n=5$ is shown in Figure 3.1(iv).

(ii) If $|a^{-1} b|=n < \infty$ then there is a spherical picture over $\mathcal{P}$ consisting of two regions with label $(a^{-1} b)^n$ and between them a layer consisting of $n$ 2-gons with label $c^{-2}$, $n$ 2-gons with label 
$d^2$ and $2n$ 4-gons with label $ad^{-1} b^{-1} c$.  The case $n=3$ is shown in Figure 3.1(v).










(iii) If $k=3$ a spherical picture is given by Figure 3.2(i) in which $A=ab^{-1} ab^{-1}$, $B=bda^{-1} c^{-1}$, $C=c^3$ and $D=d^{-3}$.  If $k=4$ a spherical picture is given by Figure 3.2(ii) in which $D=d^{-4}$.  If $k=5$ there 
is a spherical picture consisting of 30 A regions, 60 B regions, 20 C regions and 12 $D=d^{-5}$ regions, half of which is shown in Figure 3.2(iii).

\newpage
\begin{figure}
\begin{center}
\psfig{file=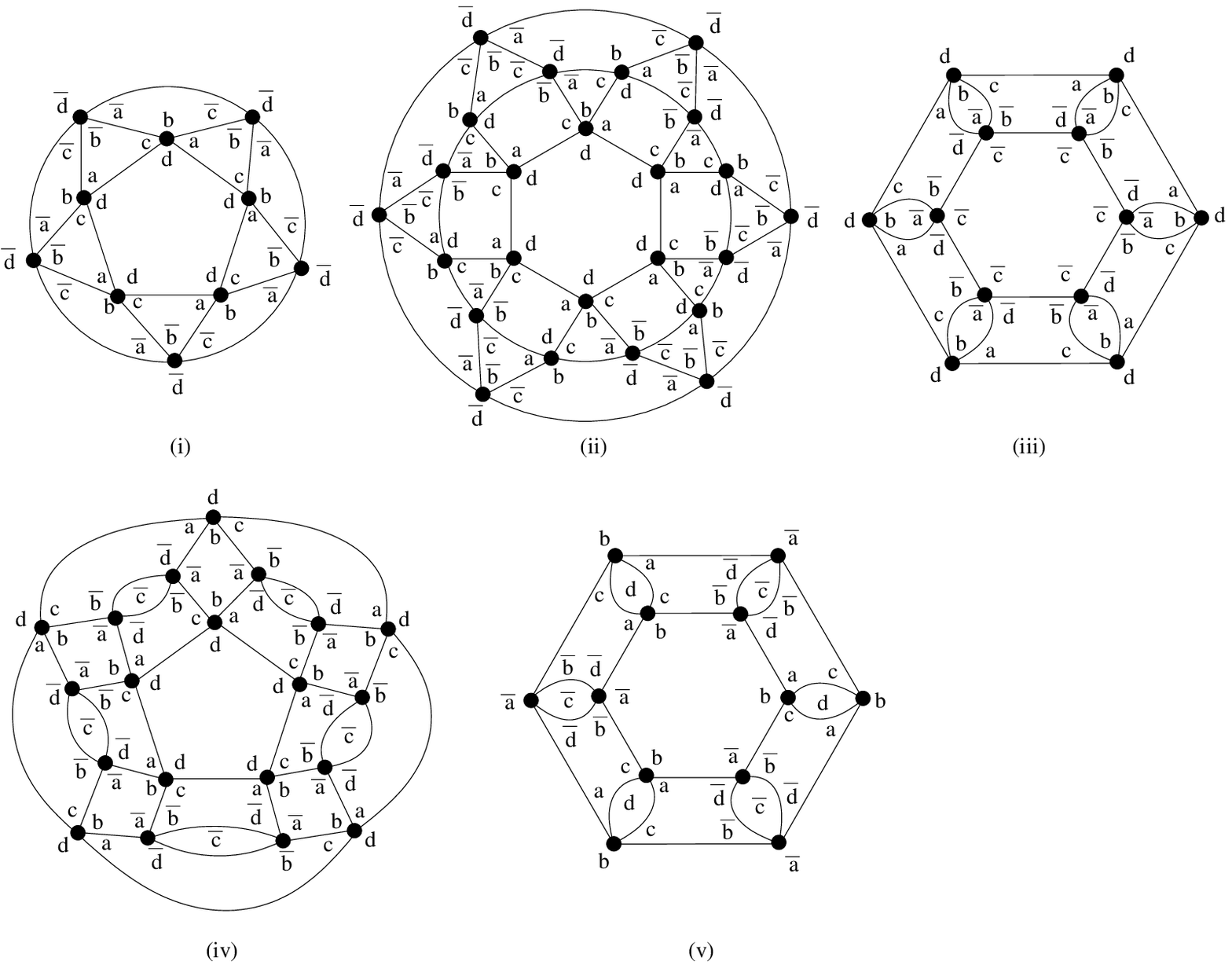}
\end{center}
\caption{spherical pictures}
\end{figure}

\begin{figure}
\begin{center}
\psfig{file=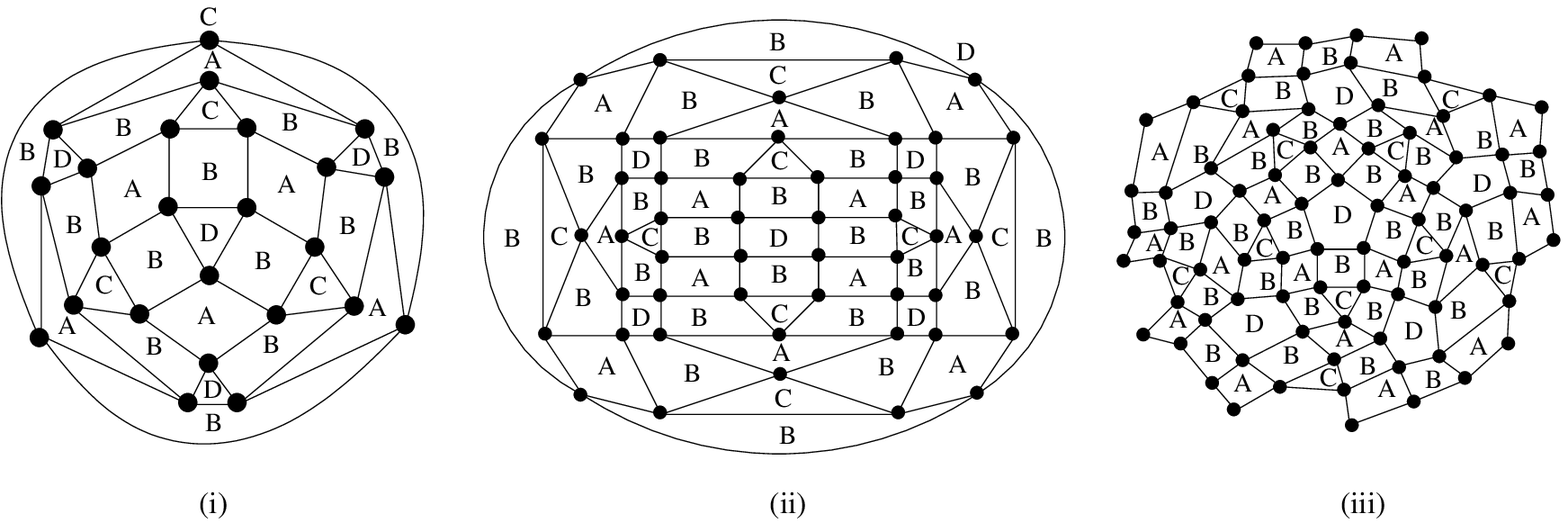}
\end{center}
\caption{spherical pictures}
\end{figure}  

\begin{figure}
\begin{center}
\psfig{file=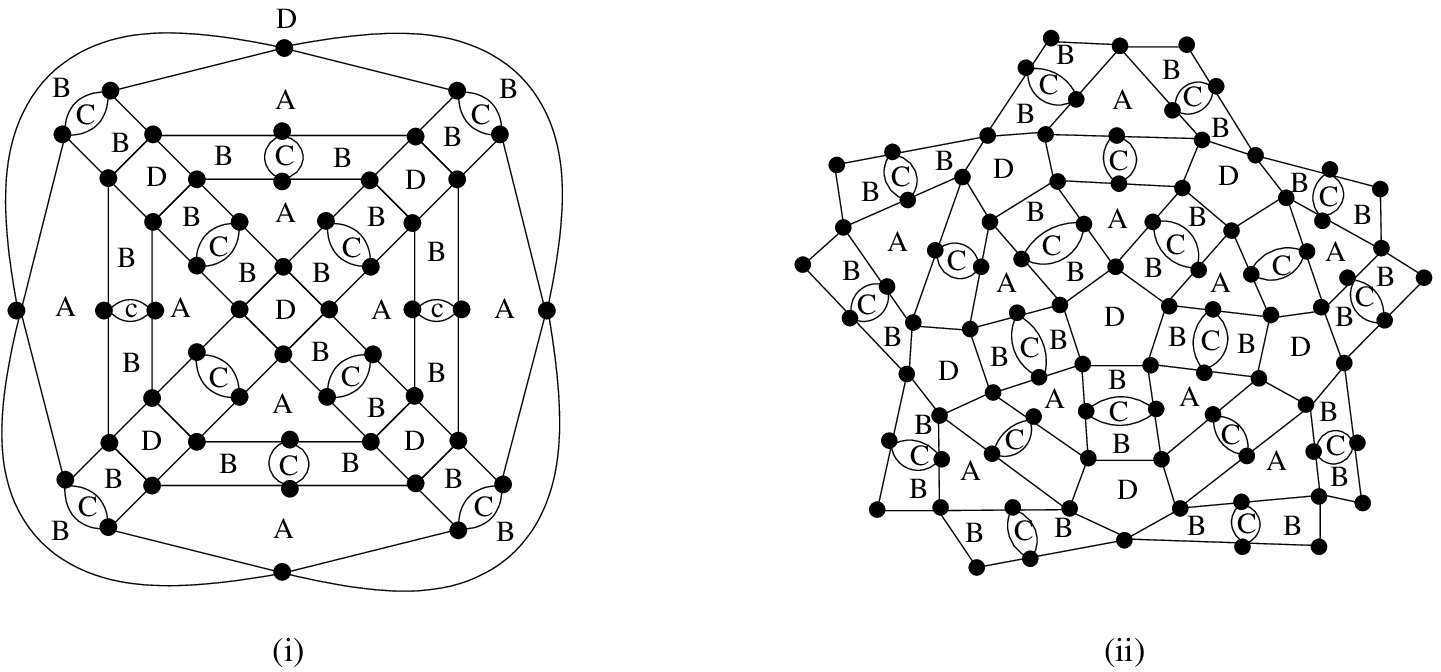}
\end{center}
\caption{spherical pictures}
\end{figure}  

\begin{figure} 
\begin{center}
\psfig{file=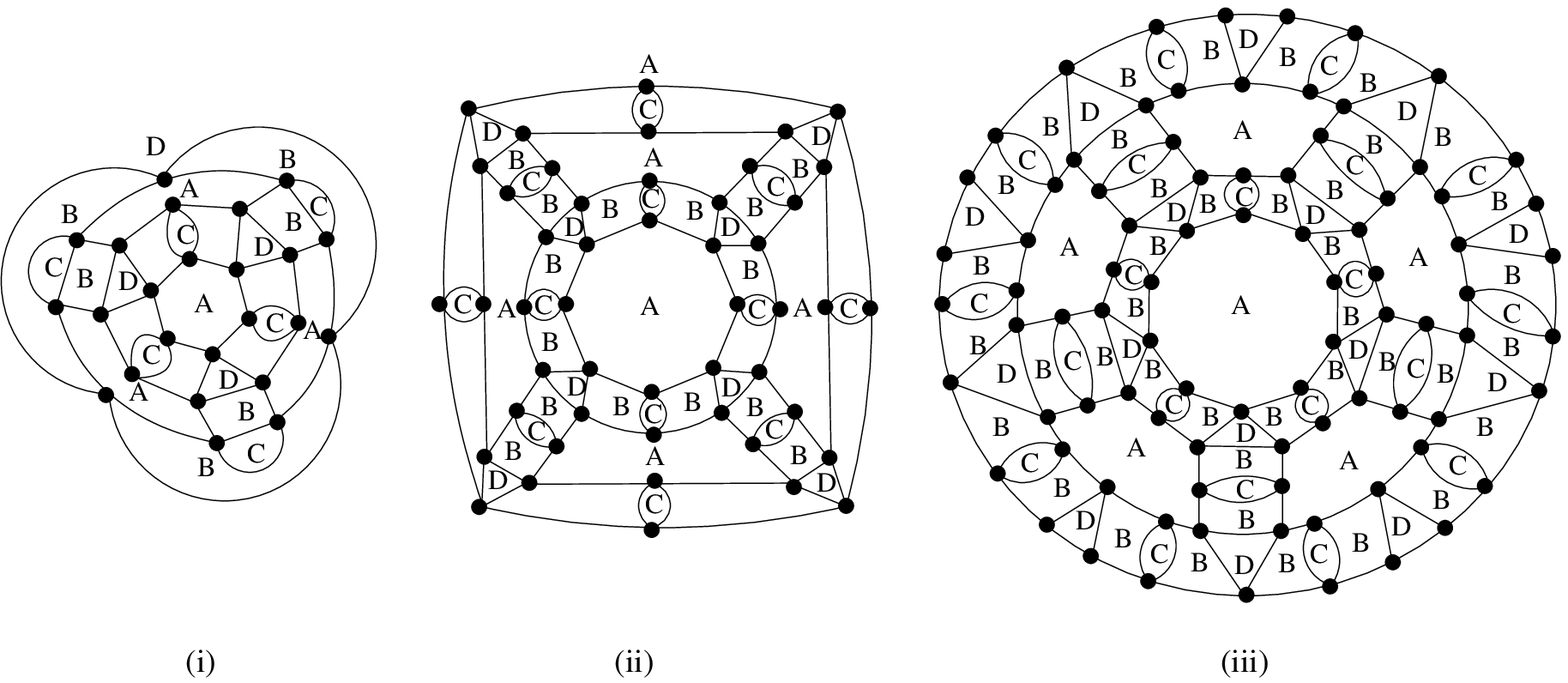}
\end{center}
\caption{spherical pictures}
\end{figure}  

\begin{figure}
\begin{center}
\psfig{file=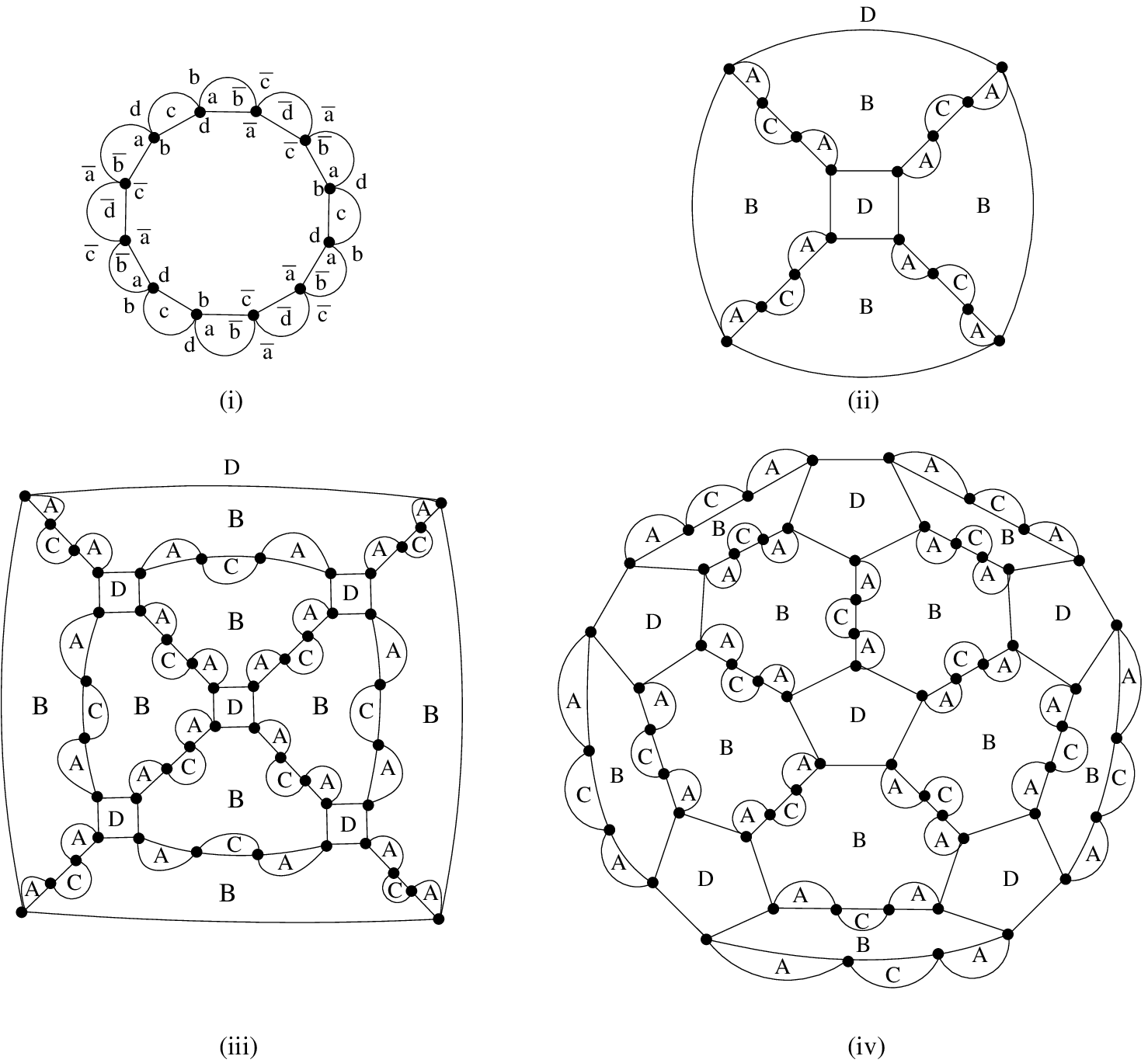}
\end{center}
\caption{spherical pictures}
\end{figure}

\begin{figure} 
\begin{center}
\psfig{file=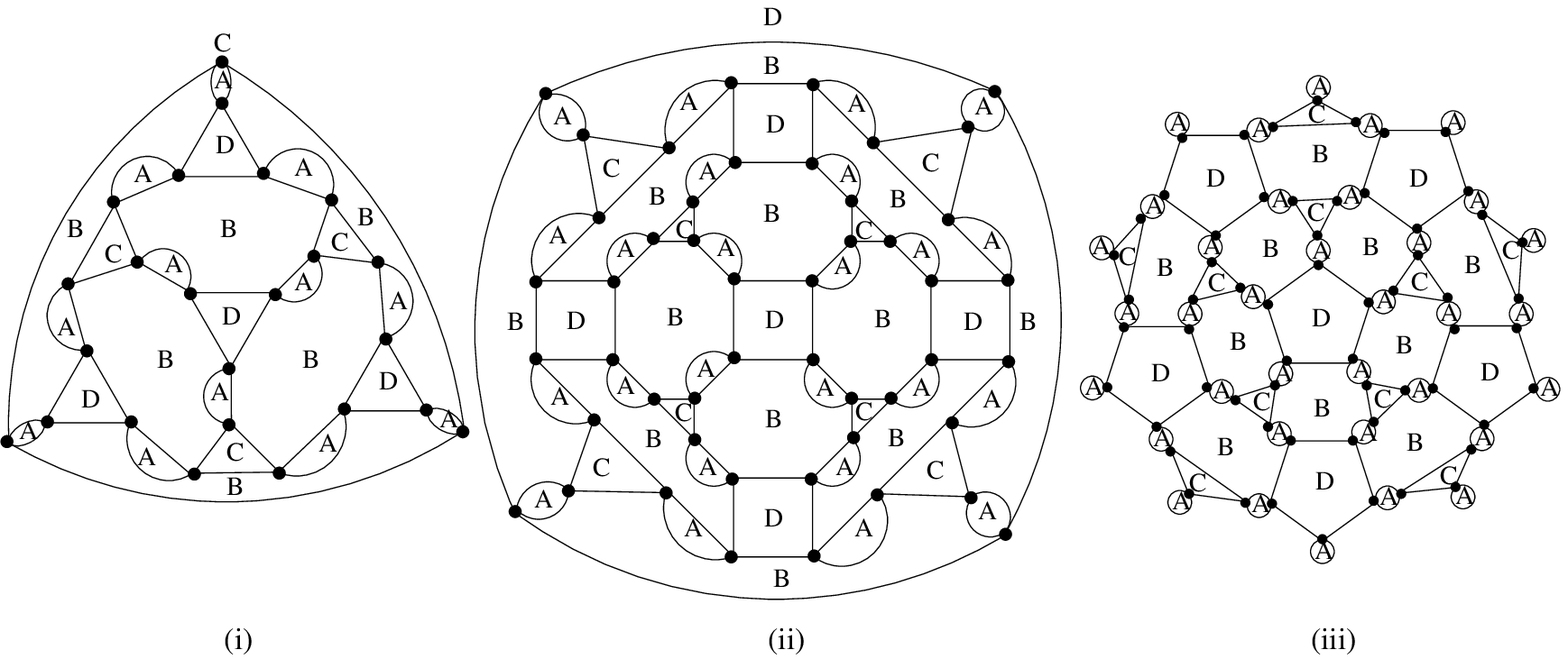}
\end{center}
\caption{spherical pictures}
\end{figure}

\begin{figure}
\begin{center}
\psfig{file=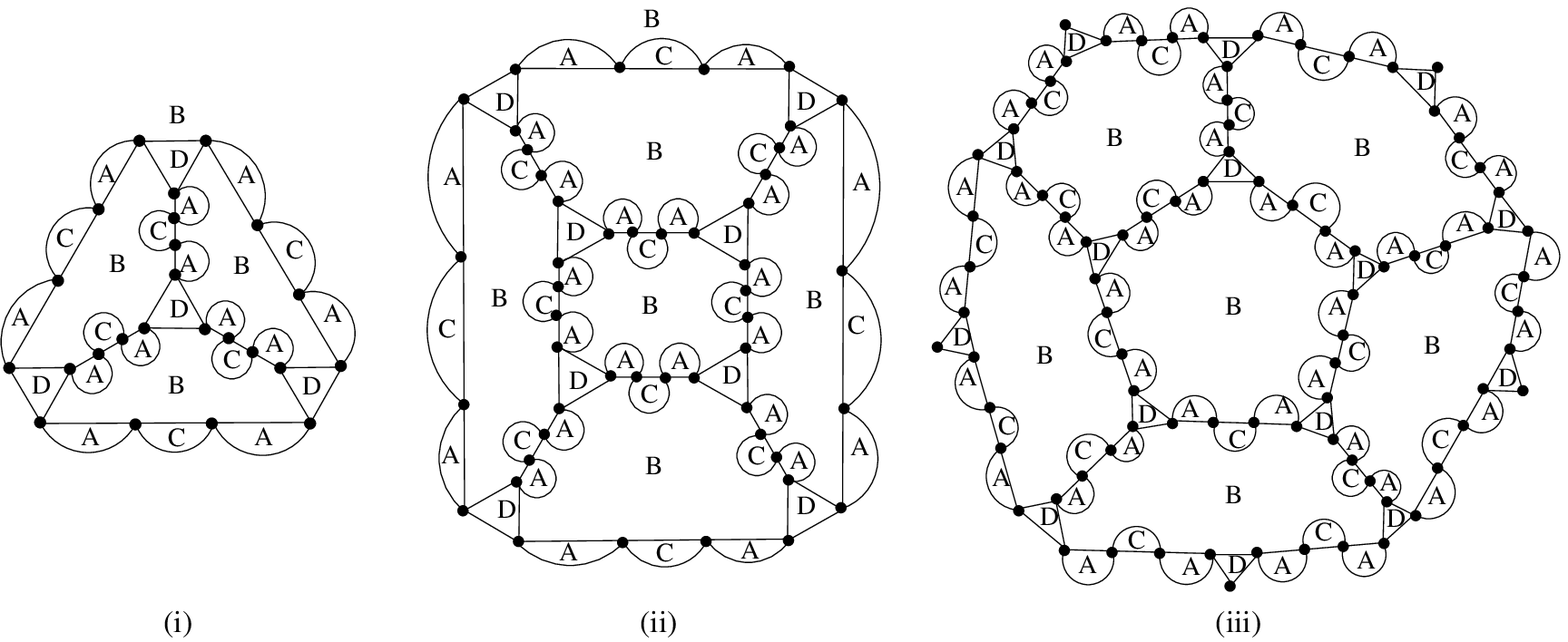}
\end{center}
\caption{spherical pictures}
\end{figure}

\begin{figure}
\begin{center}
\psfig{file=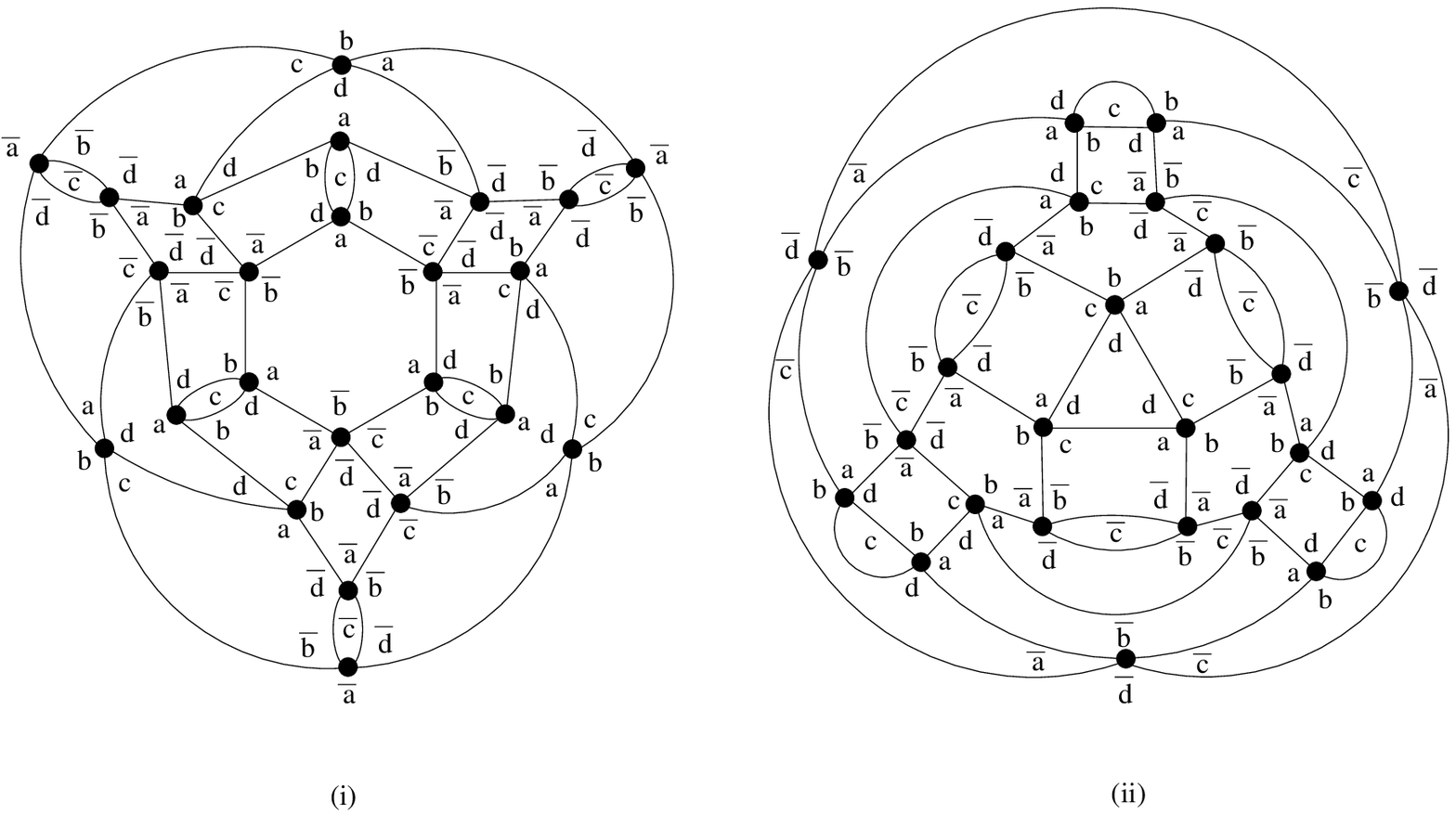}
\end{center}  
\caption{spherical pictures}
\end{figure}

\newpage
(iv) If $l=4$ a spherical picture is given by Figure 3.3(i) in which $A=(ab^{-1})^3$, $B=bda^{-1} c^{-1}$, $C=c^2$ and $D=d^{-4}$.  If $l=5$ there is a spherical picture consisting of 20 $A$ regions, 60 $B$ regions, 30 $C$ regions 
and 12 $D=d^{-5}$ regions, half of which is shown in Figure 3.3(ii).

(v) Spherical pictures for $k=3$ and 4 are given in Figure 3.4(i) and (ii) in which $A=(ab^{-1})^3$ and $(ab^{-1})^4$ (respectively), $B=bda^{-1} c^{-1}$, $C=c^2$ and
$D=d^{-3}$.  If $k=5$ there is a spherical picture consisting of 12 $A=(ab^{-1})^5$ regions, 60 $B$ regions, 30 $C$ regions and 20 $D$ regions half of which is shown in Figure 3.4(iii).


(c)(i) If $|bda^{-1} c^{-1} | = n < \infty$ then there is a spherical picture over $\mathcal{P}$ consisting of two 4$n$-gons each with label
$(bda^{-1} c^{-1})^n$ together with $n$ 2-gons with label $c^2$, $n$ 2-gons with label $d^{-2}$ and $2n$ 2-gons with label $ab^{-1}$.  The sphere for $n=3$ is shown in Figure 3.5(i).

(ii)If $|d|=n < \infty$ then there is a sphere consisting of two $n$-gons with label $D=d^{-n}$ and between them a layer consisting of $n$ 8-gons with label $B=(bda^{-1} c^{-1})^2$, $2n$ 2-gons with label $A=ab^{-1}$ and $n$ 2-gons 
with label $C=c^2$.  The sphere for $n=4$ is shown in Figure 3.5(ii).


(iii) A spherical picture for $l=4$ is given by Figure 3.5(iii) in which $A=ab^{-1}$, $B=(bda^{-1} c^{-1})^3$, $C=c^2$ and $D=d^{-4}$.  If $l=5$ there is a sphere consisting of 60 $A$ regions, 20 $B$ regions, 30 $C$ regions and 12 
$D=d^{-5}$ regions half of which is given by Figure 3.5(iv).

(iv) Spherical pictures for $k=3$ and 4 are given by Figure 3.6(i) and (ii) in which $A=ab^{-1}$, $B=(bda^{-1} c^{-1})^2$, $C=c^3$ and $D=d^{-3}$ and $d^{-4}$ (respectively). If $k=5$ there is a spherical picture consisting of 60 
$A$ regions, 30 $B$ regions, 20 $C$ regions and 12 $D=d^{-5}$ regions half of which is given by Figure 3.6(iii).

(v) Spherical pictures for $k=3$ and 4 are given by Figure 3.7(i) and (ii) in which $A=ab^{-1}$; $B=(bda^{-1}c^{-1})^3$ and $(bda^{-1} c^{-1})^4$ (respectively),
$C=c^2$ and $D=d^{-3}$.  If $k=5$ there is a spherical picture consisting of 60 $A$ regions, 12 $B=(bda^{-1}c^{-1})^5$ regions, 30 $C$ regions and 20 $D$ regions half of which is given by Figure 3.7(iii). $\Box$

It follows from Theorem 1(2) in [1] that if $|t| < \infty$ in $G(\mathcal{P})$ then $\mathcal{P}$ fails to be aspherical.  We apply this fact in the proof of the next lemma.

\textbf{Lemma 3.2}\quad \textit{If any of the following conditions hold then $\mathcal{P}$ fails to be aspherical.}
\begin{enumerate}
\item[(i)]
$|a^{-1}b| < \infty$, $|c|=|d|=2$ \textit{and} $a^{-1}cad=bdb^{-1}c=1$.
\item[(ii)]
$|d| < \infty$ \textit{and} $c^{-1}ab^{-1}=db^{-1}a=1$.
\item[(iii)]
$|d| < \infty$ \textit{and} $a^{-1}b=cada^{-1}=1$.
\item[(iv)]
$c^2 = cbda^{-1} = a^{-1} bd^{-2} =1$.
\item[(v)]
$c^2 = cbda^{-1} = (a^{-1} b)^2 d^{-1} =1$.
\item[(vi)]
$a^{-1}b=c^2=d^3=1$ \textit{and} $cada^{-1} cad^{-1} a^{-1}=1$.
\item[(vii)]
$a^{-1}b=c^2ada^{-1}=1$ \textit{and} $4 \leq |c| \leq 5$.
\item[(viii)]
$a^{-1}b=c^3 ada^{-1}=1$ \textit{and} $|c|=6$.
\end{enumerate}
\textit{Proof}
\begin{enumerate}
\item[(i)]
It is enough to show that the group $G=\langle b,d,t \mid d^2=b^k=1, bd=db, t^2 btdt^{-1}d=1 \rangle$ has order $2k(3^{2k}-1)$.
Now $G= \langle d,t \mid d^2, t^{-2} d^{-1} td^{-1} t^{-1} dtdt^{-1} dt^2 d^{-1}, (t^3 dt^{-1} d)^k \rangle$ and
$G/G' = \langle d,t \mid d^2 = t^{2k} = [d,t] = 1 \rangle$.  Let $\mathcal{K}$ denote the covering 2-complex associated with $G'$ [3].  Then $\mathcal{K}$ has edges $t_{0j}, t_{1j}, d_{j0}, d_{j1}$ ($1 \leq j \leq 2k$) and 2-cells $d_{j0} d_{j1}$, $t_{0j} d_{1-j} t_{1j}^{-1} d_{2-j}$,
$t_{i1} t_{i2} \ldots t_{i2k}$ where $1 \leq i \leq 2$ and $1 \leq j \leq 2k$ and the $d$ subscripts are mod $2k$.  Collapsing the maximal subtree whose edges are $d_{j0}$ ($1 \leq j \leq 2k$), $t_{0l}$ ($2 \leq l \leq 2k$) and using the lifts of $d^2$ shows that $G'=\langle t_{01},t_{1j} ~(1 \leq j \leq 2k) \rangle$.  Using the lifts of the second relator it is easily shown that $G'= \langle t_{01}, t_{11}, t_{12} \rangle$ where
$t_{11} t_{01}^{-1} = t_{12}^{-3^{2k-1}}$ and, finally, using the lift of the third relator $(t^3 dt^{-1} d)^k$ one can show that
$G'=\langle t_{12} \mid t_{12}^r \rangle$ where $r=\frac{1}{2} (3^{2k}-1)$.  We omit the details.
\item[(ii)]
It is enough to show that $G=\langle d,x \mid d^k, t^2 dtd^{-1} t^{-1} d \rangle$ has order $2k(1+ 4 + 4^2 + \ldots + 4^{k-1} )$.
Now $G= \langle u,t \mid (ut^{-2})^k, tut^{-1} u^{-2} \rangle$ and $G/G' = \langle u,t \mid u=t^{2k}=1 \rangle$.  Let $\mathcal{L}$ denote the covering complex associated with $G'$.  Then $\mathcal{L}$ has edges $t_j,u_j$ ($1 \leq j \leq 2k$) and 2-cells $t_1 t_2 \ldots t_{2k},u_j$ ($1 \leq j \leq 2k$).
Collapsing the maximal tree whose edges are $t_1, \ldots, t_{k-1}$ implies $G'= \langle t_{2k}, u_j ~(1 \leq j \leq 2k) \rangle$.  The lifts of
$tu t^{-1} u^{-2}$ yield the relators $u_l = u_1^{2^{l-1}}$ for $2 \leq l \leq 2k$ and $t_{2k} u_1 t_{2k}^{-1} u_1^{-4^k}$.
The lifts of $(ut^{-2})^k$ yield the relators
$t_{2k}^{-1} = \prod_{i=0}^{k-1} u_{2i+1} = \prod_{i=1}^k u_{2i}$.  It follows that $G'=\langle u_1 \mid u_1^r \rangle$ where
$r=1+4 + 4^2 + \ldots + 4^{k-1}$.
\item[(iii)]
Here $r=t^3 dt^{-1} d^{-1}$ and $|d| < \infty$ implies $|t| < \infty$.
\item[(iv)]
-- (v) A spherical picture for (iv), (v) is shown in Figure 3.8(i), (ii) (respectively).
\item[(vi)]
-- (viii) For these cases we use GAP [9].  For (vi), $r=t^3 ct^{-1} d$ together with the conditions yields $|t| \leq 12$; for (vii) $r=t^3 ct^{-1} c^{-2}$ and $|c| = 4,5$ implies $|t| \leq 8,10$ (respectively); and for (viii) 
$r=t^3 ct^{-1} c^{-3}$ and $|c|=6$ implies $|t| \leq 24$. $\Box$
\end{enumerate}

\textbf{Lemma 3.3}\quad \textit{If any of the conditions (i)--(iii) of Theorem 1.1 or (i)--(x) of Theorem 1.2 holds then $\mathcal{P}$ fails to be aspherical.}

\textit{Proof} Consider Theorem 1.1. If (i) holds then $\mathcal{P}$ fails to be aspherical by Lemma 3.1(a)(i).
If (ii) holds then $\mathcal{P}$ fails to be aspherical by Lemma 3.2(i).
This leaves condition (iii).  If $a^{-1} b \neq 1$ and $bda^{-1} c^{-1} \neq 1$ then (iii) does not hold; and if $a^{-1}b=bda^{-1} c^{-1}=1$ then $H$ is cyclic.
Let $a^{-1}b=1$.  Since $(|c|,|d|,|bda^{-1} c^{-1}|)$ is $T$-equivalent to $(|d|,|c|,|bda^{-1} c^{-1}|)$ it can be assumed without any loss that
$|c| \leq |d|$.  The resulting ten cases are dealt with by Lemma 3.1(c).  Let $bda^{-1} c^{-1}=1$.  Since $(|a^{-1}b|,|c|,|d|)$ is $T$-equivalent to
$(|a^{-1}b|,|d|,|c|)$ it can again be assumed without any loss that $|c| \leq |d|$.  The resulting ten cases are dealt with by Lemma 3.1(b).
Now consider Theorem 1.2.
If (i) holds then $\mathcal{P}$ is aspherical by Lemma 3.1(a)(i);
if (ii) holds then by Lemma 3.1(a)(ii);
if (iii) holds then by Lemma 3.2(ii);
if (iv) holds then by Lemma 3.2(i);
if (v) holds then by Lemma 3.2(iv);
if (vi) holds then by Lemma 3.2(v);
if (vii) holds then by Lemmas 3.1(a)(iii) and 3.2(iii);
if (viii) holds then by Lemma 3.2(vi);
if (ix) holds then by Lemma 3.2(vii); and
if (x) holds then by Lemma 3.2(viii). $\Box$

A \textit{weight function} $\alpha$ on the star graph $\Gamma$ of Figure 2.1(iii) is a real-valued function on the set of edges of $\Gamma$.
Denote the edge labelled $a,b,c,d$ by $e_a,e_b,e_c,e_d$ (respectively).  The function $\alpha$ is \textit{weakly aspherical} if the following two conditions are satisfied:
\begin{enumerate}
\item[(1)]
$\alpha (e_a) + \alpha (e_b) + \alpha (e_c) + \alpha (e_d) \leq 2$;
\item[(2)]
each admissible cycle in $\Gamma$ has weight at least 2.
\end{enumerate}
If there is a weakly aspherical weight function on $\Gamma$ then $\mathcal{P}$ is aspherical [2].

\medskip

\textbf{Lemma 3.4}\quad \textit{If any of the following conditions holds then $\mathcal{P}$ is aspherical.}
\begin{enumerate}
\item[(i)]
$|c|=|d|=\infty$\textit{;}
\item[(ii)]
$1 < |b| < \infty$ \textit{and} $|d| = \infty$\textit{;}
\item[(iii)]
\textit{$|c| < \infty$, $|d| < \infty$ and $|b| = \infty$.}
\end{enumerate}

\textit{Proof} The following functions $\alpha$ are weakly aspherical.
\begin{enumerate}
\item[(i)]
$\alpha (e_a)= \alpha (e_b) =1$, $\alpha (e_c) = \alpha (e_d) = 0$.
\item[(ii)]
$\alpha (e_a) = \alpha (e_b) = \frac{1}{2}$, $\alpha (e_c)=1$, $\alpha (e_d)=0$.
\item[(iii)]
$\alpha (e_a)=\alpha (e_b)=0$, $\alpha (e_c)=\alpha (e_d)=1$. $\Box$
\end{enumerate}

The following lemmas will be useful in later sections.

\textbf{Lemma 3.5}\quad \textit{Let $d(\hat{\Delta})=k$ where $\hat{\Delta}$ is a region of the spherical picture $\bs{P}$ over $\mathcal{P}$.}
\begin{enumerate}
\item[(i)]
\textit{If $\hat{\Delta}$ receives at most $\frac{\pi}{6}$ across each edge and $k \geq 6$ then $c^{\ast}(\hat{\Delta}) \leq 0$.}
\item[(ii)]
\textit{If $\hat{\Delta}$ receives at most $\frac{\pi}{6}$ across at most two-thirds of its edges, nothing across the remaining edges and $k \geq 5$ then
$c^{\ast}(\hat{\Delta}) \leq 0$.}
\item[(iii)]
\textit{If $\hat{\Delta}$ receives at most $\frac{\pi}{4}$ across each edge and $k \geq 8$ then $c^{\ast}(\hat{\Delta}) \leq 0$.}
\item[(iv)]
\textit{If $\hat{\Delta}$ receives at most $\frac{\pi}{4}$ across at most two-thirds of its edges, nothing across the remaining edges and $k \geq 6$ then $c^{\ast}(\hat{\Delta}) \leq 0$.}
\item[(v)]
\textit{If $\hat{\Delta}$ receives at most $\frac{\pi}{4}$ across at most half of its edges, nothing across the remaining edges and $k \geq 5$ then
$c^{\ast}(\hat{\Delta}) \leq 0$.}
\item[(vi)]
\textit{If $\hat{\Delta}$ receives at most $\frac{\pi}{2}$ across at most half of its edges, nothing across the remaining edges and $k \geq 7$ then
$c^{\ast}(\hat{\Delta}) \leq 0$.}
\item[(vii)]
\textit{If $\hat{\Delta}$ receives at most $\frac{\pi}{2}$ across at most three-fifths of its edges, nothing across the remaining edges and $k \geq 8$ then $c^{\ast} (\hat{\Delta}) \leq 0$.}
\end{enumerate}
\textit{Proof} Recall that $c(\hat{\Delta})=\pi \left( 2 - \frac{k}{2} \right)$.
(i) $c^{\ast}(\hat{\Delta}) \leq c(\hat{\Delta}) + k. \frac{\pi}{6} \leq 0$ for $k \geq 6$.
(ii) $c^{\ast}(\hat{\Delta}) \leq c(\hat{\Delta}) + \frac{2k}{3} . \frac{\pi}{6} \leq 0$ for $k \geq 6$.
If $k=5$ then $\hat{\Delta}$ receives at most $\frac{\pi}{6}$ across at most three edges and so $c^{\ast}(\hat{\Delta}) \leq \pi \left( 2 - \frac{5}{2} \right) + 3. \frac{\pi}{6}=0$.
(iii) $c^{\ast}(\hat{\Delta}) \leq c(\Delta) + k. \frac{\pi}{4} \leq 0$ for $k \geq 8$.
(iv) $c^{\ast}(\hat{\Delta}) \leq c(\Delta) + \frac{2k}{3} . \frac{\pi}{4} \leq 0$ for $k \geq 6$.
(v) $c^{\ast}(\hat{\Delta}) \leq c(\Delta) + \frac{k}{2} . \frac{\pi}{4} \leq 0$ for $k \geq 6$.
If $k=5$ then $\hat{\Delta}$ receives at most
$\frac{\pi}{4}$ across at most two edges and $c^{\ast}(\hat{\Delta}) \leq \pi \left( 2 - \frac{5}{2} \right) + 2. \frac{\pi}{4} = 0$.
(vi) $c^{\ast} (\hat{\Delta}) \leq c(\Delta) + \frac{k}{2} . \frac{\pi}{2} \leq 0$ for $k \geq 8$.
If $k=7$ then $\hat{\Delta}$ receives at most
$\frac{\pi}{2}$ across at most three edges and $c^{\ast}(\hat{\Delta}) \leq \pi \left( 2 - \frac{7}{2} \right) + 3 . \frac{\pi}{2} = 0$.
(vii) $c^{\ast}(\hat{\Delta}) \leq c(\Delta) + \frac{3k}{5} . \frac{\pi}{2} \leq 0$ for $k \geq 10$.  If $k=9$ then $\hat{\Delta}$ receives at most $\frac{\pi}{2}$ across five edges and $c^{\ast}(\hat{\Delta}) \leq \pi ( 2 - \frac{9}{2}) + 5 . \frac{\pi}{2} = 0$; if $k=8$ then $\hat{\Delta}$ receives across at most four edges and $c^{\ast}(\hat{\Delta}) \leq \pi (2 - \frac{8}{2}) + 4. \frac{\pi}{2} = 0$. $\Box$

\textbf{Remark}\quad We will use the above lemmas as follows.  Suppose that $\hat{\Delta}$ receives positive curvature across its edge $e_i$.
If we know that it then never receives curvatures across $e_{i-1}$ or across $e_{i+1}$ then we can apply the ``half'' results; or if we know that it receives positive curvature across at most one of $e_{i-1},e_{i+1}$ then we can apply the ``two-thirds'' results.

Let $\hat{\Delta}$ be a region of $\bs{P}$ and let $e$ be an edge of $\hat{\Delta}$.  If $\hat{\Delta}$ receives no curvature across $e$ then $e$ is called a \textit{gap}; if at most $\frac{\pi}{6}$ then $e$ is called a \textit{two-thirds gap}; and if at most $\frac{\pi}{4}$ then $e$ is called a \textit{half gap}.

\textbf{Lemma 3.6 (The Four Gaps Lemma)}\quad
\textit{If $\hat{\Delta}$ has a total of at least four gaps (in particular, four edges across which $\hat{\Delta}$ does not receive any curvature) and the most curvature that crosses any edge is $\frac{\pi}{2}$ then $c^{*}(\hat{\Delta})\le 0$}.

\textit{Proof}  It follows that $c^{*}(\hat{\Delta}) \leq \pi \left( 2 - \frac{k}{2} \right) +(k-4)\frac{\pi}{2}=0$. $\Box$

Checking the star graph shows that we will have the following LIST for the labels of regions of degree $k$ where $k\in\{2,3,4,5,6,7\}$:

If $d(\Delta)=2$ then $l(\Delta)\in\{c^{2},a^{-1}b,d^{2}\}$.

If $d(\Delta)=3$ then $l(\Delta)\in\{c^{3},cab^{-1},c^{-1}ab^{-1},db^{-1}a,d^{-1}b^{-1}a,d^{3}\}$.

If $d(\Delta)=4$ then $l(\Delta)\in\{ d^{4},d^{2}a^{-1}b,d^{2}b^{-1}a,c^{2}ab^{-1},c^{2}ba^{-1},c^{4},ab^{-1}ab^{-1},
d\{a^{-1},b^{-1}\}\\\{c,c^{-1}\}\{a,b\}\}$.

If $d(\Delta)=5$ then $l(\Delta)\in\{d^{5},d^{3}a^{-1}b,d^{3}b^{-1}a,c^{3}ab^{-1},c^{3}ba^{-1},c^{5},
cab^{-1}ab^{-1},cba^{-1}ba^{-1},\\da^{-1}ba^{-1}b,db^{-1}ab^{-1}a,d^{2}\{a^{-1},b^{-1}\}\{c,c^{-1}\}\{a,b\},
c^{2}\{a,b\}\{d,d^{-1}\}\{a^{-1},b^{-1}\}\}$.

If $d(\Delta)=6$ then $l(\Delta)\in\{d^{6},d^{4}a^{-1}b,d^{4}b^{-1}a,c^{4}ab^{-1},c^{4}ba^{-1},c^{6},ab^{-1}ab^{-1}ab^{-1},
d^{2}a^{-1}\\ba^{-1}b,d^{2}b^{-1}ab^{-1}a,c^{2}ab^{-1}ab^{-1},c^{2}ba^{-1}ba^{-1},
d^{3}\{a^{-1},b^{-1}\}\{c,c^{-1}\}\{a,b\},d^{2}\{a^{-1},b^{-1}\}\\\{c^{2},c^{-2}\}\{a,b\},
d\{a^{-1},b^{-1}\}\{c^{3},c^{-3}\}\{a,b\},c\{ab^{-1},ba^{-1}\}\{c,c^{-1}\}\{ab^{-1},ba^{-1}\},
d\{a^{-1}\\b,b^{-1}a\}\{d,d^{-1}\}\{a^{-1}b,b^{-1}a\},
c\{ab^{-1}a, ba^{-1}b \} \{ d,d^{-1} \} \{ a^{-1}, b^{-1} \}, c\{a,b\} \{d,d^{-1} \} \{ a^{-1}ba^{-1}, b^{-1} ab^{-1} \}
\}$.

Where, for example $d^2 \{ a^{-1},b^{-1}\} \{c,c^{-1} \} \{ a,b \}$ yields the eight labels $d^2 a^{-1} c^{\pm 1} a$, $d^2 a^{-1} c^{\pm 1} b$,
$d^2 b^{-1} c^{\pm 1} a$, $d^2 b^{-1} c^{\pm 1} b$.

We will use the above LIST throughout the following sections often without explicit reference.

\section{Proof of Case A} 

In this section we prove Theorems 1.1 and 1.2 for Case A, that is, we make the following assumptions:
\[
|c|>2,~ |d|>2 \text{ and } a^{-1}b\ne 1.
\]
This implies that $d(\Delta)\ge3$ for each region $\Delta$ of the spherical diagram $\bs{P}$.
If $d(\Delta)=3$ then we will fix the names of the fifteen neighbouring regions $\Delta_{i}$ ($1 \le i \le 15$) of $\Delta$ as shown in Figure 4.1(i) and we use this notation throughout the section.
We treat each of the cases (\textbf{A0})--(\textbf{A10}) in turn.

\textbf{(A0)} $\boldsymbol{|c|>3}$, $\boldsymbol{|d|>3}$, $\bs{a^{-1}b\ne1}$, $\bs{c^{\pm1}ab^{-1}\ne1}$, $\bs{d^{\pm1}b^{-1}a\ne1}$.\newline
In this case $d(\Delta)>3$ for all regions $\Delta$. Since the degree of each vertex equals $4$ it follows that $c(\Delta)=(2-d(\Delta))\pi+d(\Delta)\frac{2\pi}{4}\le 0$ and so $\mathcal{P}$ is aspherical.

\textbf{(A1)} $\boldsymbol{|c|=3}$, $\bs{|d|>3}$, $\bs{a^{-1}b\ne1}$, $\bs{c^{\pm1}ab^{-1}\ne1}$, $\bs{d^{\pm1}b^{-1}a\ne1}$.\newline
If $d(\Delta)=3$ then $c(\Delta)=c(4,4,4)=\frac{\pi}{2}$, $l(\Delta)=c^{3}$ and $\Delta$ is given by Figure 4.1(ii).
Observe that if $d(\Delta_{i})=4$ for $i\in\{1,3,5\}$ then $l(\Delta_{i})=bd\omega$ and so $l(\Delta_{i})\in\{bd^{2}a^{-1},bda^{-1}c^{\pm1},bdb^{-1}c^{\pm1}\}$. (See the LIST of Section 3.) But $bdb^{-1}c^{\pm1}$ implies $|d|=|c|$, a contradiction. Observe further that at most one of $bd^{2}a^{-1}$, $bda^{-1}c$, $bda^{-1}c^{-1}$ equals $1$ otherwise there is a contradiction to $|c|=3$ or $|d|>3$. This leaves the following cases:
(i) $bd^{2}a^{-1}\ne1$; $bda^{-1}c^{\pm1}\ne1$;
(ii) $bd^{2}a^{-1}=1$;
(iii) $bda^{-1}c=1$;
(iv) $bda^{-1}c^{-1}=1$.

Consider case (i). Here $d(\Delta_{i})>4$ for $i\in\{1,3,5\}$ of Figure 4.1(ii), so add $\frac{1}{3}c(\Delta)=\frac{\pi}{6}$ to each of $c(\Delta_{i})$ across the $bd$ edge as shown. Observe from Figure 4.1(ii) that $\Delta_{1}$, say, does 

\newpage
\begin{figure}
\begin{center}
\psfig{file=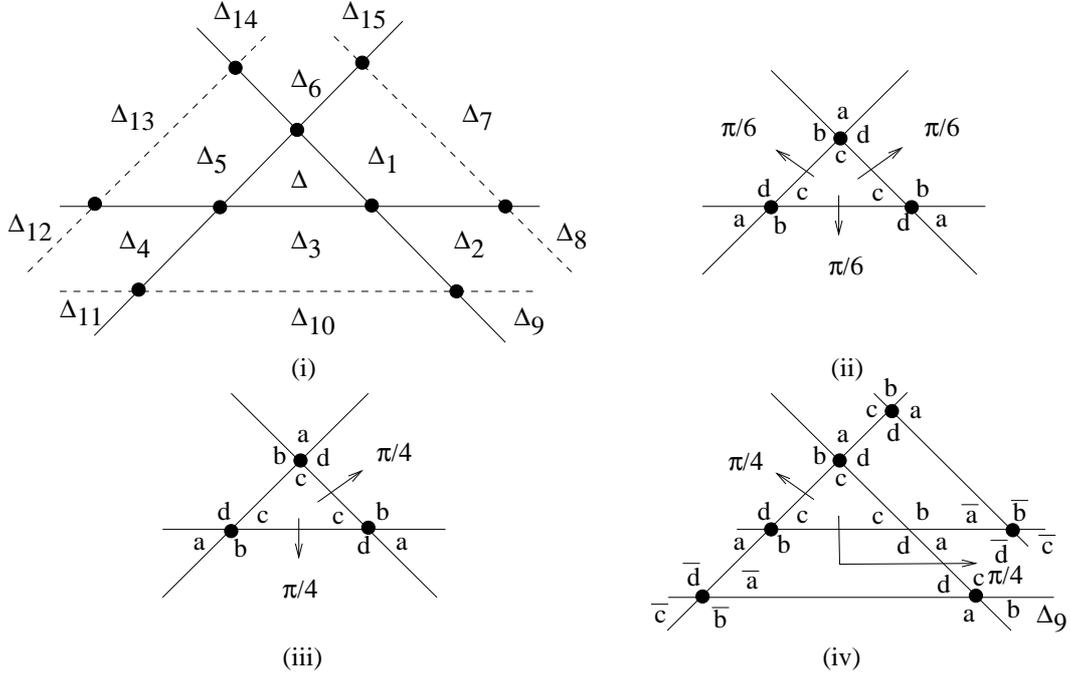}
\end{center}
\caption{the region $\Delta$ and curvature distribution for Cases (A1,2,6,7,8) (B9,12)}
\end{figure}

\noindent not receive any curvature from $\Delta_{2}$ or $\Delta_{6}$. Therefore if $\hat{\Delta}$ receives positive curvature, since $d(\hat{\Delta})\ge5$, it follows that $c^{*}(\hat{\Delta})\le 0$ by Lemma 3.5(ii).

Now consider case (ii), $bd^2 a^{-1}=1$. If $d(\Delta_{i})>4$ for at least two of the $\Delta_{i}$ where $i\in\{1,3,5\}$, say $\Delta_{1}$ and $\Delta_{3}$, then add $\frac{1}{2} c(\Delta)=\frac{\pi}{4}$ to each of $c(\Delta_{1})$ and $c(\Delta_{3})$ across the $bd$ edge as in Figure 4.1(iii).
Otherwise by symmetry it can be assumed without any loss of generality that $d(\Delta_{1})=d(\Delta_{3})=4$. Then add $\frac{1}{2} c(\Delta)=\frac{\pi}{4}$ firstly to $c(\Delta_{3})$ and then on to $c(\Delta_{2})$ across the $bd$ and $ca$ edges as shown in Figure 4.1(iv). If $d(\Delta_{5})>4$ then add the remaining $\frac{1}{2} c(\Delta)=\frac{\pi}{4}$ to $c(\Delta_{5})$ across the $bd$ edge; otherwise $d(\Delta_{5})=4$ and similarly add the $\frac{\pi}{4}$ firstly to $c(\Delta_{5})$ then on to $c(\Delta_{4})$. Now observe that in Figure 4.1(iii) $\Delta_{1}$ does 
not receive positive curvature from $\Delta_{2}$ or $\Delta_{6}$ and in Figure 4.1(iv) that $\Delta_{2}$ does not receive positive curvature from $\Delta_{1}$ or $\Delta_{9}$.
It follows that if $\hat{\Delta}$ receives positive curvature then $\hat{\Delta}$ receives curvature across at most half of its edges and so if $d(\hat{\Delta})\ge5$, it follows that $c^{*}(\hat{\Delta})\le 0$ by Lemma 3.5(v).
This leaves the case when $d(\hat{\Delta})=4$ and $l(\hat{\Delta})=cad^{-1}\omega$. Therefore $l(\hat{\Delta})\in\{cad^{-1}b^{-1},cad^{-1}a^{-1}\}$. In each case $l(\hat{\Delta})$ together with $bd^{2}a^{-1}=c^{3}=1$ implies $|d|=3$, a contradiction.

Now consider (iii), $bda^{-1}c=1$. As before, if $d(\Delta_{1})>4$ and $d(\Delta_{3})>4$ then add $\frac{1}{2} c(\Delta)=\frac{\pi}{4}$ to each of $c(\Delta_{1})$ and $c(\Delta_{3})$ across the $bd$ edge as in Figure 4.1(iii). If $d(\Delta_{1})= d(\Delta_{3})=4$ then add $\frac{1}{2} c(\Delta)=\frac{\pi}{4}$ firstly to $c(\Delta_{1})$ and then on to $c(\Delta_{2})$ across the $bd$ and $ad$ edges as shown in Figure 4.2(i). The remaining $\frac{1}{2} c(\Delta)=\frac{\pi}{4}$ is distributed either to $\Delta_5$ when $d(\Delta_5) > 4$ or to $\Delta_6$ via $\Delta_5$ when $d(\Delta_5)=4$. Now observe that in Figure 4.1(iii) $\Delta_{1}$ does 

\newpage
\begin{figure}
\begin{center}
\psfig{file=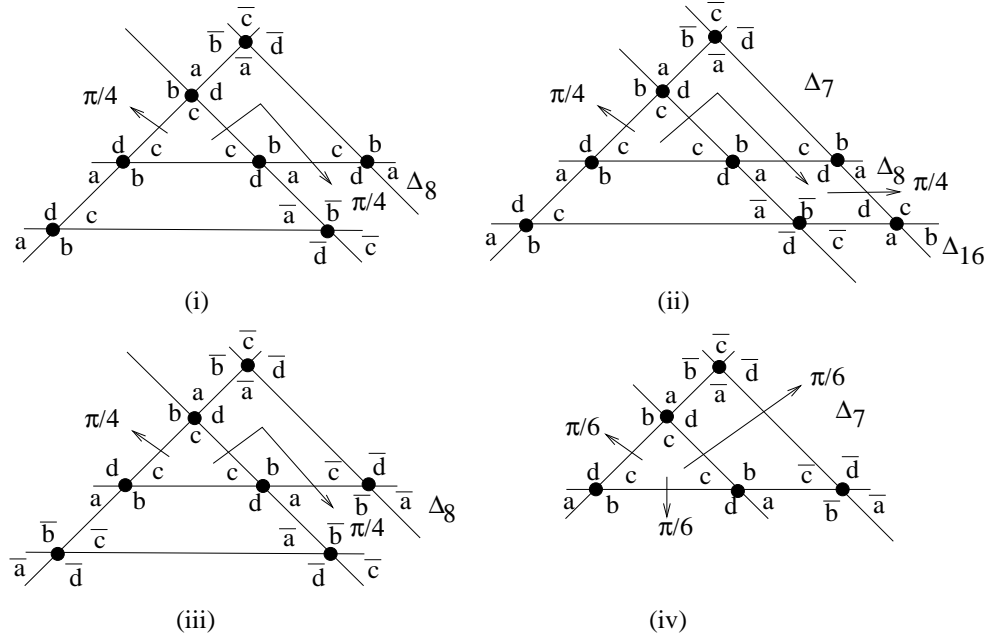}
\end{center}  
\caption{curvature distribution for Cases (A1) and (A7)}
\end{figure}

\noindent not receive positive curvature from $\Delta_{2}$ or $\Delta_{6}$; and in Figure 4.2(i) that $\Delta_{2}$ does not receive positive curvature from $\Delta_{3}$ or $\Delta_{8}$.
It follows that if the region $\hat{\Delta}$ receives positive curvature then it does so across at most half of its edges and therefore if $d(\hat{\Delta})\ge5$ then $c^{*}(\hat{\Delta})\le 0$ by Lemma 3.5(v).
This leaves the case when $ d(\hat{\Delta})=4$ and $l(\hat{\Delta})=b^{-1}ad\omega$. Therefore $l(\hat{\Delta})=b^{-1}ad^{2}$ and $H$ is cyclic.
If this occurs then add $\frac{1}{2} c(\Delta)=\frac{\pi}{4}$ to $c(\Delta_8)$ as shown in Figure 4.2(ii).
If $d(\Delta_8)=4$ then $l(\Delta_8) \in \{ cad^{\pm 1} b^{-1}, cad^{\pm 1} a^{-1}, cab^{-1} c \}$ which, together with $b^{-1} add$, contradicts the (\textbf{A1}) assumptions.
Also note that in Figure 4.2(ii), $\Delta_8$ does not receive positive curvature from $\Delta_7$ or from $\Delta_{16}$ and so if
$\hat{\Delta}$ receives positive curvature it does so across at most half its edges and $c^{\ast}(\hat{\Delta}) \leq 0$.

Finally consider case (iv), $bda^{-1}c^{-1}=1$. First assume that $b^{2}\ne1$. Again if $d(\Delta_{1})>4$ and $d(\Delta_{3})>4$ then add $\frac{1}{2} c(\Delta)=\frac{\pi}{4}$ to each of $c(\Delta_{1})$ and $c(\Delta_{3})$ across the $bd$ edge as in Figure 4.1(iii). If $d(\Delta_{1})=d(\Delta_{3})=4$ then add $\frac{1}{2} c(\Delta)=\frac{\pi}{4}$ firstly to $c(\Delta_{1})$ and then on to $c(\Delta_{2})$ across the $bd$ and $ab^{-1}$ edges as shown in Figure 4.2(iii) and add the remaining $\frac{1}{2}c(\Delta)=\frac{\pi}{4}$ to $\Delta_5$ or $\Delta_6$ as in the above. Now observe that in Figure 4.1(iii) $\Delta_{1}$ does not receive positive curvature from $\Delta_{2}$ or $\Delta_{6}$ and in Figure 4.2(iii) that $\Delta_{2}$ does not receive positive curvature from $\Delta_{3}$ or $\Delta_{8}$.
Thus if $\hat{\Delta}$ receives positive curvature then it does so across at most half of its edges so $d(\hat{\Delta})\ge5$ implies $c^{*}(\hat{\Delta})\le0$ by Lemma 3.5(v).
This leaves the case when $ d(\hat{\Delta})=4$ and $l(\hat{\Delta})=b^{-1}ab^{-1}\omega$. Therefore $l(\hat{\Delta})=b^{-1}ab^{-1}a$ and there is a contradiction to $b^{2}\ne1$.
Now assume that $b^{2}=1$. The distribution in this case is different in that we add $\frac{1}{3} c(\Delta)=\frac{\pi}{6}$ to $c(\Delta_{j})$ for $j\in\{1,3,5\}$ as in Figure 4.1(ii) across the $bd$ edge. If $d(\Delta_{1})=4$, say, then distribute the $\frac{\pi}{6}$ further on to $c(\Delta_{7})$ across the $d^{-2}$ edge 

\newpage
\begin{figure}
\begin{center}
\psfig{file=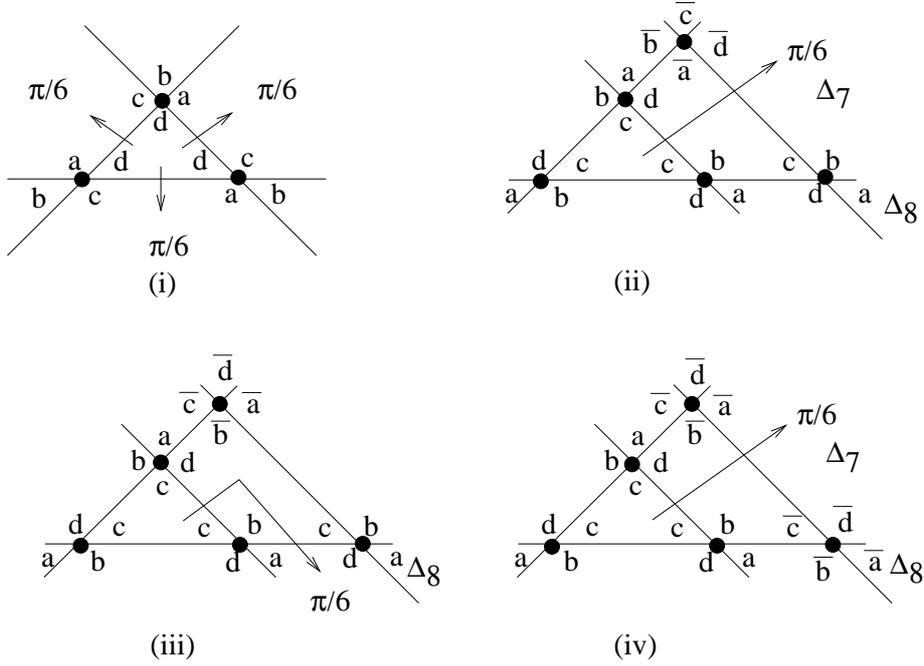}
\end{center}
\caption{curvature distribution for Cases (A2,8) and (B3,6,11,12)}
\end{figure}

\noindent as shown in Figure 4.2(iv).
Therefore if $\hat{\Delta}$ receives positive curvature and if $d(\hat{\Delta})\ge6$, it follows that $c^{*}(\hat{\Delta})\le 0$ by Lemma 3.5(i).
This leaves the case when $4 \le d(\hat{\Delta})\le 5$ and $l(\hat{\Delta})\in\{bd\omega,d^{-2}\omega\}$. If $d(\hat{\Delta})=4$ then $l(\hat{\Delta})=d^{-2}\omega$ and
$l(\hat{\Delta})\in\{d^{-4},d^{-2}a^{-1}b,d^{-2}b^{-1}a\}$ therefore
$d^{4}=1$ and there is a sphere by Lemma 3.1(b)(iii). So let $d(\hat{\Delta})=5$. If $\hat{\Delta}$ receives curvature across at most three edges then $c^{*}(\hat{\Delta})\le \pi ( 2 - \frac{5}{2}) + 3.\frac{\pi}{6}=0 $. Checking (the LIST) now shows that $c^{*}(\hat{\Delta})\le0$ except when $l(\hat{\Delta})=d^{-5}$ and in this case we obtain a sphere by Lemma 3.1(b)(iii).

In conclusion $\mathcal{P}$ fails to be aspherical in this case if and only if $b^2=bda^{-1} c^{-1} = 1$ and $|d| \in \{ 4,5 \}$.  Note that if these conditions hold then $H$ is non-cyclic for otherwise we would obtain $d^6=1$, a contradiction.

\textbf{(A2)} $\boldsymbol{|c|=3}$, $\bs{|d|=3}$, $\bs{a^{-1}b\ne1}$, $\bs{c^{\pm1}ab^{-1}\ne1}$, $\bs{d^{\pm1}b^{-1}a\ne1}$.\newline
If $d(\Delta)=3$ then $\Delta$ is given by Figures 4.1(ii) and 4.3(i). If $d(\Delta)=4$ and $l(\Delta)\in\{bd\omega,ca\omega\}$ then $l(\Delta)\in S = \{bda^{-1}c^{\pm1},bdb^{-1}c^{\pm 1},cad^{\pm1}a^{-1},cad^{\pm1}b^{-1}\}$ otherwise there is a contradiction to one of the assumptions.
(Throughout this case unless otherwise stated this means one of the (\textbf{A2}) assumptions.)

The cases to be considered are (where for example case (ii) means $bda^{-1}c$ is the only member of $S$ to equal 1):
(i) $bda^{-1}c^{\pm1}\ne1$, $bdb^{-1}c^{\pm1}\ne1$, $cad^{\pm1}a^{-1} \ne1$, $cadb^{-1}\ne1$;
(ii) $bda^{-1}c=1$;
(iii) $bda^{-1}c^{-1}=1$;
(iv) $bdb^{-1}c=1$;
(v) $bdb^{-1}c^{-1}=1$;
(vi) $cada^{-1}=1$;
(vii) $cad^{-1}a^{-1}=1$;
(viii) $cadb^{-1}=1$;
(ix) $bdb^{-1}c=1$, $cada^{-1}=1$;
(x) $bdb^{-1}c=1$, $cad^{-1}a^{-1}=1$;
(xi) $bdb^{-1}c^{-1}=1$, $cada^{-1}=1$;
(xii) $bdb^{-1}c^{-1}=1$, $cad^{-1}a^{-1}=1$.

\newpage
\begin{figure}
\begin{center}
\psfig{file=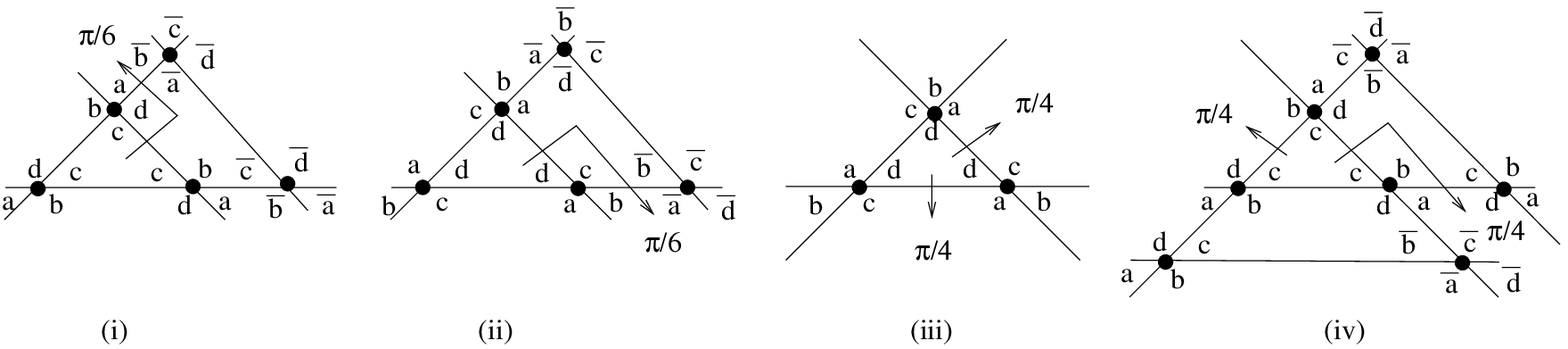}
\end{center}
\caption{curvature distribution for Case (A2)}
\end{figure}  

Note that any other combination of these conditions gives a contradiction to one of the assumptions. Moreover, (ii) is T-equivalent to (viii); (iv) is T-equivalent to (vi); (v) is T-equivalent to (vii); and (x) is T-equivalent to (xi). So it remains to consider (i), (ii), (iii), (iv), (v), (ix), (x) and (xii).

Now let $c(\Delta)>0$ and so $l(\Delta)\in\{c^{3},d^{3}\}$.
In cases (i), (ii), (iv) and (v) add $\frac{1}{3} c(\Delta)=\frac{\pi}{6}$ to $c(\Delta_{i})$ for $i\in\{1,3,5\}$ as shown in Figures 4.1(ii) and 4.3(i).
If $d(\Delta_{i})>4$ then no further distribution takes place. Suppose without any loss of generality that $d(\Delta_{1})=4$. This cannot happen in case (i); in case (ii) $\Delta_{1}$ is given by Figure 4.3(ii) and so add the $\frac{\pi}{6}$ from $c(\Delta)$ to $c(\Delta_{7})$ across the $bd$ and $bd^{-1}$ edges noting that $d(\Delta_{7})>4$ otherwise $l(\Delta_{7})\in\{bd^{-2}a^{-1},bd^{-1}a^{-1}c^{\pm1},bd^{-1}b^{-1}c^{\pm1}\}$ which contradicts one of the assumptions; in case (iv) $\Delta_{1}$ is given by Figure 4.3(iii) and so add the $\frac{\pi}{6}$ from $c(\Delta)$ to $c(\Delta_{2})$ across the $bd$ and $ad$ edges noting that $d(\Delta_{2})>4$ otherwise $l(\Delta_{2})\in\{ad^{2}b^{-1},ada^{-1}c^{\pm1},ad^{-1}b^{-1}c^{\pm1}\}$ which contradicts one of the assumptions or yields case (ix) or (x); and in case (v) $\Delta$ is given by Figure 4.3(iv) and so add the $\frac{\pi}{6}$ from $c(\Delta)$ to $c(\Delta_{7})$ across the $bd$ and $d^{-1}a^{-1}$ edges noting $d(\Delta_{7})>4$ otherwise $l(\Delta_{7})\in\{d^{-2}a^{-1}b,d^{-1}a^{-1}c^{\pm1}a,d^{-1}a^{-1}c^{\pm1}b\}$ which contradicts one of the assumptions or yields case (xi) or (xii).
Therefore if the region $\hat{\Delta}$ receives positive curvature then it receives $\frac{\pi}{6}$ across each edge and so if $d(\hat{\Delta})\ge6$ then $c^{*}(\hat{\Delta})\le 0$ by Lemma 3.5(i).
This leaves the case when $d(\hat{\Delta})=5$. After checking for vertex labels that contain the sublabels $(bd)$, $(ca)$, $(ad)$ and $(bd^{-1})$ corresponding to the edges crossed in Figures 4.1(ii) and 4.3(i)--(iv) we obtain
$c^{\ast}(\hat{\Delta})\le \pi \left( 2 - \frac{5}{2} \right) +3.\frac{\pi}{6}=0$. This completes cases (i), (ii), (iv) and (v).

Consider case (iii), $bda^{-1}c^{-1}=1$. If $b^{2}=1$ then we obtain a sphere by Lemma 3.1(b)(iii).
Note also that $H$ is non-cyclic in this case otherwise we obtain $b=1$, a contradiction.
Suppose then that $b^{2}\ne1$. First let $l(\Delta)=c^{3}$. Add $\frac{1}{3}c(\Delta)=\frac{\pi}{6}$ to $c(\Delta_{i})$ for $i\in\{1,3,5\}$ as in Figure 4.1(ii).
If $d(\Delta_{i})>4$ then no further distribution takes place. Suppose that $d(\Delta_{1})=4$. Then add $\frac{\pi}{6}$ from $c(\Delta)$ to $c(\Delta_{6})$ across the $bd$ and $ab^{-1}$ edges as shown in Figure 4.4(i) noting that $d(\Delta_{6})>4$, otherwise $l(\Delta_{6})\in\{b^{-1}ab^{-1}a,b^{-1}ad^{\pm2}\}$ which contradicts $b^{2}\ne1$ or one of the assumptions; and if $d(\Delta_{3})=4$ or $d(\Delta_{5})=4$ in Figure 4.1(ii) then similarly add $\frac{\pi}{6}$ to $\Delta_2$ or $\Delta_4$.
Secondly, let $l(\Delta)=d^{3}$. Add $\frac{1}{3}c(\Delta)=\frac{\pi}{6}$ to $c(\Delta_{i})$ for $i\in\{1,3,5\}$ as in Figure 4.3(i).
If $d(\Delta_{i})>4$ then no further distribution takes place. 

\newpage
\begin{figure}
\begin{center}
\psfig{file=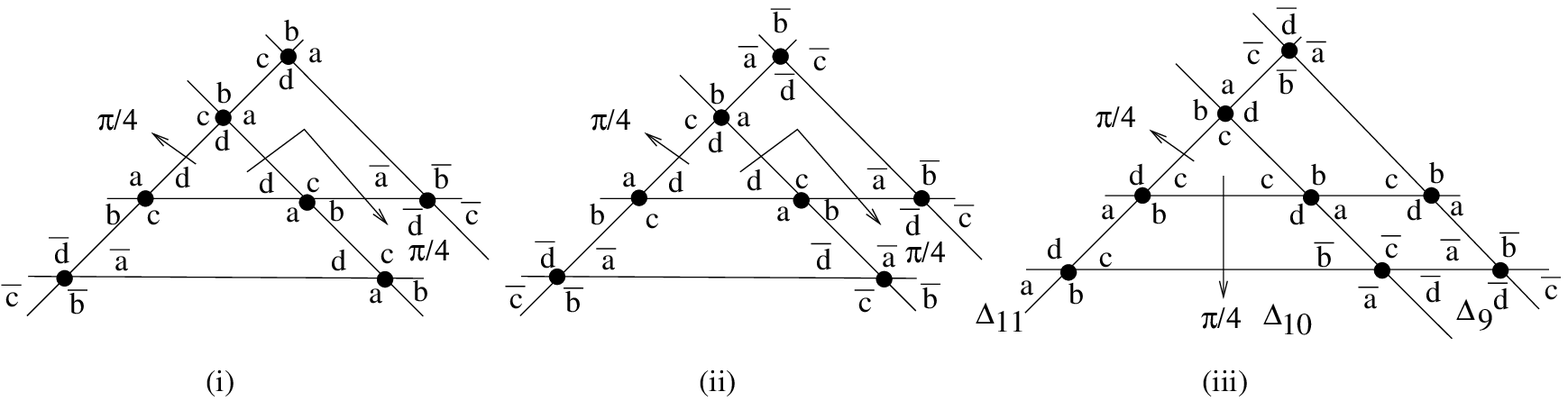}
\end{center}  
\caption{curvature distribution for Case (A2)}
\end{figure}

\noindent Suppose without any loss of generality that $d(\Delta_{1})=4$. Then add $\frac{\pi}{6}$ from $c(\Delta)$ to $c(\Delta_{2})$ across the $ca$ and $ba^{-1}$ edges as shown in Figure 4.4(ii), noting $d(\Delta_{2})>4$, 
otherwise $l(\Delta_{2})\in\{ba^{-1}ba^{-1},ba^{-1}c^{\pm2}\}$ which contradicts $b^{2}\ne1$ or one of the assumptions. If $d(\Delta_{3})=4$ or $d(\Delta_{5})=4$ then
similarly add $\frac{\pi}{6}$ to $\Delta_4$ or $\Delta_6$.
If $\hat{\Delta}$ receives positive curvature and $d(\hat{\Delta})\ge6$,it follows by Lemma 3.5(i) that $c^{*}(\hat{\Delta})\le0$. It remains to check $d(\hat{\Delta})=5$. After checking for vertex labels that contain the sublabels $(bd)$,$(ca)$,$(b^{-1}a)$ and $(ba^{-1})$ corresponding to the edges crossed in Figures 4.1(ii), 4.3(i) and 4.4(i)--(ii) we obtain
$c^{\ast}(\hat{\Delta})\le \pi \left( 2 - \frac{5}{2} \right) +3.\frac{\pi}{6}=0$ or $l(\hat{\Delta})\in\{bda^{-1}ba^{-1},cab^{-1}ab^{-1}\}$ and this contradicts
one of the assumptions. Therefore $c^{*}(\hat{\Delta})\le0$.

Consider case (ix), $bdb^{-1}c=1$ and $cada^{-1}=1$. If $d(\Delta_{i})>4$ for at least two of $\Delta_{i}$ where $i\in\{1,3,5\}$, say $\Delta_{1}$ and $\Delta_{3}$,then add $\frac{1}{2}c(\Delta)=\frac{\pi}{4}$ to each of $c(\Delta_{1})$ and $c(\Delta_{3})$ across the $bd$ and $ca$ edges as shown in Figures 4.1(iii) and 4.4(iii).
By symmetry it can be assumed that $d(\Delta_{1})=d(\Delta_{3})=4$. The two possibilities are given in Figures 4.4(iv) and 4.5(i) and in both cases add $\frac{1}{2} c(\Delta)=\frac{\pi}{4}$ to $c(\Delta_{2})$ across the $bd$, $ad$ or $ca$, $bd^{-1}$ edges as shown. If $d(\Delta_{5})>4$ then add the remaining $\frac{1}{2} c(\Delta)=\frac{\pi}{4}$ to $c(\Delta_{5})$; or if $d(\Delta_{5})=4$ then apply the above to $\Delta_{1}$ and $\Delta_{5}$ to distribute the remaining $\frac{1}{2} c(\Delta)=\frac{\pi}{4}$ similarly to $c(\Delta_{6})$. Now observe that if $\Delta_{1}$ receives positive curvature from $\Delta$ then it does not receive positive curvature from $\Delta_{2}$; and if $\Delta_{2}$ receives positive curvature from $\Delta$ (as in Figures 4.4(iv) and 4.5(i)) then it does not receive positive curvature from $\Delta_{3}$.
It follows that if the region $\hat{\Delta}$ receives positive curvature then it does so across at most two-thirds of its edges and therefore if $\hat{\Delta}$ receives positive curvature and if $d(\hat{\Delta})\ge6$ then $c^{*}(\hat{\Delta})\le 0$ by Lemma 3.5(iv).
Note that $d(\Delta_{2})>4$ in Figures 4.4(iv) and 4.5(i) otherwise $l(\Delta_{2})\in\{c^{-1}ada^{-1},c^{-1}adb^{-1},cbd^{-1}a^{-1},cbd^{-1}b^{-1}\}$ which contradicts one of the assumptions. So there remains the case $d(\hat{\Delta})=5$ and $l(\hat{\Delta})\in\{bd\omega,ca\omega, c^{-1}ad\omega,cbd^{-1}\omega\}$.
Checking shows that in each case $c^{*}(\hat{\Delta})\le \pi (2-\frac{5}{2})+2.\frac{\pi}{4}=0$.

Consider (x), $bdb^{-1}c=1$ and $cad^{-1}a^{-1}=1$. First consider $l(\Delta)=d^{3}$. If at least two of the $\Delta_{i}$ where $i\in\{1,3,5\}$ have degree greater than four, say $\Delta_{1}$ and $\Delta_{3}$, then add $\frac{1}{2} c(\Delta)=\frac{\pi}{4}$ to $c(\Delta_{1})$ and $c(\Delta_{3})$ as shown in Figure 4.4(iii). So assume otherwise and without any loss 

\newpage
\begin{figure}
\begin{center}
\psfig{file=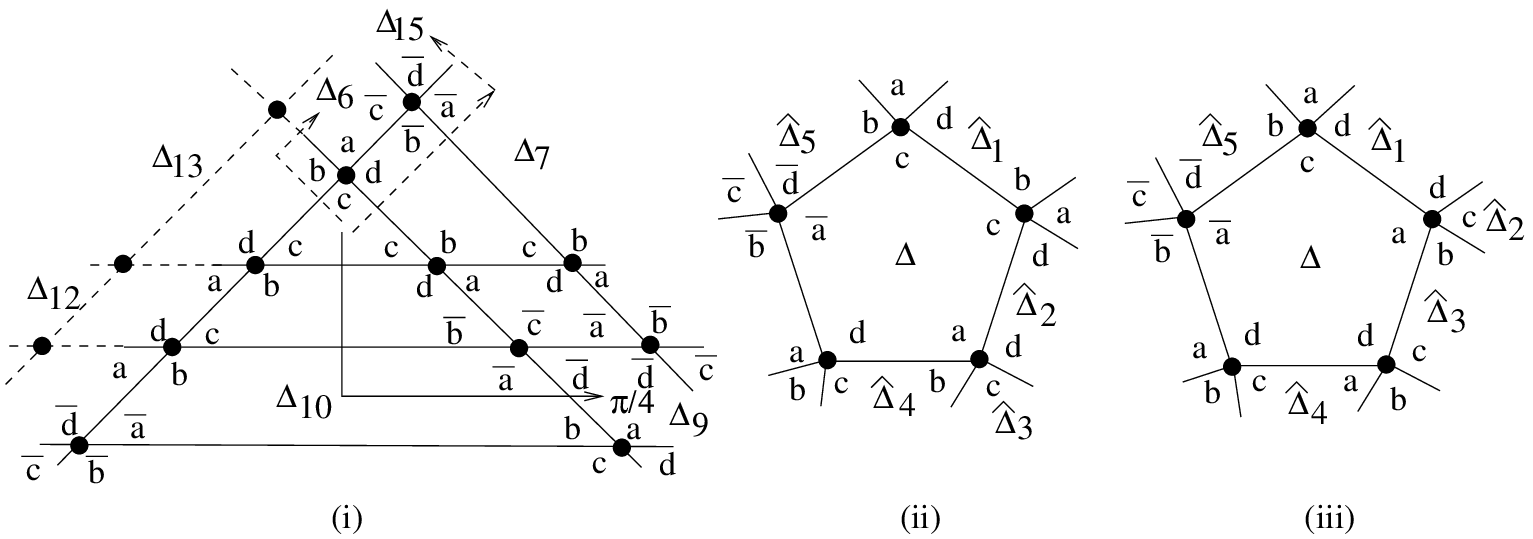}
\end{center}
\caption{curvature distribution for Case (A2)}
\end{figure}

\noindent of generality let $d(\Delta_{1})=d(\Delta_{3})=4$ as shown in Figure 4.5(ii) where $d(\Delta_{2})>4$ otherwise $l(\Delta_{2})=a^{-1}bd^{-2}$ which contradicts
$d^{-1} b^{-1} a \neq 1$. So add $\frac{1}{2} c(\Delta)=\frac{\pi}{4}$ to $c(\Delta_{2})$ as shown in Figure 4.5(ii). If $d(\Delta_{5})>4$ in Figure 4.5(ii) add the remaining $\frac{1}{2} c(\Delta)=\frac{\pi}{4}$ to $c(\Delta_{5})$ 
otherwise use the same argument as above for $\Delta_{1}$ and $\Delta_{5}$ and add $\frac{1}{2} c(\Delta)=\frac{\pi}{4}$ to $c(\Delta_{6})$.
Now consider $l(\Delta)=c^3$.  If at least two of the $\Delta_i$ where $i \in \{ 1,3,5 \}$ have degree $>4$, say, $\Delta_1$ and $\Delta_3$, then add
$\frac{1}{2} c(\Delta)=\frac{\pi}{4}$ to $c(\Delta_1)$ and $c(\Delta_3)$ as in Figure 4.1(iii).
Suppose exactly two of the $\Delta_i$ have degree $=4$, say, $\Delta_1$ and $\Delta_3$.
Add $\frac{1}{2} c(\Delta)=\frac{\pi}{4}$ to $c(\Delta_5)$.
If $d(\Delta_2)>4$ then add the remaining $\frac{1}{2} c(\Delta)=\frac{\pi}{4}$ to $c(\Delta_2)$ as in Figure 4.4(iv).
If $d(\Delta_2)=4$ then add $\frac{1}{2} c(\Delta)=\frac{\pi}{4}$ to $c(\Delta_{10})$ as in Figure 4.5(iii).
If now $d(\Delta_{10})=4$ then $l(\Delta_{10})=ba^{-1} ba^{-1}$ and so add the $\frac{1}{2} c(\Delta)=\frac{\pi}{4}$ to $c(\Delta_9)$ as in Figure 4.6(i).
Observe that $l(\Delta_9)=ad^{-1} d^{-1} w$ forces $d(\Delta_9) > 4$ otherwise there is a contradiction to $d^{-1} b^{-1} a \neq 1$.
Finally suppose that $d(\Delta_i)=4$ for $i \in \{1,3,5 \}$.  Then $\frac{1}{2} c(\Delta)=\frac{\pi}{4}$ is added to either $\Delta_2$, $\Delta_{10}$ or $\Delta_9$ exactly as above; similarly $\frac{1}{2} c(\Delta)= \frac{\pi}{4}$ is added to $\Delta_6$, $\Delta_7$ or $\Delta_{15}$ as shown in Figure 4.6(i).

Now observe that in
Figures 4.4(iii) and 4.1(iii) $\Delta_1$ does not receive positive curvature from $\Delta_2$; in
Figures 4.5(ii) and 4.4(iv) $\Delta_{2}$ does not receive positive curvature from $\Delta_{3}$;
in Figure 4.5(iii) $\Delta_{10}$ does not receive positive curvature from $\Delta_{11}$;
and in Figure 4.6(i) $\Delta_{9}$ does not receive positive curvature from $\Delta_{2}$.
Observe that if $\hat{\Delta}$ receives positive curvature then $d(\hat{\Delta})\ge5$.
It follows from Lemma 3.5(iv) that if $d(\hat{\Delta})\ge6$ then $c^{\ast}(\hat{\Delta})\le0$ so let $d(\hat{\Delta})=5$.
If $\hat{\Delta}$ receives across at most two edges then $c^{*}(\hat{\Delta})\le0$ so it remains to check if $\hat{\Delta}$ receives curvature from more than two edges.
From the above we see that positive curvature is transferred across $(ca),(bd),(bd^{-1}),(ad),(ba^{-1}),(ad^{-1})$-edges.
The only two labels that contain more than two such sublabels and do not yield a contradiction are $a^{-1}ccad$ and $a^{-1}cadd$ as shown in Figures 4.6(ii)--(iii).
Let $l(\Delta)=a^{-1}ccad=1$ as in Figure 4.6(ii).
Here $\Delta$ receives nothing from $\hat{\Delta}_{1}$ or $\hat{\Delta}_{5}$.
If $d(\hat{\Delta}_{2})>3$ then $\Delta$ receives nothing from $\hat{\Delta}_{2}$ and so $c^{*}(\hat{\Delta})\le0$.
If $d(\hat{\Delta}_{2})=3$ then $d(\hat{\Delta}_{3})>3$ and $\Delta$ receives nothing from $\hat{\Delta}_{3}$ via $\hat{\Delta}_{4}$ as in Figure 4.4(iv) and again $c^{*}(\hat{\Delta})\le0$.
Let $l(\Delta)=a^{-1}cadd=1$ as in Figure 4.6(iii).
Here $\Delta$ receives nothing from $\hat{\Delta}_{4}$ or $\hat{\Delta}_{5}$.
If $d(\hat{\Delta}_{1})>3$ then $\Delta$ receives nothing from $\hat{\Delta}_{1}$ 

\newpage
\begin{figure}
\begin{center}
\psfig{file=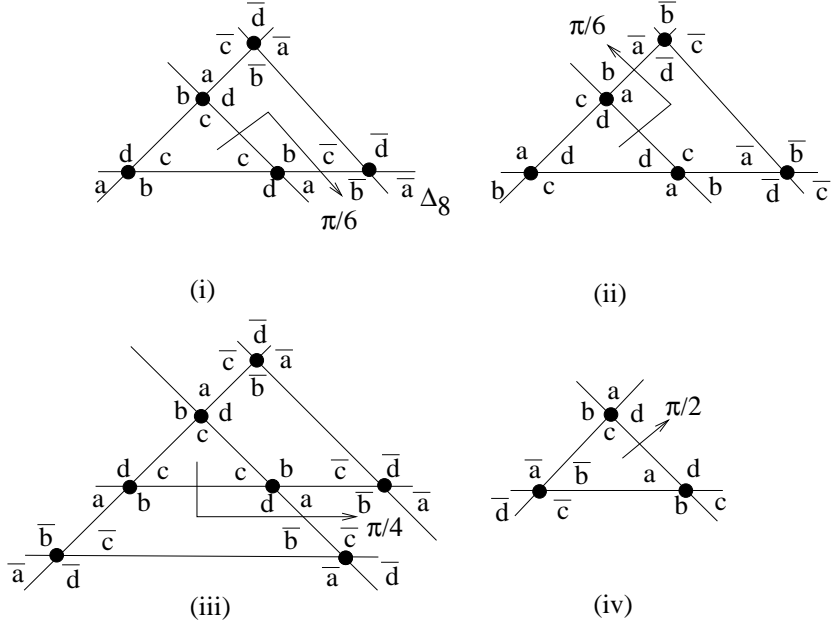}
\end{center}
\caption{curvature distribution for Cases (A2) and (A3)}
\end{figure}

\noindent and so $c^{*}(\hat{\Delta})\le0$.
If $d(\hat{\Delta}_{1})=3$ then $d(\hat{\Delta}_{2})>3$ and $\Delta$ receives nothing from $\hat{\Delta}_{2}$ via $\hat{\Delta}_{3}$ and again $c^{*}(\hat{\Delta})\le0$.

Finally consider case (xii), $bdb^{-1}c^{-1}=cad^{-1}a^{-1}=1$.
Then $c=d$ and so $bdb^{-1}d^{2}=1$.
If now $b^{2}=1$ then we obtain $bdbd^{2}=1$ and $H=\langle bd \rangle$ is cyclic.  Assume first that $H$ is non-cyclic so, in particular, $|b|>2$.
Add $\frac{1}{3}c(\Delta)=\frac{\pi}{6}$ to $c(\Delta_{i})$ for $i\in\{1,3,5\}$ as in Figures 4.1(ii) and 4.3(i) across the $bd$ and $ca$ edges.
If, say, $d(\Delta_{1})>4$ then no further distribution takes place.
If $d(\Delta_{1})=4$ then the $\frac{1}{3} c(\Delta)=\frac{\pi}{6}$ is added to $c(\Delta_{2})$ if $l(\Delta)=c^{3}$ across the $bd$ and $ab^{-1}$ edges, or to $c(\Delta_{6})$ if $l(\Delta)=d^{3}$ across the $ca$ and $a^{-1}b$ edges as shown in Figures 4.7(i)--(ii).
Observe that $d(\Delta_{2})>4$ and $d(\Delta_{6})>4$. If $\hat{\Delta}$ receives positive curvature and $d(\hat{\Delta})\ge6$, it follows by Lemma 3.5(i) that $c^{*}(\hat{\Delta})\le0$. It remains to check $d(\hat{\Delta})=5$. After checking for vertex labels that contain the sublabels $(bd)$,$(ca)$,$(ab^{-1})$ and $(a^{-1}b)$ corresponding to the edges crossed in Figures 4.1(ii), 4.3(i) and 4.7(i)--(ii)
it follows either that $\hat{\Delta}$ receives at most $3.\frac{\pi}{6}$ and so $c^{*}(\hat{\Delta})\le0$ or $l(\hat{\Delta})\in \{bda^{-1}ba^{-1},cab^{-1}ab^{-1}\}$ which in each case yields a contradiction to $H$ non-cyclic.
Now assume that $H$ is cyclic.
If at least two of the $\Delta_i$ where $i \in \{ 1,3,5 \}$ have degree greater than four, say $\Delta_1$ and $\Delta_3$, then add $\frac{1}{2} c(\Delta)=\frac{\pi}{4}$ to each of $c(\Delta_1)$ and $c(\Delta_3)$ as shown in Figures 4.1(iii) and 4.4(iii).
By symmetry assume then that $d(\Delta_1)=d(\Delta_3)=4$.
The two possibilities are in Figures 4.7(iii) and 4.5(ii) and in each case add $\frac{1}{2}c(\Delta)=\frac{\pi}{4}$ to $c(\Delta_2)$ as shown.
If $d(\Delta_5) > 4$ then add the remaining $\frac{1}{2} c(\Delta)=\frac{\pi}{4}$ to $c(\Delta_5)$; or if $d(\Delta_5)=4$ then similarly distribute the remaining $\frac{1}{2} c(\Delta)=\frac{\pi}{4}$ 
via $\Delta_5$ to $\Delta_4$ or $\Delta_6$.  Now observe that if $\Delta_1$ receives positive curvature from $\Delta$ it does not receive curvature from $\Delta_2$; 
and if $\Delta_2$ receives positive curvature from $\Delta$ it does not receive positive curvature from $\Delta_1$ in Figure 4.7(iii) or from $\Delta_3$ in Figure 4.5(ii).   

\newpage
\begin{figure}
\begin{center}
\psfig{file=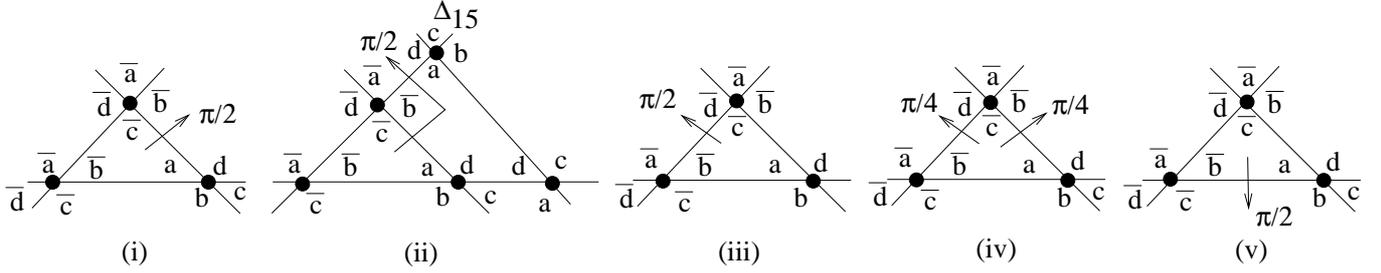}
\end{center}
\caption{curvature distribution for Cases (A4), (A6) and (A8)}
\end{figure}

\noindent It follows that if $d(\hat{\Delta}) \geq 6$ then $c^{\ast}(\Delta) \leq 0$ by Lemma 3.5(iv).  Now $d(\Delta_2) > 4$ in Figures 4.7(iii) and 4.5(ii) so there remains the case $d(\hat{\Delta})=5$ and $d(\hat{\Delta}) \in \{ 
bdw, caw, c^{-1} 
aw, bd^{-1} w \}$.  But checking shows that in all cases $c^{\ast} (\hat{\Delta}) \leq \pi \left( 2- \frac{5}{2} \right) + 2. \frac{\pi}{4} = 0$.

In conclusion $\mathcal{P}$ fails to be aspherical in this case if and only if $b^2=bda^{-1} c^{-1}=1$ (and $H$ is non-cyclic).

\textbf{(A3)} $\boldsymbol{|c|>3}$, $\bs{|d|>3}$, $\bs{a^{-1}b\ne1}$, $\bs{cab^{-1}=1}$, $\bs{c^{-1}ab^{-1}\ne1}$, $\bs{d^{\pm1}b^{-1}a\ne1}$.\newline
If $d(\Delta)=3$ then $l(\Delta)= cab^{-1}$ and $\Delta$ is given by Figure 4.7(iv).
Add $c(\Delta)=\frac{\pi}{2}$ to $c(\Delta_{1})$ across the $d^{2}$ edge. If $\Delta_{5}$ receives no curvature across at least four edges then there are at least four \textit{gaps} as defined in Section 3 and $c^{*}(\Delta_{5})\le 0$ by Lemma 3.6. This leaves the case when $l(\Delta_{1})\in\{d^{n},d^{n}a^{-1}b,d^{n}b^{-1}a\}$. If $l(\Delta_{5})=d^{n}$ then there is a sphere by Lemma 3.1(a)(i); and if $d^n (a^{-1}b)^{\pm 1}=1$ then $H$ is cyclic.
Therefore $\mathcal{P}$ is aspherical if and only if $|d| = \infty$.

\textbf{(A4)} $\boldsymbol{|c|>3}$, $\bs{|d|>3}$, $\bs{a^{-1}b\ne1}$, $\bs{cab^{-1}\ne1}$, $\bs{c^{-1}ab^{-1}=1}$, $\bs{d^{\pm1}b^{-1}a\ne1}$.\newline
If $d(\Delta)=3$ then $l(\Delta)=c^{-1}ab^{-1}$ and $\Delta$ is given by Figure 4.8(i).
First assume that $H$ is non-cyclic.
Add $c(\Delta)=\frac{\pi}{2}$ to $c(\Delta_{1})$ across the $db^{-1}$ edge and note that $l(\Delta_{1})=db^{-1}\omega$ forces $d(\Delta_{1})\ge5$, otherwise $l(\Delta_{1})\in\{db^{-1}ad,db^{-1}c^{\pm1}a,db^{-1}c^{\pm1}b\}$ which contradicts $|d|>3$ or $H$ non-cyclic. Observe also that $\Delta_{1}$ does not receive curvature from $\Delta_{2}$ or $\Delta_{6}$ and so if $d(\Delta_1)\ge7$, it follows that $c^{*}(\Delta_{1})\le 0$ by Lemma 3.5(vi); if $d(\Delta_{1})=5$ then $l(\Delta_{1})$ contains at most one occurrence of $(db^{-1})^{\pm1}$ and so $c^{*}(\Delta_1)\le \pi \left( 2 - \frac{5}{2} \right)+\frac{\pi}{2}=0$; and if $d(\Delta_{1})=6$ then $l(\Delta_{1})$ contains at most two occurrence of $(db^{-1})^{\pm1}$ and so $c^{*}(\Delta_1)\le \pi \left( 2 - \frac{6}{2} \right)+2.\frac{\pi}{2}=0$.  Now let $H$ be cyclic.
If $d(\Delta_{1})>4$ in Figure 4.8(i) then add $c(\Delta)=\frac{\pi}{2}$ to $c(\Delta_{1})$. Assume otherwise and $d(\Delta_{1})=4$. Then $l(\Delta_{1})\in \{db^{-1}ad,db^{-1} c^{\pm 1} a,db^{-1}c^{\pm 1} b\}$ which yields a contradiction to the (\textbf{A4}) assumptions except when $l(\Delta_{1})\in\{db^{-1}ad,db^{-1}ca\}$.
If $db^{-1} ca=1$ then there is a sphere by Lemma 3.1(a)(ii) so assume otherwise.
Thus if $d(\Delta_{1})=4$ then $l(\Delta_{1})=db^{-1}ad$.
In this case as shown in Figure 4.8(ii) add $c(\Delta)=\frac{\pi}{2}$ firstly to $c(\Delta_{1})$ then on to $c(\Delta_{6})$. Observe that $d(\Delta_{6})>4$ otherwise
$l(\Delta_{6})\in\{da^{-1}bd,da^{-1}c^{\pm1}a,da^{-1}cb,da^{-1}c^{-1}b\}$ which either contradicts $|c|>3$, $|d|>3$ or $d^{\pm1}b^{-1}a\ne1$, or implies the exceptional case (\textbf{E1}) when $da^{-1}c^{-1}b=1$.
Now observe that if $\Delta_{1}$ receives positive curvature from $\Delta$ then it does not receive positive curvature from $\Delta_{2}$ or $\Delta_{6}$ in Figure 4.8(i); and if $\Delta_{6}$ receives positive curvature from $\Delta$ then it does not receive positive curvature from $\Delta_{5}$ or $\Delta_{15}$ in Figure 

\newpage
\begin{figure}
\begin{center}
\psfig{file=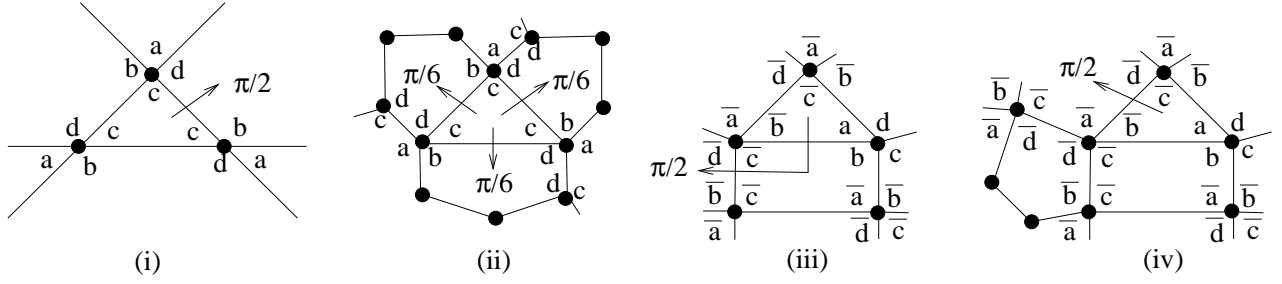}
\end{center}
\caption{curvature distribution for Case (A6)}
\end{figure}

\noindent 4.8(ii).
Therefore if $\hat{\Delta}$ receives positive curvature then it does so across at most half of its edges which implies that if $d(\hat{\Delta})\ge7$, then  $c^{*}(\hat{\Delta})\le 0$ by Lemma 3.5(vi).
If $d(\hat{\Delta})=5$ then checking shows that $l(\hat{\Delta})$ contains at most one occurrence of $(db^{-1})^{\pm1}$ or $(da^{-1})^{\pm1}$. It follows that $c^{*}(\hat{\Delta})\le 0$. If $d(\hat{\Delta})=6$ then checking shows that $l(\hat{\Delta})$ contains at most two occurrences of $(db^{-1})^{\pm1}$ or $(da^{-1})^{\pm1}$ and so $c^{*}(\hat{\Delta})\le 0$.
In conclusion in this case $\mathcal{P}$ fails to be aspherical if and only if $db^{-1}ca=1$ (and so $H$ is cyclic).

\textbf{(A5)}  $\boldsymbol{|c|=3}$, $\bs{|d|>3}$, $\bs{a^{-1}b\ne1}$, $\bs{cab^{-1}=1}$, $\bs{c^{-1}ab^{-1}\ne1}$, $\bs{d^{\pm1}b^{-1}a\ne1}$.\newline
If $|d|<\infty$ then there is a sphere by Lemma 3.1(a)(i) and if $|d| = \infty$ then $\mathcal{P}$ is aspherical by Lemma 3.4(ii).

\textbf{(A6)}  $\boldsymbol{|c|=3}$, $\bs{|d|>3}$, $\bs{a^{-1}b\ne1}$, $\bs{cab^{-1}\ne1}$, $\bs{c^{-1}ab^{-1}=1}$, $\bs{d^{\pm1}b^{-1}a\ne1}$.\newline
If $d(\Delta)=3$ then $\Delta$ is given by Figure 4.1(ii) and 4.8(iii). First assume that either $H$ is non-cyclic or $H$ is cyclic but $b \neq d^{\pm 2}$.
If $l(\Delta)=c^{3}$ then add $\frac{1}{3} c(\Delta)=\frac{\pi}{6}$ to $c(\Delta_{i})$ for $i\in\{1,3,5\}$ across the $bd$ edge. Note that if $d(\Delta_{i})=4$ then $l(\Delta_{i})\in\{bdda^{-1},bda^{-1}c^{\pm1},bdb^{-1}c^{\pm1}\}$ which contradicts $|d| > 3$ or $b \neq d^{\pm 2}$ and so $d(\Delta_{i})\ge5$.
If $l(\Delta)=c^{-1}ab^{-1}$ then add $c(\Delta)=\frac{\pi}{2}$ to $c(\Delta_{5})$ across the $d^{-1}a^{-1}$ edge and note that if $d(\Delta_{5})=4$ then $l(\Delta_{5})\in\{d^{-1}a^{-1}bd^{-1},d^{-1}a^{-1}c^{\pm1}a,d^{-1}a^{-1}c^{\pm1}b\}$ which contradicts $|d| > 3$ or $b \neq d^{\pm 2}$ and so $d(\Delta_{5})\ge 5$.
Observe now that $\Delta_{1}$ does not receive curvature from $\Delta_{2}$ or $\Delta_{6}$ in Figure 4.1(ii); and $\Delta_{5}$ does not receive curvature from $\Delta_{4}$ or $\Delta_{6}$ in Figure 4.8(iii).
Therefore if $\hat{\Delta}$ receives positive curvature then it does so across at most half of its edges which implies that if $d(\hat{\Delta})\ge7$, then  $c^{*}(\hat{\Delta})\le 0$ by Lemma 3.5(vi).
It remains to consider $5\le d(\hat{\Delta})\le6$. If $d(\hat{\Delta})=5$ then checking shows that $l(\hat{\Delta})$ contains at most one occurrence of $(bd)^{\pm1}$ or $(d^{-1}a^{-1})^{\pm1}$. It follows that $c^{*}(\hat{\Delta})\le 0$. If $d(\hat{\Delta})=6$ then checking shows that $l(\hat{\Delta})$ contains at most two occurrences of $(bd)^{\pm1}$ or $(d^{-1}a^{-1})^{\pm1}$ and so $c^{*}(\hat{\Delta})\leq 0$ and $\mathcal{P}$ is aspherical.

Now let $b=d^{-2}$.  This implies $c=d^2$ and $d^6=1$.  It follows that if
$d(\Delta) \in \{ 4,5 \}$ then $l(\Delta) \in \{ bdda^{-1}, ccab^{-1}, cab^{-1} ab^{-1}, dda^{-1} c^{-1} a,
ddb^{-1} ca, ddb^{-1} c^{-1} b \}$.  If $l(\Delta)=c^{-1} ab^{-1}$ then add $\frac{1}{2} c(\Delta)=\frac{\pi}{4}$ to
$c(\Delta_1)$ and to $c(\Delta_5)$ as shown in Figure 4.8(iv).
Let $l(\Delta)=c^3$.  If at least two of $d(\Delta_i) > 4$, say $i=1,3$, then add $\frac{1}{2}c(\Delta)=\frac{\pi}{4}$ to $c(\Delta_1)$ and to $c(\Delta_3)$ as shown in Figure 4.1(iii); if say $d(\Delta_1)=d(\Delta_3)=4$ and $d(\Delta_5)>4$ then add $\frac{1}{2}c(\Delta)=\frac{\pi}{4}$ to $c(\Delta_5)$ and $\frac{1}{2}c(\Delta)=\frac{\pi}{4}$ to $c(\Delta_2)$ via $\Delta_3$ as shown in Figure 4.1(iv); and if $d(\Delta_i)=4$ for $i \in \{ 1,3,5 \}$ then add $\frac{1}{3}c(\Delta)=\frac{\pi}{6}$ to each of $c(\Delta_j)$ for $j \in \{ 2,4,6 \}$ similarly.  Therefore if $\hat{\Delta}$ receives  

\newpage
\begin{figure}
\begin{center}
\psfig{file=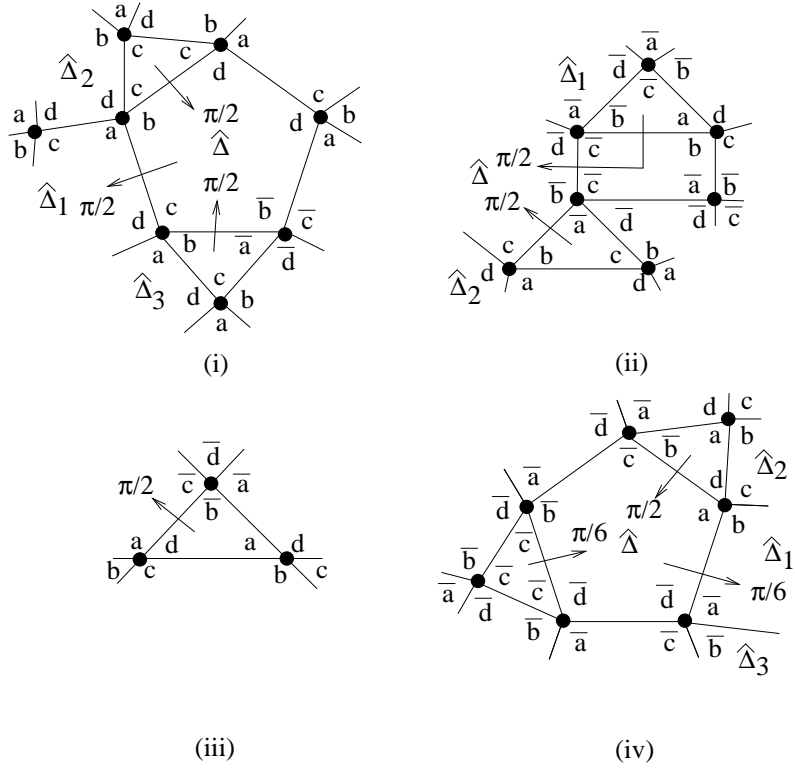}
\end{center}
\caption{curvature distribution for Cases (A6), (A7), (B5) and (B6)}
\end{figure}

\noindent positive curvature then $d(\hat{\Delta}) \geq 5$.  Observe that $\Delta_1$ does not receive positive curvature from $\Delta_2$ in Figure 4.1(iii); $\Delta_1,\Delta_5$ does not receive positive curvature from 
$\Delta_6,\Delta_4$ 
respectively in Figure 4.8(iv); and $\Delta_2$ does not receive positive curvature from $\Delta_9$ in Figure 4.1(iv).  It follows that if $\hat{\Delta}$ receives positive curvature then it does so across at most two-thirds of its edges.  Therefore $d(\hat{\Delta}) \geq 6$ implies $c^{\ast}(\hat{\Delta}) \leq 0$ by Lemma 3.5(iv). Since there are at most two occurrences of $(ad)^{\pm 1}$, $(db^{-1})^{\pm 1}$, $(bd)^{\pm 1}$, $(cad^{-1})^{\pm 1}$ in each of $cab^{-1} ab^{-1}$, $dda^{-1} c^{-1} a$, $ddb^{-1} ca$, $ddb^{-1} c^{-1} b$ it follows that if $d(\hat{\Delta})=5$ then again $c^{\ast}(\hat{\Delta}) \leq 0$.

Finally let $b=d^2$.  This implies $c=d^{-2}$ and $d^6=1$.  It follows that if $d(\Delta) \in \{4,5\}$ then
$l(\Delta) \in \{ d^2 b^{-1} a, c^2 ab^{-1}, bd^2 b^{-1} c, bd^2 a^{-1} c^{-1}, d^2 a^{-1} ca, cab^{-1} ab^{-1} \}$.
Let $c(\Delta)=c^3$.  If at least one of $\Delta_i$ where $i \in \{ 1,3,5 \}$ has degree $>5$, say $d(\Delta_1)>5$, then add $c(\Delta)=\frac{\pi}{2}$ to $c(\Delta_1)$ as shown in Figure 4.9(i); or if $d(\Delta_i)=5$ for $i\in \{1,3,5\}$ then add $\frac{1}{3}c(\Delta)=\frac{\pi}{6}$ to each $c(\Delta_i)$ as shown in Figure 4.9(ii) in which
$l(\Delta_i) \in \{ bd^2 b^{-1} c, bd^2 a^{-1} c^{-1} \}$.  Let $c(\Delta)=c^{-1}ab^{-1}$.  If $d(\Delta_3) > 4$ then add $c(\Delta)=\frac{\pi}{2}$ to $c(\Delta_3)$ as in Figure 4.8(v); if $d(\Delta_3)=4$ and $d(\Delta_4)>5$ then add
$c(\Delta)=\frac{\pi}{2}$ to $c(\Delta_4)$ via $\Delta_3$ as in Figure 4.9(iii); and if $d(\Delta_3)=4$ and $d(\Delta_4)=5$ then add $c(\Delta)=\frac{\pi}{2}$ to $c(\Delta_5)$ as in Figure 4.9(iv).  Finally if $d(\hat{\Delta})=5$, $l(\hat{\Delta})=
bd^2 b^{-1} c$ and $\hat{\Delta}$ receives positive curvature from two neighbouring regions then add $\frac{\pi}{2}$ from $c^{\ast}(\hat{\Delta})$ to $c(\hat{\Delta}_1)$ as shown in Figure 4.10(i).  Now observe that $\Delta_1$ does not receive positive curvature from $\Delta_6$ in Figure 4.9(i); $\Delta_1$, for example, does not receive from $\Delta_6$ or $\Delta_2$ in Figure 4.9(ii); $\Delta_3$ does not receive from $\Delta_4$ in Figure 4.8(v); $\Delta_4$ does not receive from $\Delta_5$ in Figure 4.9(iii); $\Delta_5$ does not receive from $\Delta_4$ or $\Delta_6$ in Figure 4.9(iv); and $\hat{\Delta}_1$ does not receive from $\hat{\Delta}_2$ or $\hat{\Delta}_3$ in Figure 4.10(i).
It follows that if $\hat{\Delta}$ receives positive curvature across two consecutive edges then $d(\hat{\Delta})>5$ and $\hat{\Delta}$ is given by Figure 4.10(ii) in which $\hat{\Delta}$ does not receive from $\hat{\Delta}_1$ or $\hat{\Delta}_2$.  Since two sublabels of the form $(c^{-1} bd)^{\pm 1}$ must be separated by a sublabel of length at least 2 it follows that if $\hat{\Delta}$ receives positive curvature then it does so across at most $\frac{3}{5}$ of its edges.  Therefore if $d(\hat{\Delta}) \geq 8$ then $c^{\ast}(\hat{\Delta}) \leq 0$ by Lemma 3.5(vii); and if $l(\hat{\Delta})$ does not involve $(c^{-1}bd)^{\pm 1}$ then Lemma 3.5(vi) applies and $c^{\ast}(\hat{\Delta})\leq 0$ for $d(\hat{\Delta}) \geq 7$.  If $d(\hat{\Delta})=5$ then either $c^{\ast}(\hat{\Delta}) \leq 0$ or $\hat{\Delta}$ is given by Figure 4.10(i) and $c^{\ast}(\hat{\Delta}) - \frac{\pi}{2} \leq 0$.  If $l(\hat{\Delta})$ does not involve
$(c^{-1} bd)^{\pm 1}$ it remains to consider $d(\hat{\Delta})=6$.  But checking shows that each candidate for
$l(\hat{\Delta})$ contains at most two occurrences of $(bd)^{\pm 1}$, $(b^{-1} c)^{\pm 1}$ or $(cad)^{\pm 1}$ so
$c^{\ast}(\hat{\Delta}) \leq 0$.  Finally if $l(\hat{\Delta})=c^{-1} bdw$ and $d(\hat{\Delta})=6$ then $l(\hat{\Delta})=c^{-1}bd^2 b^{-1} c^{-1}$ and $c^{\ast}(\hat{\Delta}) \leq 0$; or if $d(\hat{\Delta})=7$ then $l(\hat{\Delta}) \in \{
c^{-1} bd^4 b^{-1}, c^{-1} bda^{-1} bdb^{-1}, c^{-1} bda^{-1} bd^{-1} a^{-1}, c^{-1} bdb^{-1} ad^{-1} b^{-1} \}$ and
$c^{\ast} (\hat{\Delta}) \leq c(\hat{\Delta}) + 3( \frac{\pi}{2} )=0$.
In conclusion $\mathcal{P}$ is aspherical for this case.

\textbf{(A7)} $\boldsymbol{|c|=3}$, $\bs{|d|>3}$, $\bs{a^{-1}b\ne1}$, $\bs{c^{\pm1}ab^{-1}\ne1}$, $\bs{db^{-1}a=1}$, $\bs{d^{-1}b^{-1}a\ne1}$.\newline
If $d(\Delta)=3$ then $l(\Delta) \in \{ c^3, db^{-1} a \}$ and $\Delta$ is given by Figures 4.1(ii) and 4.10(iii).
Note that if $l(\Delta_i)=bdw$, $c^{-1} aw$ and $d(\Delta_i)=4$ then $l(\Delta_i)=bda^{-1}c^{\pm 1}$,
$c^{-1} ad^{-1} b^{-1}$ (respectively) otherwise there is a contradiction to $|d| > 3$.  First assume that
$bda^{-1} c^{\pm 1} \neq 1$.  If $l(\Delta)=c^3$ then add $\frac{1}{3} c(\Delta)=\frac{\pi}{6}$ to $c(\Delta_i)$ for $i \in \{ 1,3,5 \}$ across the $bd$ edge as shown in Figure 4.1(ii); and if $l(\Delta)=db^{-1} a$ then add
$c(\Delta)=\frac{\pi}{2}$ to $c(\Delta_5)$ as in Figure 4.10(iii).  Observe then that $\Delta_1$ does not receive positive curvature from $\Delta_2$ or $\Delta_6$ in Figure 4.1(ii); and $\Delta_5$ does not receive from $\Delta_4$ or $\Delta_6$ in Figure 4.10(iii).  It follows that $c^{\ast}(\hat{\Delta}) \leq 0$ for $d(\hat{\Delta}) \geq 7$ by Lemma 3.5(vi).
If $d(\hat{\Delta})=5$ then checking shows that either $l(\hat{\Delta})$ contains at most one occurrence of $(bd)^{\pm1}$ or $(c^{-1}a)^{\pm1}$ and so  $c^{*}(\hat{\Delta})\le 0$ or $l(\hat{\Delta}) = c^{-1}ad^{-2}b^{-1}$.
If $d(\hat{\Delta})=6$ then checking shows that $l(\hat{\Delta})$ contains at most two occurrences of $(bd)^{\pm1}$ or $(c^{-1}a)^{\pm1}$ and so $c^{*}(\hat{\Delta})\le 0$.
Thus if $c^{\ast}(\hat{\Delta}) > 0$ then $\hat{\Delta}$ is given by Figure 4.10(iv) and add $c^{\ast}(\hat{\Delta}) = \frac{\pi}{6}$ to $c(\hat{\Delta}_1)$ as shown.  Now $\hat{\Delta}_1$ does not receive positive curvature from $\hat{\Delta}_2$ or $\hat{\Delta}_3$ in Figure 4.10(iv) and the above statements still hold for $\Delta_i$ of Figure 4.1(ii) and $\Delta_5$ of Figure 4.10(iii).
Therefore if $\hat{\Delta}$ receives positive curvature it does so across at most half of its edges and so
$c^{\ast}(\hat{\Delta}) \leq 0$ for $d(\hat{\Delta}) \geq 7$.
Let $d(\hat{\Delta}) < 7$.
If $l(\hat{\Delta})$ does not involve $a^{-1} b$ then checking shows that $c^{\ast}(\hat{\Delta}) \leq 0$ so assume $l(\hat{\Delta})=a^{-1}bw$.
But if now $d(\hat{\Delta}) < 6$ then $l(\hat{\Delta})$ together with $c^{-1} ad^{-2} b^{-1}$ contradicts $|d| =9$; and if $d(\hat{\Delta})=6$ then the only labels with more than two occurrences of $(a^{-1}b)^{\pm 1}$, $(bd)^{\pm 1}$ or
$(c^{-1} a)^{\pm 1}$ that do not yield a contradiction are $a^{-1} ba^{-1} cbd$, $a^{-1} bda^{-1} cb$ and
$a^{-1} b d^{-1} b^{-1} ad$ in which case $c^{\ast} (\hat{\Delta}) \leq c(\Delta) + 2(\frac{\pi}{6}) + \frac{\pi}{2} < 0$.
If $bda^{-1} c=1$ then, since $db^{-1} a \leftrightarrow c^{-1} ab^{-1}$ and $bda^{-1}c \leftrightarrow a^{-1} c^{-1} bd^{-1}$, there is a sphere by Lemma 3.1(a)(ii).  This leaves the case $bda^{-1} c^{-1}=1$ which implies $d^6=1$ and $c=d^2$.  If $l(\Delta)=db^{-1} a$ then again add $c(\Delta)=\frac{\pi}{2}$ for $c(\Delta_5)$ as in Figure 4.10(iii).
Let $l(\Delta)=c^3$.  If at least two $d(\Delta_i) \geq 5$, say $\Delta_1$ and $\Delta_3$, then add
$\frac{1}{2} c(\Delta)=\frac{\pi}{4}$ to $c(\Delta_1)$ and to $c(\Delta_3)$ as in Figure 4.1(iii).  If, say,
$d(\Delta_1)=d(\Delta_3)=4$ then add $\frac{1}{2} c(\Delta)=\frac{\pi}{4}$ to $c(\Delta_2)$ via $\Delta_1$ as in Figure 4.2(iii) and $\frac{1}{2} c(\Delta)=\frac{\pi}{4}$ to $c(\Delta_5)$ if $d(\Delta_5) \geq 5$ or to $c(\Delta_6)$ via
$\Delta_5$ if $d(\Delta_5)=4$ also.  Observe that $\Delta_5$ does not receive positive curvature from $\Delta_4$ or $\Delta_6$ in Figure 4.10(iii); $\Delta_1$ does not receive any from $\Delta_6$ in Figure 4.1(iii); 
and $\Delta_2$ does not receive any from $\Delta_3$ in  

\newpage
\begin{figure}
\begin{center}
\psfig{file=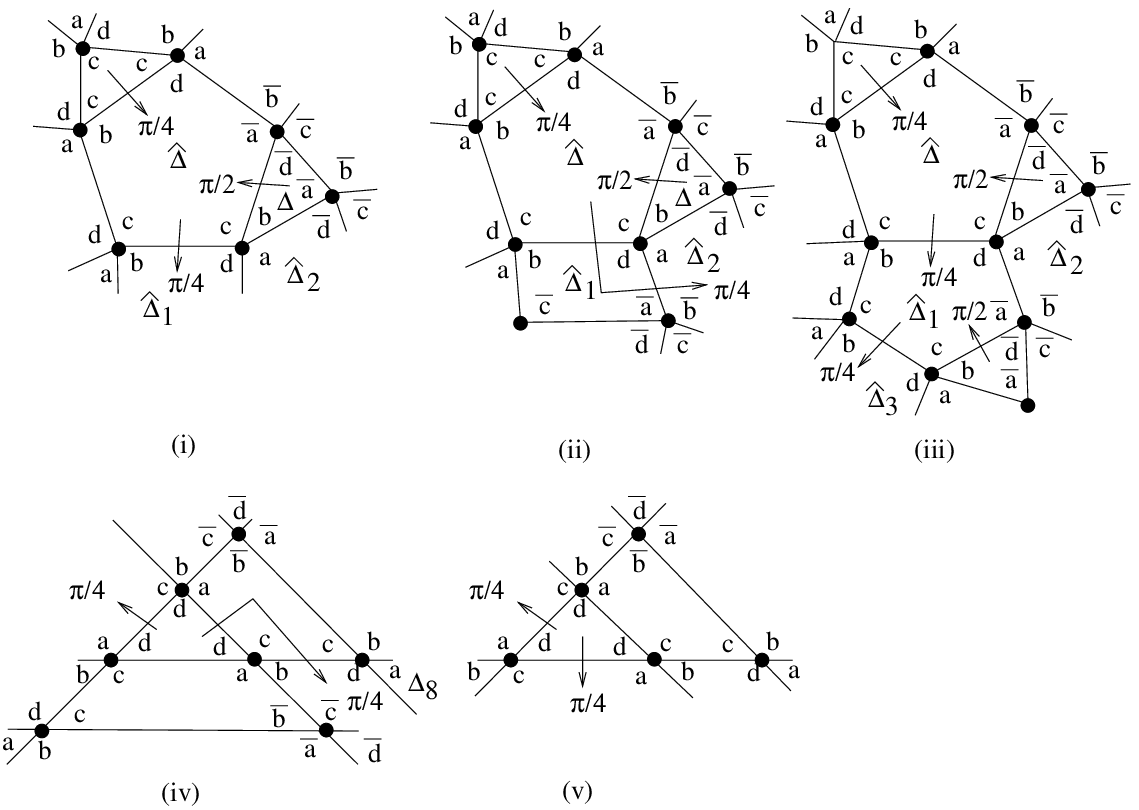}
\end{center}
\caption{curvature distribution for Cases (A7) and (A8)}
\end{figure}

\noindent Figure 4.2(iii).  It follows from Lemma 3.5(v) that if $d(\hat{\Delta}) \geq 7$ then $c^{\ast}(\hat{\Delta}) \leq 0$.  If $d(\hat{\Delta})=6$ then checking the allowable labels shows that
$l(\hat{\Delta})$ involves at most two occurrences of $(bd)^{\pm 1}$, $(c^{-1}a)^{\pm 1}$ or $(b^{-1}ab)^{\pm 1}$ and so $c^{\ast}(\hat{\Delta}) \leq 0$.  If $d(\hat{\Delta})=5$ then $c^{\ast}(\hat{\Delta}) \leq 0$ except when $\hat{\Delta}$ is given by Figure 4.11(i) in which case add $c^{\ast}(\hat{\Delta}) = \frac{\pi}{4}$ to $c(\hat{\Delta}_1)$ across the $bd$ edge as shown; and note that $\hat{\Delta}_1$ does not receive positive curvature from $\hat{\Delta}_2$ (and so none of the above is affected).  If $c^{\ast}(\hat{\Delta}_1) \leq 0$ then we are done, otherwise
$l(\hat{\Delta}_1) \in \{ bda^{-1} c^{-1}, bda^{-1} c^2 \}$.  Let $l(\hat{\Delta}_1)=bda^{-1}c^{-1}$.  Then add
$c^{\ast}(\hat{\Delta})=\frac{\pi}{4}$ to $\hat{\Delta}_2$ via $\hat{\Delta}_1$ as shown in Figure 4.11(ii) noting that $\hat{\Delta}_2$ does not receive positive curvature from the region $\Delta$.  It follows that if $c^{\ast}(\hat{\Delta}) > 0$ then $l(\hat{\Delta})=b^{-1} ad^{-1} w$ and $d(\hat{\Delta}) \leq 6$.  But this implies that $l(\hat{\Delta}) \in \{ b^{-1} ad^{-1} a^{-1} bd, b^{-1} ad^{-1} a^{-1} ca \}$ and $c^{\ast}(\hat{\Delta}) \leq 0$.  Finally if $\hat{\Delta}_1$ is given by Figure 4.11(iii) then add $c^{\ast}(\hat{\Delta}_1)=\frac{\pi}{4}$ to $\hat{\Delta}_3$ as shown and repeat the above argument.  This procedure must terminate at a region $\hat{\Delta}_k$, say, where either $c^{\ast}(\hat{\Delta}_k) \leq 0$ or $d(\hat{\Delta}_k)=4$ as in Figure 4.11(ii), in which case
$c^{\ast}(\hat{\Delta}_k)=\frac{\pi}{4}$ is added to $c(\hat{\Delta}_{k+1})$ where $l(\hat{\Delta}_k)=b^{-1}ad^{-1}w$ and so $c^{\ast}(\hat{\Delta}_k) \leq 0$.
In conclusion $\mathcal{P}$ is aspherical if and only if $bda^{-1}c \neq 1$.  (Note that if $bda^{-1}c=1$ then $H$ is cyclic.)

\textbf{(A8)} $\boldsymbol{|c|=3}$, $\bs{|d|=3}$, $\bs{a^{-1}b\ne1}$, $\bs{cab^{-1}\ne1}$, $\bs{c^{-1}ab^{-1}=1}$, $\bs{d^{\pm1}b^{-1}a\ne1}$.\newline
If $d(\Delta)=3$ then $l(\Delta)\in\{c^{3}, d^{3}, c^{-1}ab^{-1}\}$ and $\Delta$ is given by
Figures 4.1(ii), 4.3(i) and 4.8(iii).
If $l(\Delta)=c^{3}$ then add $\frac{1}{3} c(\Delta)=\frac{\pi}{6}$ to $c(\Delta_{i})$ for $i\in\{1,3,5\}$ across the $bd$ edge as shown in Figure 4.1(ii). Note that if $d(\Delta_{i})=4$ then $l(\Delta_{i})\in\{bdda^{-1},bda^{-1}c^{\pm1},bdb^{-1}c^{\pm1}\}$ which contradicts
$|d|=3$ or $d^{\pm 1} b^{-1} a \neq 1$ and so $d(\Delta_{i})\ge5$.
If $l(\Delta)=c^{-1}ab^{-1}$ then add $c(\Delta)=\frac{\pi}{2}$ to $c(\Delta_{5})$ across the $d^{-1}a^{-1}$ edge as shown in Figure 4.8(iii) and note that if $d(\Delta_{5})=4$ then $l(\Delta_{5})\in\{d^{-1}a^{-1}bd^{-1},d^{-1}a^{-1}c^{\pm1}a,d^{-1}a^{-1}c^{\pm1}b\}$ which again contradicts $|d|=3$ or $d^{\pm 1} b^{-1} a \neq 1$ and so $d(\Delta_{5})\ge 5$.
If $l(\Delta)=d^{3}$ then there are three subcases.
First assume that $d(\Delta_{i})\ge5$ for $i\in\{1,3,5\}$. Then add $\frac{1}{3}c(\Delta)= \frac{\pi}{6}$ to each $c(\Delta_{i})$ across the $ca$ edge as shown in Figures 4.3(i). Now suppose that at least two of the $\Delta_{i}$ have degree four, without loss of generality $\Delta_{1}$ and $\Delta_{3}$. Then $l(\Delta_{1})=l(\Delta_{3})=cab^{-1}c$ (otherwise we will have a contradiction to $d^{\pm 1} b^{-1} a \neq 1$) which in turn forces $l(\Delta_{2})=c^{-1}bd\omega$ and so $d(\Delta_{2})>4$ otherwise $l(\Delta_{2})\in\{c^{-1}bda^{-1},c^{-1}bdb^{-1}\}$ which contradicts $|d|=3$ or $d^{\pm 1} b^{-1} a \neq 1$. In this case add $\frac{1}{2} c(\Delta)=\frac{\pi}{4}$ to $c(\Delta_{2})$ and $c(\Delta_{5})$ across the $ca$ and $bd$ edges as shown in Figure 4.11(iv)
and if $d(\Delta_5)=4$ then in the same way add $\frac{1}{2} c(\Delta)=\frac{\pi}{4}$ to $c(\Delta_6)$. This leaves the case where exactly one of $\Delta_{i}$ has degree four and without loss assume $d(\Delta_{1})=4$. Then add $\frac{1}{2} c(\Delta)=\frac{\pi}{4}$ to each of $c(\Delta_{3})$ and $c(\Delta_{5})$ as shown in Figure 4.11(v).

Observe that in Figure 4.1(ii) $\Delta_{1}$ does not receive any positive curvature from $\Delta_{2}$ or $\Delta_{6}$; in Figure 4.3(i) $\Delta_{1}$ does not receive positive curvature from $\Delta_{2}$; in Figure 4.8(iii) $\Delta_{5}$ does not receive any positive curvature from $\Delta_{6}$; in Figure 4.11(iv) $\Delta_{2}$ does not receive any positive curvature from $\Delta_{3}$ or $\Delta_{8}$ and $\Delta_{5}$ does not receive any positive curvature from $\Delta_{6}$; and in Figure 4.11(v) $\Delta_{3}$, $\Delta_{5}$ (respectively) does not receive any positive curvature from $\Delta_{4}$, $\Delta_{6}$ (respectively).
Let $\hat{\Delta}$ receive positive curvature and first suppose that $l(\hat{\Delta})$ does not involve $(ad)^{\pm1}$.
It follows that $\hat{\Delta}$ receives positive curvature across at most half of its edges and so, since $d(\hat{\Delta})\ge5$, $c^{*}(\hat{\Delta})\le 0$ by Lemma 3.5(v).
Now assume that $l(\hat{\Delta})=ad\omega$ and put $d(\hat{\Delta})=k$. Since positive curvature crosses $(bd)$, $(ca)$ and $(ad)$ edges only it follows that if $l(\hat{\Delta})$ involves at least two occurrences of $(ad)^{\pm1}$ then there are at least four gaps and $c^{*}(\hat{\Delta})\le 0$ by Lemma 3.6, so assume otherwise. Moreover, $l(\hat{\Delta})=ad\omega$ forces $\hat{\Delta}$ to have at least two gaps and so $c^{*}(\hat{\Delta})\le
\pi \left( 2 - \frac{k}{2} \right)+ \frac{\pi}{2}+(k-3).\frac{\pi}{4}$ which implies that if $k\ge 7$ then $c^{*}(\hat{\Delta})\le 0$.
If $d(\hat{\Delta})=5$ then $l(\hat{\Delta})\in\{ad^{3}b^{-1},adb^{-1}ab^{-1},ad^{2}a^{-1}c^{\pm1},
ad^{2}b^{-1}c^{\pm1},ada^{-1}c^{\pm2},\\adb^{-1}c^{\pm2}\}$ which contradicts $|b|=3$ or $d^{\pm 1} b^{-1} a \neq 1$. If $d(\hat{\Delta})=6$ then checking shows that $\hat{\Delta}$ receives curvature across at most two edges and so $c^{*}(\hat{\Delta})\le 0$.  In conclusion in this case $\mathcal{P}$ is aspherical.

\textbf{(A9)} $\boldsymbol{|c|=3}$, $\bs{|d|=3}$, $\bs{a^{-1}b\ne1}$, $\bs{cab^{-1}\ne1}$, $\bs{c^{-1}ab^{-1}=1}$, $\bs{db^{-1}a=1}$, $\bs{d^{-1}b^{-1}a \neq 1}$.\newline
Since $H$ is cyclic $\mathcal{P}$ is aspherical by Lemma 3.2(ii).

\textbf{(A10)} $\boldsymbol{|c|>3}$, $\bs{|d|>3}$, $\bs{a^{-1}b\ne1}$, $\bs{cab^{-1}\ne1}$, $\bs{c^{-1}ab^{-1}=1}$, $\bs{db^{-1}a=1}$, $\bs{d^{-1}b^{-1}a \neq 1}$.\newline
Since $H$ is cyclic, $|d| < \infty$ and so $\mathcal{P}$ is aspherical by Lemma 3.2(ii).

It follows from (\textbf{A0})--(\textbf{A10}) above that either $\mathcal{P}$ is aspherical or modulo $T$-equivalence one of the conditions from Theorem 1.1(i), (iii) or Theorem 1.2(i), (ii), (iii) is satisfied and so Theorems 1.1 and 1.2 are proved for Case A.



\section{Proof of Case B} 

In this section we prove Theorems 1.1 and 1.2 for Case B, that is, we make the following assumption:
at least one of $c^2$, $d^2$, $a^{-1} b$ equals 1 in $H$.

\newpage
\begin{figure}
\begin{center}
\psfig{file=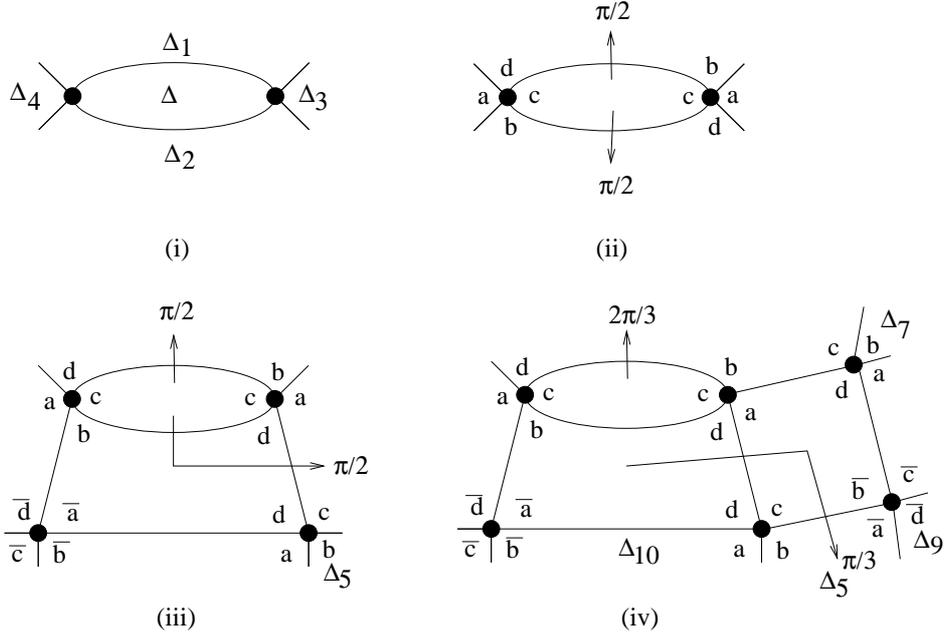}
\end{center}
\caption{the region $\Delta$ and curvature distribution for Cases (B1,2,3,5,6,8,10)}
\end{figure}

If $d(\Delta)=2$ then we will fix the names of the four neighbouring regions $\Delta_{i}$ ($1\le i\le 4$) of $\Delta$ as shown in Figure 5.1(i) and we use this notation throughout the section.

\textbf{Remark}\quad Recall that if $c^2=1$ then $c^{\pm 2}$ cannot be a proper sublabel.  This fact will be used often without explicit reference.

We treat each of the cases (\textbf{B1}) to (\textbf{B12}) in turn.

\textbf{(B1)} $\boldsymbol{|c|=2}$, $\bs{|d|>3}$, $\bs{a^{-1}b\ne1}$, $\bs{c^{\pm1}ab^{-1}\ne1}$, $\bs{d^{\pm1}b^{-1}a\ne 1}$.\newline
If $d(\Delta)=2$ then $\Delta$ is given by Figure 5.1(ii).
Observe that if $d(\Delta_{i})=4$ for $i\in\{1,2\}$ then $l(\Delta_{i})=\{bdda^{-1},bda^{-1}c^{\pm1},bdb^{-1}c^{\pm1}\}$. But $bdb^{-1}c^{\pm1}=1$ implies $|d|=|c|$, a contradiction. Observe further that at most one of $bd^{2}a^{-1},bda^{-1}c^{\pm1}$ equals 1 otherwise there is a contradiction to $|d| >3$.
This leaves the following cases:
\begin{enumerate}
\item[(i)]
$bd^{2}a^{-1}\ne1$, $bda^{-1}c^{\pm1}\ne1$;
\item[(ii)]
$bd^{2}a^{-1}=1$, $bda^{-1}c^{\pm1}\ne1$;
\item[(iii)]
$bda^{-1}c^{\pm1}=1$, $bd^{2}a^{-1}\ne1$.
\end{enumerate}

\begin{enumerate}
\item[(i)]
In this case $d(\Delta_{i})>4$ for $\Delta_{i}$ ($1 \leq i \leq 2$) of Figure 5.1(ii) so
add $\frac{1}{2}c(\Delta)=\frac{\pi}{2}$ to each of $c(\Delta_{i})$ ($1\le i\le 2$). Observe from Figure 5.1(ii) that $\Delta_{i}$ does not receive positive curvature from $\Delta_{j}$ for $j\in\{3,4\}$.
It follows that if $\hat{\Delta}$ receives positive curvature then it does so across at most half of its edges and so $d(\hat{\Delta})\ge 7$ implies that 

\newpage  
\begin{figure}
\begin{center}
\psfig{file=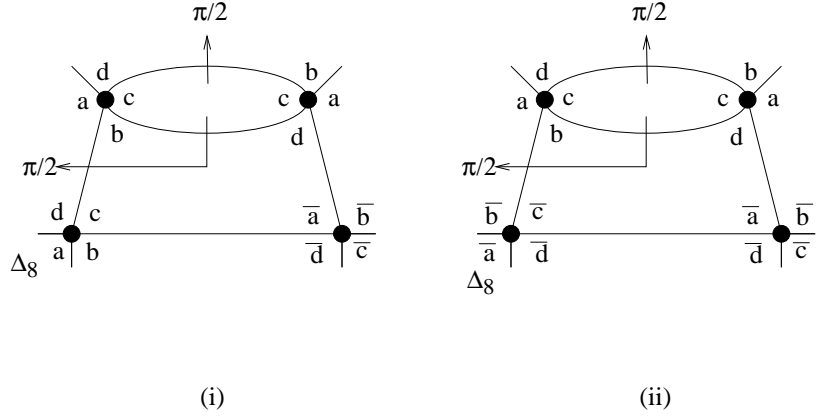}
\end{center}
\caption{curvature distribution for Case (B1) }
\end{figure}

\newpage
$c^{*}(\hat{\Delta})\le0$ by Lemma 3.5(vi). Checking (the LIST of Section 3) 
shows that 
if $d(\hat{\Delta})=5$ then $\hat{\Delta}$ receives positive curvature across at most one edge and so $c^{*}(\hat{\Delta})\le0$. Also if $d(\hat{\Delta})=6$ then checking shows $\hat{\Delta}$ receives positive curvature across at most two edges and so $c^{*}(\hat{\Delta})\le0$.

\item[(ii)]
Suppose that $l(\Delta)=bd^{2}a^{-1}=1$, $bda^{-1}c^{\pm1}\ne1$. If $d(\Delta_{1})>4$ and $d(\Delta_{2})>4$ then add  $\frac{1}{2}c(\Delta) =\frac{\pi}{2}$ to $c(\Delta_{1})$ and $c(\Delta_{2})$ as in Figure 5.1(ii). If say $d(\Delta_{2})=4$ as in Figure 5.1(iii) then $l(\Delta_{2})=bdda^{-1}$ which forces $l(\Delta_{3})=ca\omega$.
First assume that $cadb^{-1} \neq 1$.  Then $d(\Delta_3) > 3$ and so
add $\frac{\pi}{2}$ to $c(\Delta_{3})$ via $\Delta_2$. If $d(\Delta_{1})=4$ then add $\frac{\pi}{2}$ to $c(\Delta_4)$ via $\Delta_1$ in a similar way.
Observe that $\Delta_1$ does not receive positive curvature from $\Delta_{3}$ or $\Delta_{4}$ in Figure 5.1(ii); and $\Delta_3$ does not receive positive curvature from $\Delta_{1}$ or $\Delta_{5}$ in Figure 5.1(iii). It follows that if $\hat{\Delta}$ receives positive curvature then it does so across at most half of its edges and so $d(\hat{\Delta})\ge 7$ implies that $c^{*}(\hat{\Delta})\le0$ by Lemma 3.5(vi).
It remains to study $5 \le d(\hat{\Delta})\le6$.
Checking shows that if $d(\hat{\Delta})=5$ then either the label contradicts $|c|\ne1$ or $\hat{\Delta}$ receives positive curvature across at most one edge and so $c^{*}(\hat{\Delta})\le0$. Also if $d(\hat{\Delta})=6$ then
checking shows that $\hat{\Delta}$ receives positive curvature across at most two edges and so $c^{*}(\hat{\Delta})\le0$.
Now assume that $cadb^{-1}=1$, in which case $c=d^3$, $b=d^{-2}$ and $|d|=6$.  The distribution of curvature is exactly the same except when $d(\Delta_3)=4$ in Figure 5.1(iii).  In this case add $\frac{2}{3} c(\Delta)=\frac{2\pi}{3}$ to $c(\Delta_1)$ and $\frac{1}{3}c(\Delta)=\frac{\pi}{3}$ to $c(\Delta_5)$ via $\Delta_3$ as shown in Figure 5.1(iv).  Together with the observations above (which still hold) we also have that $\Delta_1$ does not receive positive curvature from $\Delta_3$, $\Delta_4$ or $\Delta_7$ and that $\Delta_5$ does not receive positive curvature from $\Delta_9$ or $\Delta_{10}$ in Figure 5.1(iv).  An argument similar to those given in the proof of Lemma 3.5 now shows that if $d(\hat{\Delta}) \geq 8$ then $c^{\ast}(\hat{\Delta}) \leq 0$; and that if $l(\hat{\Delta})$ does not involve $(cbd)^{\pm 1}$ then $c^{\ast}(\hat{\Delta}) \leq 0$ for $d(\hat{\Delta}) \geq 7$ by Lemma 3.5(vi).  The conditions on $b$, $c$ and $d$ imply that if $2 < d(\hat{\Delta}) < 6$ then $l(\hat{\Delta}) \in \{ d^2 a^{-1} b, db^{-1} c^{\pm 1} a \}$ so if $l(\hat{\Delta})$ does not involve $(cbd)^{\pm 1}$ it remains to consider $d(\hat{\Delta})=6$.  But checking shows that $l(\hat{\Delta})$ will then either involve at most two non-adjacent occurrences of $(bd)^{\pm 1}$, $(ca)^{\pm 1}$ or 

\begin{figure}  
\begin{center}
\psfig{file=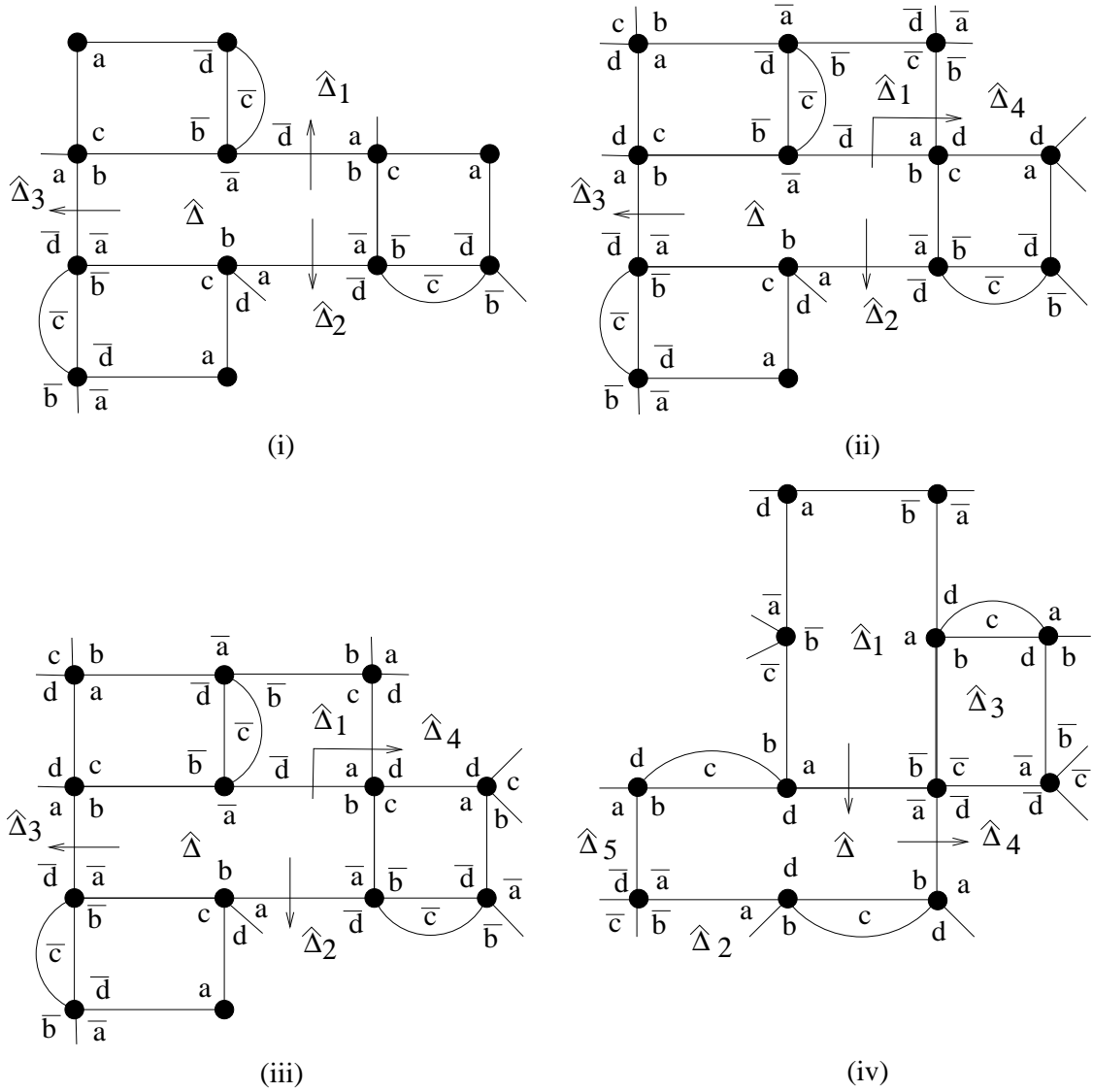}
\end{center}
\caption{curvature distribution for Case (B1)}
\end{figure}

\noindent $(ba^{-1})^{\pm 1}$ and so
$c^{\ast}(\hat{\Delta}) \leq 0$ or
$l(\hat{\Delta})=(ab^{-1})^3$ in which case $c^{\ast}(\hat{\Delta}) \leq c(\hat{\Delta}) + 3 ( \frac{\pi}{3} )=0$.
Finally if $l(\hat{\Delta})=cbd \omega$ and $d(\hat{\Delta}) \leq 7$ then
$l(\hat{\Delta}) \in \{ cbda^{-1} ba^{-1}, cbdb^{-1} ab^{-1}, cbd^3 b^{-1} \}$ and $c^{\ast}(\hat{\Delta}) \leq c(\hat{\Delta}) + \frac{2 \pi}{3} + \frac{\pi}{3} = 0$.

\item[(iii)]
Now suppose that $l(\Delta)=bda^{-1}c^{\pm1}=1$, $bd^{2}a^{-1}\ne1$. First assume that $|b|\ge3$.
Add $\frac{1}{2}c(\Delta)=\frac{\pi}{2}$ to $c(\Delta_{1})$ and $c(\Delta_{2})$ as in Figure 5.1(ii). If say $d(\Delta_{2})=4$ then $l(\Delta_{2})=bda^{-1}c^{\pm1}$. First let $l(\Delta_{2})=bda^{-1}c$ as in Figure 5.2(i). This forces $l(\Delta_{4})=ad\omega$ and so
$d(\Delta_{4})=4$ forces $l(\Delta_4)=addb^{-1}$.  But if $b=d^2$ then $c=d^3$ and there is a sphere by Lemma 3.2(iv), so it can be assumed that $d(\Delta_4) > 4$.
So add $\frac{\pi}{2}$ to $c(\Delta_4)$ via $\Delta_2$ as shown.
Suppose now that $l(\Delta_{2})=bda^{-1}c^{-1}$ as in Figure 5.2(ii). This forces $l(\Delta_{4})=ab^{-1}\omega$ and so $d(\Delta_{4})>4$, otherwise there is a contradiction to $|b|\ge3$. So add $\frac{\pi}{2}$ to $c(\Delta_4)$ via $\Delta_2$ as shown.
Similarly add $\frac{1}{2} c(\Delta)=\frac{\pi}{2}$ to $c(\Delta_3)$ if $d(\Delta_1)=4$.
Observe that $\Delta_{1}$ does not receive positive curvature from $\Delta_{3}$ or $\Delta_{4}$ in Figure 5.1(ii); $\Delta_{2}$ does not receive positive curvature from $\Delta_{3}$ or $\Delta_{4}$ in Figure 5.1(ii);
and $\Delta_{4}$ does not receive positive curvature from $\Delta_{1}$ or $\Delta_{8}$ in Figures 5.2(i), (ii). It follows that if $\hat{\Delta}$ receives positive curvature then it does so across at most half of its edges and so $d(\hat{\Delta})\ge 7$ implies that $c^{*}(\hat{\Delta})\le0$ by Lemma 3.5(vi).
It remains to study $5 \le d(\hat{\Delta})\le6$.
\end{enumerate}

If $|b|>3$ then checking shows that if $d(\hat{\Delta})=5$ then either the label contradicts 
one of the (\textbf{B1}) assumptions or $\hat{\Delta}$ receives positive curvature across at most one edge and so $c^{*}(\hat{\Delta})\le0$. Checking shows that if  $d(\hat{\Delta})=6$ then $\hat{\Delta}$ receives 
positive curvature across at most two edges
or $l(\hat{\Delta})=(ab^{-1})^3$ contradicting $|b| > 3$, and so $c^{*}(\hat{\Delta})\le0$.

Let $|b|=3$.  If $H$ is cyclic then $bdc=1$ implies $b=d^2$ and we obtain a sphere as before, so assume otherwise.
If $|b|=3$ and $|d|\in\{4,5\}$ then we obtain a sphere by Lemma 3.1(b)(iv). So let $|b|=3$, $|d|\ge6$.
Distribute curvature from $\Delta$ as above and as shown in Figure 5.2.
Checking shows that if $d(\hat{\Delta})=5$ then either the label contradicts $|b|=3$ or $|d|\ge6$ or $\hat{\Delta}$ receives positive curvature across at most one edge and so $c^{*}(\hat{\Delta})\le0$. Checking shows that if $d(\hat{\Delta})=6$ then $\hat{\Delta}$ receives positive curvature across at most two edges and so $c^{*}(\hat{\Delta})\le0$ except when $l(\hat{\Delta})=ba^{-1}ba^{-1}ba^{-1}$.
This case is shown in Figure 5.3(i) where $c^{*}(\hat{\Delta})=\frac{\pi}{2}$ and so add $\frac{1}{3} c^{*}(\hat{\Delta})=\frac{\pi}{6}$ to $c(\hat{\Delta}_{i})$ for $i\in\{1,2,3\}$ across the edge $ad^{-1}$. If $d(\hat{\Delta}_{i})=4$ then $l(\hat{\Delta}_{i})\in\{ad^{-1}b^{-1}c^{-1},ad^{-1}b^{-1}c\}$. Suppose that $l(\hat{\Delta}_{1})=ad^{-1}b^{-1}c^{-1}$ as in Figure 5.3(ii).  Then $l(\hat{\Delta}_{4})=d^{2}b^{-1}\omega$ and $d(\hat{\Delta}_{4})>4$ otherwise there is a contradiction to $H$ non-cyclic. So add $\frac{\pi}{6}$ to $c(\hat{\Delta}_{4})$ across the edge $db^{-1}$. If $l(\hat{\Delta}_{1})=ad^{-1}b^{-1}c$ as in Figure 5.3(iii) then $l(\hat{\Delta}_{4})=d^{3}\omega$ and $d(\hat{\Delta}_{4})>4$ otherwise there is a contradiction to $|d|\ge6$. So add $\frac{\pi}{6}$ to $c(\hat{\Delta}_{4})$ across the edge $d^{2}$.
Observe that if $\hat{\Delta}$ receives positive curvature then it receives $\frac{\pi}{2}$ across the edges $ab^{-1}$, $ad$ or $bd$; and receives $\frac{\pi}{6}$ across the edges $ad^{-1}$, $db^{-1}$ or $dd$.
Thus there is a gap preceding $c^{\pm1}$, $b$, $a$ and a gap after $c^{\pm1}$, $b^{-1}$, $a^{-1}$ and there is a two-thirds gap across the edges $ad^{-1}$, $db^{-1}$ and $dd$. Also there is always a gap between two $d$'s (other than when the subword is $d^{\pm 2}$).  Now since $c^2$ cannot be a proper sublabel it follows that if there are at least two occurrences of $c^{\pm1}$ then we obtain four gaps. Suppose now that there is at most one occurrence of $c^{\pm1}$. If there is exactly one occurrence of $c$ and either no occurrences of $b$ or no occurrences of $d$ then $H$ is cyclic, a contradiction; and if there 

\begin{figure}
\begin{center}
\psfig{file=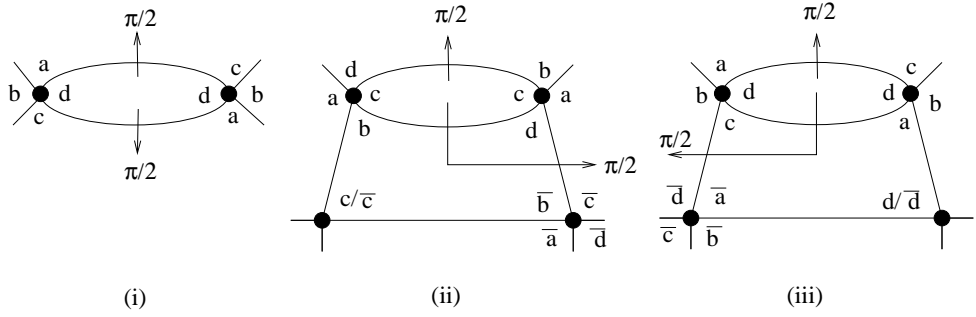}
\end{center}
\caption{curvature distribution for Case (B2)}
\end{figure}

\noindent are no occurrences of $c$ and exactly one occurrence of $d$ or of $b$ then again $H$ is cyclic, a contradiction. So assume otherwise. It follows that $l(\hat{\Delta})$ contains at least four gaps or 
$l(\hat{\Delta})\in\{c^{\pm1}ad^{\pm1}a^{-1}ba^{-1},c^{\pm1}bd^{\pm1}a^{-1}ba^{-1},c^{\pm1}ad^{\pm1}b^{-1}ab^{-1},
c^{\pm1}bd^{\pm1}b^{-1}ab^{-1},\\d(a^{-1}b)^{\pm1}d^{\pm1}(a^{-1}b)^{\pm1}\}$. But since $c=bd$ each of these labels contradicts $H$ non-cyclic, $|d|\ge6$ or $|b|=3$ except when $l(\hat{\Delta})=da^{-1}bda^{-1}b$. In this case if $c^{*}(\hat{\Delta})>0$ then it can be assumed without loss of generality that $\hat{\Delta}$ is given by Figure 5.3(iv) and $\hat{\Delta}$ receives $\frac{1}{3}c^{*}(\hat{\Delta}_{1})=\frac{\pi}{6}$ from $c(\hat{\Delta}_{1})$. This implies that $l(\hat{\Delta}_{3})=a^{-1}c^{-1}bd$ and $l(\hat{\Delta}_{4})=ad^{-1}d^{-1}\omega$. So add $\frac{\pi}{6}$ from $c(\hat{\Delta})$ to $c(\hat{\Delta}_{4})$. Since this $\frac{\pi}{6}$ is across an $ad^{-1}$ edge and since $l(\hat{\Delta}_{4})=ad^{-2}\omega$ it follows from the above that $c^{*}(\hat{\Delta}_{4})\le 0$. If $l(\hat{\Delta}_{2})=b^{-1}ab^{-1}ab^{-1}$ in Figure 5.3(iv) then a similar argument applies to $c(\hat{\Delta}_{5})$.

Finally let $|b|=2$.
In particular, $H$ is non-cyclic for otherwise $d^2=1$, a contradiction.
If $|d|<\infty$ then we obtain a sphere by Lemma 3.1(b)(i) and if $|d| = \infty$ then $\mathcal{P}$ is aspherical by Lemma 3.4(ii).

In conclusion $\mathcal{P}$ is aspherical in this case except when $H$ is non-cyclic, $bda^{-1}c^{\pm 1}=1$ and either $|b|=3$, $|d| \in \{4,5\}$ or $|b|=2$, $|d| < \infty$; or when $H$ is cyclic, $b=d^2$, $c=d^3$ and $|d|=6$.

\textbf{(B2)} $\boldsymbol{|c|=2}$, $\bs{|d|=2,~ a^{-1}b\ne1, ~c^{\pm1}ab^{-1}\ne1,~ d^{\pm1}b^{-1}a\ne1}$.\newline
If $H$ is cyclic then $c=d$.  Therefore $bd^{-1}b^{-1}c=cada^{-1}=1$ and there is a sphere by Lemma 3.2(i).
So suppose from now on that $H$ is non-cyclic.
If $d(\Delta)=2$ then $\Delta$ is given by Figures 5.1(ii) and 5.4(i).
Moreover, if $d(\Delta_{i})=4$ for $i\in\{1,2\}$ then
$l(\Delta_{i})\in\{bd^{2}a^{-1},bda^{-1}c^{\pm1}, bdb^{-1}c^{\pm1}, cad^{\pm1}a^{-1},cad^{\pm1}b^{-1}\}$. But $bdda^{-1}=1$ implies that $b=1$ which is a contradiction. This leaves the following cases:
\begin{enumerate}
\item[(i)]
$bda^{-1}c^{\pm1}\ne1$, $bdb^{-1}c^{\pm1}\ne1$, $cad^{\pm1}a^{-1}\ne1$;
\item[(ii)]
$bda^{-1}c^{\pm1}=1$, $bdb^{-1}c^{\pm1}\ne1$, $cad^{\pm1}a^{-1}\ne1$;
\item[(iii)]
$bdb^{-1}c^{\pm1}=1$, $bda^{-1}c^{\pm1}\ne1$, $cad^{\pm1}a^{-1}\ne1$;
\item[(iv)]
$cad^{\pm1}a^{-1}=1$, $bda^{-1}c^{\pm1}\ne1$, $bdb^{-1}c^{\pm1}\ne1$;
\item[(v)]
$bdb^{-1}c^{\pm1}=1$, $cad^{\pm1}a^{-1}=1$, $bda^{-1}c^{\pm1}\ne1$.
\end{enumerate}
(Note that any other combination implies $b=1$.)
\begin{enumerate}
\item[(i)]
In this case $d(\Delta_{1})>4$ and $d(\Delta_{2})>4$ in Figures 5.1(ii) and 5.4(i) so add $\frac{1}{2}c(\Delta) =\frac{\pi}{2}$ to each of $c(\Delta_{1})$ and $c(\Delta_{2})$.
Observe that $\Delta_{1}$ and $\Delta_{2}$ do not receive positive curvature from $\Delta_{3}$ or $\Delta_{4}$ in Figures 5.1(ii) and 5.4(i). It follows that if $\hat{\Delta}$ receives positive curvature then it does so across at most half of its edges and so $d(\hat{\Delta})\ge7$ implies that $c^{*}(\hat{\Delta})\le0$ by Lemma 3.5(vi).
It remains to study $5 \le d(\hat{\Delta})\le6$.
Checking shows that if $d(\hat{\Delta})=5$ then either the label contradicts $c^{\pm1} ab^{-1} \neq 1$ or $\hat{\Delta}$ receives positive curvature across at most one edge and so $c^{*}(\hat{\Delta})\le0$. Also if $d(\hat{\Delta})=6$ then $\hat{\Delta}$ receives positive curvature across at most two edges and so $c^{*}(\hat{\Delta})\le0$.

\item[(ii)]
In this case the labels $bda^{-1}c^{\pm1}$ can occur. If $|b|< \infty$ then we obtain spheres by Lemma 3.1(b)(ii) and if $|b|=\infty$ then $\mathcal{P}$ is aspherical by Lemma 3.4(iii).

\item[(iii)]
In this case the labels $bdb^{-1}c^{\pm1}$ can occur and
$H$ is non-Abelian. If $d(\Delta_{1})>4$ and $d(\Delta_{2})>4$ as in Figures 5.1(ii) and 5.4(i) then add $\frac{1}{2}c(\Delta) =\frac{\pi}{2}$ to $c(\Delta_{1})$ and $c(\Delta_{2})$.
If say $d(\Delta_{2})=4$ as in Figure 5.4(ii) then $l(\Delta_{3})=c^{-1}a\omega$ and so $d(\Delta_{3})>4$, otherwise there is a contradiction to $c\ne d$ or $b\ne1$. So add $\frac{\pi}{2}$ to $c(\Delta_{3})$ via $\Delta_2$ as shown.
If $d(\Delta_1)=4$ then similarly add $\frac{1}{2}c(\Delta)=\frac{\pi}{2}$ to $c(\Delta_4)$ via $\Delta_1$.
We see from Figures 5.1(ii) and 5.4(i)--(ii) that if $\hat{\Delta}$ receives positive curvature then it does so across the edges $b d, ca$ or $ c^{-1}a$. It follows that the only word of length 3 that contains no gaps is $a^{-1}c^{\pm1}a$. Suppose that $l(\hat{\Delta})=\omega_{1} a^{-1}c^{\pm1}a\omega_{2}\omega$ where $\omega_{1}$ and $\omega_{2}$ have length 2 and $\omega$ has length at least 0. Then there is a gap preceding $a^{-1}$ in $\omega_{1}a^{-1}$ and after $a$ in $a\omega_{2}$. Moreover, if $\omega_{1}$ does not contain a gap then $\omega_{1}=bd$ and if $\omega_{2}$ does not contain a gap then $\omega_{2}=d^{-1}b^{-1}$. It follows that if $l(\omega)>0$ then $l(\hat{\Delta})$ contain at least 4 gaps and $c^{*}(\hat{\Delta})\le0$. If $l(\omega)=0$ and $l(\hat{\Delta})=\omega_{1}a^{-1}c^{\pm1}ad^{-1}b^{-1}$ then $l(\hat{\Delta})\in\{ c^{\pm1}ba^{-1}c^{\pm1}ad^{-1}b^{-1}, ad^{\pm1}a^{-1}c^{\pm1}ad^{-1}b^{-1}\}$ and there are 4 gaps; if $l(\omega)=0$ and $l(\hat{\Delta})=bda^{-1}c^{\pm1}a\omega_{2}$ then $l(\hat{\Delta})\in\{bda^{-1}c^{\pm1}ad^{\pm1}a^{-1},bda^{-1}c^{\pm1}ab^{-1}c^{\pm1}\}$ and again there are 4 gaps.
If now $l(\hat{\Delta})=a^{-1} c^{\pm 1} aw$ and $d(\hat{\Delta}) \leq 6$ then either $l(\hat{\Delta}) = a^{-1} c^{\pm 1} ad^{\pm 2}$, a contradiction, or there are 4 gaps. Now suppose that $l(\hat{\Delta})$ does not contain the subword $a^{-1}c^{\pm1}a$. Then $\hat{\Delta}$ receives positive curvature across at most half of its edges and so $d(\hat{\Delta})\ge7$ implies $c^{*}(\hat{\Delta})\le0$ by Lemma 3.5(vi). So let $d(\hat{\Delta})\le6$ and $l(\hat{\Delta})\in\{ca\omega,bd\omega,c^{-1}a\omega\}$. If $d(\hat{\Delta})<6$ then checking shows that there is a contradiction to $H$ non-Abelian; and if $d(\hat{\Delta})=6$ then checking shows that $\hat{\Delta}$ receives positive curvature across at most two edges and so $c^{*}(\hat{\Delta})\le0$.

\item[(iv)]
In this case the labels $cad^{\pm1}a^{-1}$ can occur and $H$ is non-Abelian.
If $d(\Delta_{1})>4$ and $d(\Delta_{2})>4$ as in Figures 5.1(ii) and 5.4(i) then add $\frac{1}{2}c(\Delta) =\frac{\pi}{2}$ to $c(\Delta_{1})$ and $c(\Delta_{2})$.
If say $d(\Delta_{2})=4$ then $l(\Delta_{2})=cad^{\pm1}a^{-1}$ which forces $l(\Delta_{4})=bd^{-1}\omega$ and so $d(\Delta_{3})>4$, otherwise there is a contradiction to $b\ne1$ or $bd^{-1}b^{-1}c^{\pm1}\ne1$. So add $\frac{\pi}{2}$ to $c(\Delta_{4})$ via $\Delta_2$
as shown in Figure 5.4(iii).  If $d(\Delta_1)=4$ then similarly add $\frac{1}{2}c(\Delta)= \frac{\pi}{2}$ to $c(\Delta_3)$ via $\Delta_1$.
We see from Figures 5.1(ii), 5.4(i) and 5.4(iii) that if $\hat{\Delta}$ receives positive curvature then it does so across the edges $bd$, $bd^{-1}$ or $ca$. It follows that the only word of length 3 that contains no gaps is $bd^{\pm1}b^{-1}$. Suppose that $l(\hat{\Delta})=\omega_{1} bd^{\pm1}b^{-1}\omega_{2}\omega$ where $\omega_{1}$ and $\omega_{2}$ have length 2 and $\omega$ has length at least 0. Then there is a gap preceding $b$ in $\omega_{1}b$ and after $b^{-1}$ in $b^{-1}\omega_{2}$. Moreover, if $\omega_{1}$ does not contain a gap then $\omega_{1}=a^{-1}c^{-1}$ and if $\omega_{2}$ does not contain a gap then $\omega_{2}=ca$. It follows that if $l(\omega)>0$ then $l(\hat{\Delta})$ contain at least 4 gaps and $c^{*}(\hat{\Delta})\le0$. If $l(\omega)=0$ and $l(\hat{\Delta})=\omega_{1}bd^{\pm1}b^{-1}ca$ then $l(\hat{\Delta})\in\{ d^{\pm1}a^{-1}bdb^{-1}ca, bc^{\pm1}bdb^{-1}ca\}$ and there are 4 gaps; if $l(\omega)=0$ and $l(\hat{\Delta})=a^{-1}c^{-1}bd^{\pm1}b^{-1}\omega_{2}$ then $l(\hat{\Delta})\in\{a^{-1}c^{-1}bdb^{-1}ad^{\pm1},a^{-1}c^{-1}bdb^{-1}c^{\pm1}b\}$ and again there are 4 gaps.
If $l(\hat{\Delta})=bd^{\pm1} b^{-1}\omega$ and $d(\hat{\Delta})\le6$ then there are 4 gaps. Now suppose that $l(\hat{\Delta})$ does not contain the subword $bd^{\pm1}b^{-1}$. Then $\hat{\Delta}$ receives positive curvature across at most half of its edges and so $d(\hat{\Delta})\ge7$ implies $c^{*}(\hat{\Delta})\le0$ by Lemma 3.5(vi). So let $d(\hat{\Delta})\le6$ and $l(\hat{\Delta})\in\{bd\omega,bd^{-1}\omega,ca\omega\}$. If $d(\hat{\Delta})<6$ then checking shows that there is a contradiction to $H$ non-Abelian; and if $d(\hat{\Delta})=6$ then checking shows that $\hat{\Delta}$ receives positive curvature across at most two edges and so $c^{*}(\hat{\Delta})\le0$.

\item[(v)]
In this case the labels $bdb^{-1}c^{\pm1}$ and $cad^{\pm1}a^{-1}$ can occur.
If $|b| = \infty$ then $\mathcal{P}$ is aspherical by Lemma 3.4(iii); and if $|b| < \infty$ then there is a sphere by Lemma 3.2(i).
\end{enumerate}

In conclusion $\mathcal{P}$ fails to be aspherical in this case either when $H$ is cyclic or when $H$ is non-cyclic and when $bda^{-1} c^{-1} =1$, $|b| < \infty$ or when
$bdb^{-1} c^{-1}= a^{-1} cad =1$, $|b| < \infty$.

\textbf{(B3)} $\boldsymbol{|c|=2}$, $\bs{|d|=3,~ a^{-1}b\ne1,~ c^{\pm1}ab^{-1}\ne1,~ d^{\pm1}b^{-1}a\ne1}$.\newline
If $d(\Delta)=2$ then $\Delta$ is given by
Figure 5.1(ii). If $d(\Delta)=3$ then $\Delta$ is given by Figure 4.3(i).
Moreover, if $d(\Delta_{i})=4$ and $l(\Delta_i)=bdw$ or $caw$ then\newline $l(\Delta_{i})\in\{bd^{2}a^{-1},bda^{-1}c^{\pm1}, bdb^{-1}c^{\pm1}, cad^{\pm1}a^{-1},cad^{\pm1}b^{-1}\}$. But each of $bd^2 a^{-1} =1$,\newline
$bdb^{-1}c^{\pm1}=1$ and $cad^{\pm1}a^{-1}=1$ implies a contradiction to one of the (\textbf{B3}) assumptions. Thus we have the following cases:
\begin{enumerate}
\item[(i)]
$bda^{-1}c^{\pm 1}\ne1$, $cadb^{-1}\ne1$;
\item[(ii)]
$bda^{-1}c^{\pm1}=1$, $cadb^{-1}\ne1$;
\item[(iii)]
$cadb^{-1}=1$, $bda^{-1}c^{\pm 1}\ne1$.
\end{enumerate}

\begin{enumerate}
\item[(i)]
In this case $d(\Delta_{1})>4$ and $d(\Delta_{2})>4$ in Figure 5.1(ii), so add $\frac{1}{2}c(\Delta) =\frac{\pi}{2}$ to each of $c(\Delta_{1})$ and $c(\Delta_{2})$. In Figure 4.3(i) $d(\Delta_{1})>4$, $d(\Delta_{3})>4$ and $d(\Delta_{5})>4$ so add $\frac{1}{3}c(\Delta) =\frac{\pi}{6}$ to each of $c(\Delta_{1})$, $c(\Delta_{3})$ and $c(\Delta_{5})$.
Observe that $\Delta_{1}$ and $\Delta_{2}$ do not receive positive curvature from $\Delta_{3}$ or $\Delta_{4}$ in Figure 5.1(ii). Also $\Delta_{1}$, $\Delta_{3}$ and $\Delta_{5}$ do not receive positive curvature from $\Delta_{m}$ for $m\in\{2,4,6\}$ in Figure 4.3(i). It follows that if $\hat{\Delta}$ receives positive curvature then it does so across at most half of its edges and so $d(\hat{\Delta})\ge7$ implies that $c^{*}(\hat{\Delta})\le0$ by Lemma 3.5(vi). It remains to study $5 \le d(\hat{\Delta})\le6$.
Checking shows that if $d(\hat{\Delta})=5$ then either the label contradicts 

\begin{figure}
\begin{center}
\psfig{file=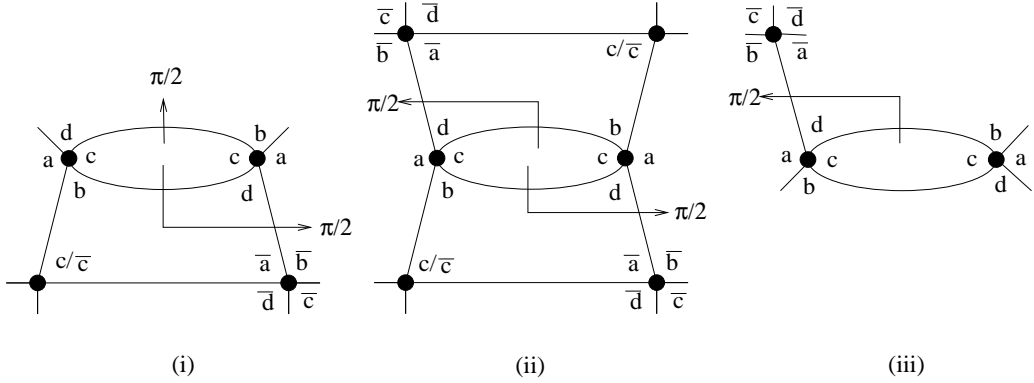}
\end{center}
\caption{curvature distribution for Case (B3)}
\end{figure}

\noindent $cab^{-1}\ne1$ or $\hat{\Delta}$ receives positive curvature across at most one edge and so $c^{*}(\hat{\Delta})\le0$. Also if $d(\hat{\Delta})=6$ then $\hat{\Delta}$ receives positive curvature across at most two edges 
and so $c^{*}(\hat{\Delta})\le0$.

\item[(ii)]
In this case the labels $bda^{-1}c^{\pm1}$ can occur.
If $H$ is cyclic then $d=b^2$, $c=b^3$ and there is a sphere by Lemma 3.2(v), so assume that $H$ is non-cyclic.
If $|b|\in\{2,3,4,5\}$ then we obtain spheres by Lemma 3.1(b)(i), (v),
so assume that $|b|\ge6$. In Figure 5.1(ii) if $d(\Delta_{1})>4$ and $d(\Delta_{2})>4$ then add $\frac{1}{2}c(\Delta) =\frac{\pi}{2}$ to $c(\Delta_{1})$ and $c(\Delta_{2})$ as shown.
If say $d(\Delta_{1})>4$ and $d(\Delta_{2})=4$ as in Figure 5.5(i) then add $\frac{\pi}{2}$ to $c(\Delta_{1})$. This implies that $l(\Delta_{2})=bda^{-1}c^{\pm1}$ as shown. This forces $l(\Delta_{3})=b^{-1}a\omega$ and so $d(\Delta_{3})>4$, otherwise there is a contradiction to $|b|\ge6$ so add $\frac{\pi}{2}$ to $c(\Delta_{3})$ via $\Delta_2$ as shown.
If $d(\Delta_{1})=d(\Delta_{2})=4$ then add $\frac{\pi}{2}$ to $c(\Delta_{j})$ for $j\in\{3,4\}$ as shown in Figure 5.5(ii).
The one exception to the above is when $l(\Delta_{1})=bda^{-1}\omega$ and $d(\Delta_1) > 4$. Then $d(\Delta_{4})>4$ and in this situation add the $\frac{\pi}{2}$ from $c(\Delta)$ to $c(\Delta_{4})$ via $\Delta_1$ as shown in Figure 5.5(iii). The same applies to $\Delta_{2}$.
In Figure 4.3(i) if $d(\Delta_{1})>4$, $d(\Delta_{3})>4$ and $d(\Delta_{5})>4$ then add $\frac{1}{3}c(\Delta) =\frac{\pi}{6}$ to each of $c(\Delta_{1})$, $c(\Delta_{2})$ and $c(\Delta_{5})$. If say $d(\Delta_{1})=4$, $d(\Delta_{3})>4$ and $d(\Delta_{5})>4$ then $l(\Delta_{2})=ba^{-1}\omega$ and $d(\Delta_{2})>4$ otherwise there is a contradiction to $|b|\ge6$ so add $\frac{\pi}{6}$ to $c(\Delta_{2})$, $c(\Delta_{3})$ and $c(\Delta_{5})$ as shown in Figure 5.6(i). Now suppose that $d(\Delta_{1})=4$ and $d(\Delta_{3})=4$. This implies that $l(\Delta_{2})=l(\Delta_{4})=ba^{-1}\omega$ as shown in Figure 5.6(ii). So add $\frac{\pi}{6}$ to $c(\Delta_{2})$, $c(\Delta_{4})$ and $c(\Delta_{5})$. If $d(\Delta_{1})=d(\Delta_{3})=d(\Delta_{5})=4$ then similarly add $\frac{\pi}{6}$ to $c(\Delta_{m})$ for $m\in\{2,4,6\}$.

We now see that if $\hat{\Delta}$ receives positive curvature then it receives at most $\frac{\pi}{2}$ across $(bd)^{\pm1}$ and $(b^{-1}a)^{\pm1}$; and it receives at most $\frac{\pi}{6}$ across $(ca)^{\pm1}$ and $(ab^{-1})^{\pm1}$. Thus there is always a gap immediately preceding $c$ and $d^{-1}$; and there is a gap immediately after $c^{-1}$ and $d$. This implies that if there are at least four occurrences of $c^{\pm1}$ or $d^{\pm1}$ then $l(\hat{\Delta})$ contains at least four gaps and so $c^{*}(\hat{\Delta})\le 0$ by Lemma 3.6.
Suppose that there are at most three occurrences of $c^{\pm1}$ or $d^{\pm1}$ in $l(\hat{\Delta})$.
Observe that in addition to the four gaps mentioned above the following sublabels yield gaps: $(cb)^{\pm1}$ and $(ad)^{\pm1}$ each yields a gap;
$(bda^{-1})^{\pm 1}$ yields two gaps (see Figure 5.5(iii)); and $(ca)^{\pm1}$ and $(ab^{-1})^{\pm1}$ 

\begin{figure}
\begin{center}
\psfig{file=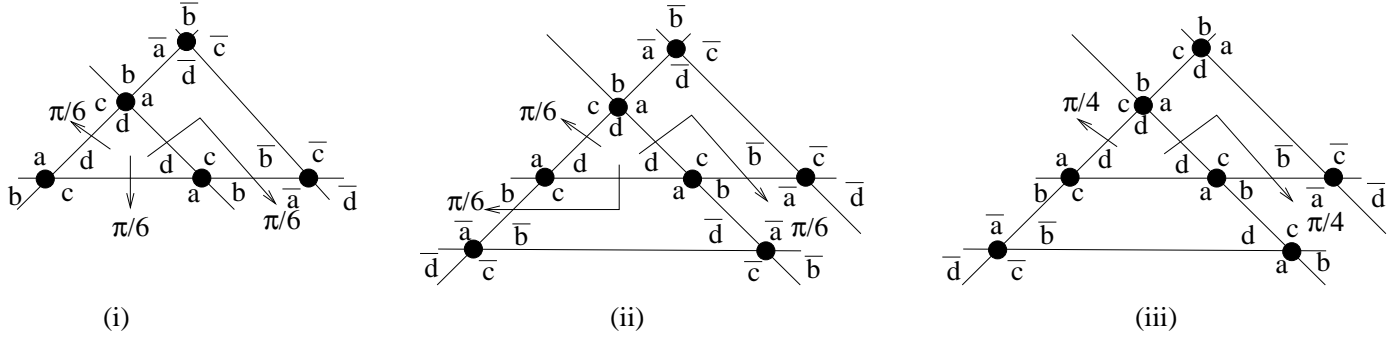}
\end{center}
\caption{curvature distribution for Case (B3)}
\end{figure}

\noindent each yields the equivalent of a two-thirds gap. If $l(\hat{\Delta})= (b^{-1}a)^{\pm n}$ where $n\ge1$ then $l(\hat{\Delta})$ obtains at least four gaps since $|b|\ge6$. If $l(\hat{\Delta})\in\{d^{\pm 1}(b^{-1}a)^{\pm 
n},(ab^{-1})^{\pm n}c^{\pm1}\}$ then $H$ is cyclic so it can be assumed that $l(\hat{\Delta})$ involves either two or three occurrences of $c^{\pm1}$ or $d^{\pm1}$.
It follows that if there are three occurrences then $c^{\ast}(\hat{\Delta}) \leq 0$; or if exactly two occurrences then either
$c^{\ast}(\hat{\Delta}) \leq 0$ or\newline
$l(\hat{\Delta})\in\{d^{\pm1}b^{-1}ada^{-1}b,d^{\pm1}b^{-1}ada^{-1}ba^{-1}b,d^{\pm1}b^{-1}ab^{-1}ada^{-1}b,
cba^{-1}bd^{\pm1}b^{-1},cbd^{\pm1}b^{-1}ab^{-1}\}$. But each of these labels forces $H$ cyclic or $|b|<6$ or a (\textbf{B3}) contradiction, therefore $c^{*}(\hat{\Delta})\le 0$ by Lemma 3.6.

\item[(iii)]
In this case the label $cadb^{-1}$ can occur.  First assume that $H$ is non-cyclic. In Figure 5.1(ii) $d(\Delta_{1})>4$ and $d(\Delta_{2})>4$ otherwise there is a contradiction to $|c|=2$ or $|d|=3$, so add $\frac{1}{2}c(\Delta)=\frac{\pi}{2}$ to $c(\Delta_{1})$ and $c(\Delta_{2})$. In Figure 4.3(i) if $d(\Delta_{1})>4$, $d(\Delta_{3})>4$ and $d(\Delta_{5})>4$ add $\frac{1}{3}c(\Delta) =\frac{\pi}{6}$ to each of $c(\Delta_{1})$, $c(\Delta_{2})$ and $c(\Delta_{5})$. If say $d(\Delta_{1})=4$ only then add $\frac{\pi}{4}$ to $c(\Delta_{3})$ and $c(\Delta_{5})$. Now suppose that $d(\Delta_{1})=d(\Delta_{3})=4$. This implies that their label is $cadb^{-1}$ which forces $l(\Delta_{2})=cba^{-1}\omega$ as shown in Figure 5.6(iii). So add $\frac{\pi}{4}$ to $c(\Delta_{2})$ via $\Delta_1$ and to $c(\Delta_{5})$. Finally if $d(\Delta_{1})=d(\Delta_{3})=d(\Delta_{5})=4$ then in a similar way add $\frac{\pi}{6}$ to $c(\Delta_{m})$ for $m\in\{2,4,6\}$.
Observe that $\Delta_{1}$ does not receive positive curvature from $\Delta_{3}$ and $\Delta_{2}$ does not receive positive curvature from $\Delta_{4}$ in Figure 5.1(ii). In Figure 4.3(i) $\Delta_{1}$ does not receive positive curvature from $\Delta_{2}$. In Figure 5.6(iii) $\Delta_{2}$ does not receive positive curvature from $\Delta_{3}$. Since $\hat{\Delta}$ receives $\frac{\pi}{2}$ across the $bd$ edge and $\frac{\pi}{4}$ across the $ca$, $ba^{-1}$ edges it follows that $\hat{\Delta}$ receives an average of $\frac{\pi}{4}$ across each of its edges, so $d(\hat{\Delta})\ge 8$ implies that $c^{*}(\hat{\Delta})\le0$.
It remains to study $5 \le d(\hat{\Delta})\le7$.
Checking shows that if $d(\hat{\Delta})=5$ then either the label contradicts $|d|=3$ or $H$ non-cyclic or $\hat{\Delta}$ receives positive curvature across at most one edge and so $c^{*}(\hat{\Delta})\le0$, except when $l(\hat{\Delta})=bdda^{-1}c^{-1}$ as in Figure 5.7(i).
In this case $\hat{\Delta}$ receives $\frac{\pi}{2}$ from $c(\hat{\Delta}_{1})$ and $\frac{\pi}{4}$ from $c(\hat{\Delta}_{2})$. If $|b|>2$ then this implies that $d(\hat{\Delta}_{3})>4$ and so add $\frac{\pi}{4}$ to $c(\hat{\Delta}_{3})$ noting that this is a similar edge to the one crossed in Figure 5.6(iii) so there is no change to the above argument and $c^{*}(\hat{\Delta})\le0$ in this case.
Suppose now that $|b|=2$ and $l(\hat{\Delta}_{3})=ab^{-1}ab^{-1}$ as in Figure 5.7(ii). If  

\begin{figure}
\begin{center}
\psfig{file=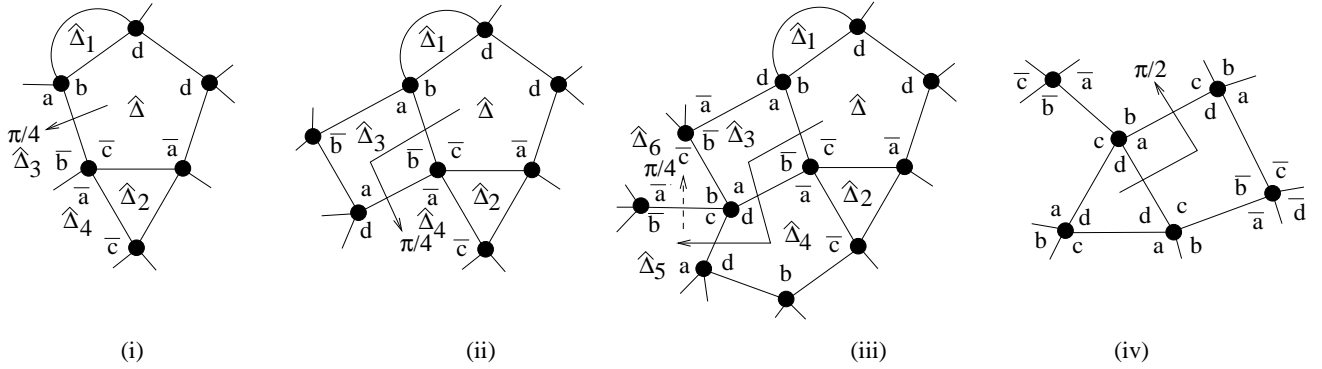}
\end{center}  
\caption{curvature distribution for Case (B3)}
\end{figure}

\noindent $d(\hat{\Delta}_{4})>5$ then add $\frac{\pi}{4}$ to $c(\hat{\Delta}_{4})$ across  the $da^{-1}$ edge. If $d(\hat{\Delta}_{4})=5$ then $l(\hat{\Delta}_{4})=da^{-1}c^{-1}bd$ which implies that $l(\hat{\Delta}_{5})=ca\omega$ 
and so if $d(\hat{\Delta}_{5})>5$ then add 
$\frac{\pi}{4}$ to $c(\hat{\Delta}_{5})$ as in Figure 5.7(iii). If $d(\hat{\Delta}_{5}) \in \{ 4,5 \}$ then $l(\hat{\Delta}_{5}) \in \{ cadb^{-1}, cad^{-2} b^{-1} \}$ and this forces $l(\hat{\Delta}_{6})=c^{-1}ba^{-1}\omega$ and so $d(\hat{\Delta}_{6})>5$ otherwise there is a contradiction to $|c|\ne1$. So add $\frac{\pi}{4}$ to $c(\hat{\Delta}_{6})$ again as shown in Figure 5.7(iii).

Observe that $\hat{\Delta}_{4}$ in Figure 5.7(ii) can now receive $\frac{\pi}{4}$ from $c(\hat{\Delta})$, however it receives no positive curvature from $\hat{\Delta}_{3}$ or any other region across the $da^{-1}$ edge. Moreover, it is clear from Figure 5.7(iii) that $\hat{\Delta}_{5}$ receives only the $\frac{\pi}{4}$ from $\hat{\Delta}_{4}$ across its $ca$ edge; and $\hat{\Delta}_{6}$ receives only the $\frac{\pi}{4}$ from $\hat{\Delta}_{5}$ across its $ba^{-1}$ edge. Finally observe that Figures 5.7(ii)--(iii) do not alter the fact that $\Delta_{1}$ does not receive positive curvature from $\Delta_{3}$ and $\Delta_{2}$ does not receive positive curvature from $\Delta_{4}$ in Figure 5.1(ii). Therefore the average positive curvature that $\hat{\Delta}$ receives across each edge is still $\frac{\pi}{4}$ and so if $d(\hat{\Delta})\ge 8$ then $c^{*}(\hat{\Delta})\le0$. It remains to check $6 \le d(\hat{\Delta}) \le 7$ for the sublabels $(bd)^{\pm 1} ( \pi /2)$ and $(ca)^{\pm 1}, (ab^{-1})^{\pm 1}$, $(da^{-1} c^{-1})^{\pm 1} (\pi /4)$.
Checking shows that if $d(\hat{\Delta})=6$ then the most curvature that $\hat{\Delta}$ can receive is either $2(\frac{\pi}{2})$ or
$\frac{\pi}{2} + 2( \frac{\pi}{4})$ or $4(\frac{\pi}{4})$ and so $c^{\ast} (\hat{\Delta}) \leq 0$.  If $d(\hat{\Delta})=7$ then the most curvature received is $3( \frac{\pi}{2})$ or $2( \frac{\pi}{2}) + 2(\frac{\pi}{4})$ or $\frac{\pi}{2} + 4 ( \frac{\pi}{4})$ or $6(\frac{\pi}{4})$ and $c^{\ast}(\hat{\Delta})\leq 0$ except for $l(\hat{\Delta})=da^{-1} c^{-1} bda^{-1} b$; but this implies $cd=1$, a contradiction.
\end{enumerate}

Now let $H$ be cyclic.  Then $d=b^4$ and $c=b^3$.  Again add $\frac{1}{2} c(\Delta)=\frac{\pi}{2}$ to each of
$c(\Delta_1)$, $c(\Delta_2)$ as in Figure 5.1(ii).  In Figure 4.3(i) if say $d(\Delta_1)>5$ then add
$c(\Delta)=\frac{\pi}{2}$ to $c(\Delta_1)$ and so it can be assumed that $d(\Delta_i) \leq 5$ for $i \in \{1,3,5\}$ in which case $l(\Delta_i) \in \{ cad b^{-1}, cad^{-2} b^{-1} \}$.  If say $d(\Delta_1)=4$ then add
$c(\Delta)=\frac{\pi}{2}$ to $c(\Delta_6)$ via $\Delta_1$ as shown in Figure 5.7(iv).  It can be assumed then that
$d(\Delta_i)=5$ for $i \in \{ 1,3,5 \}$ in which case add $\frac{1}{3} c(\Delta)=\frac{\pi}{6}$ to each $c(\hat{\Delta})$ via $\Delta_i$ where $i \in \{ 1,3,5 \}$ as shown in Figure 5.8(i).  If say $\hat{\Delta}=\hat{\Delta}_1$ and $d(\hat{\Delta}_1)=5$ then repeat the above, that is, add the $\frac{\pi}{6}$ from $c(\Delta)$ across another $ca$ edge and continue in this way until $\frac{\pi}{6}$ is eventually added to a region $\hat{\Delta}_k$ where either $d(\hat{\Delta}_k)>5$ (and so the

\begin{figure}
\begin{center}
\psfig{file=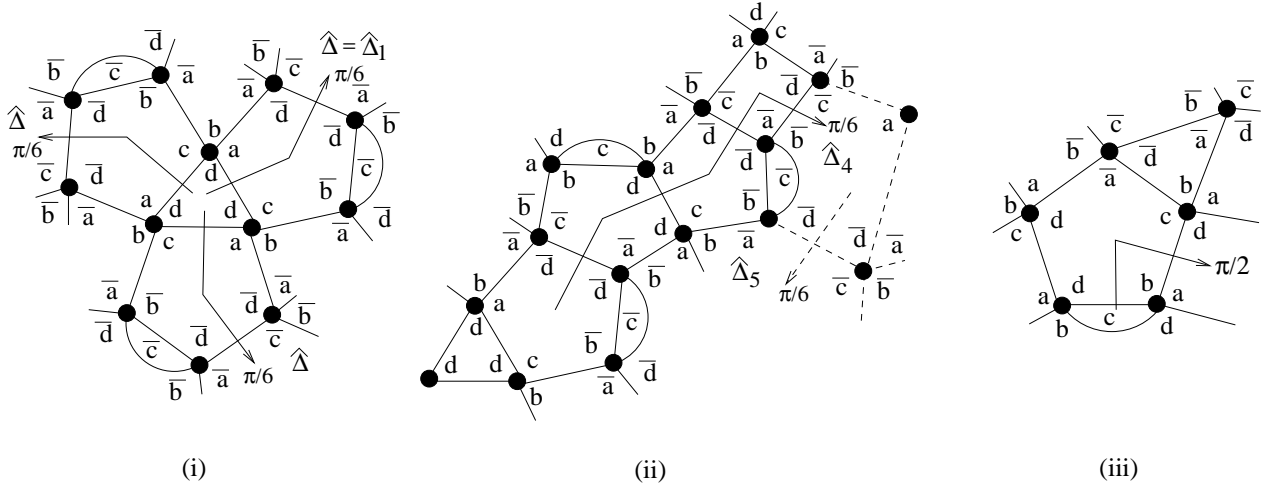}
\end{center}  
\caption{curvature distribution for Cases (B3) and (B5)}
\end{figure}

\noindent process terminates) or
$d(\hat{\Delta}_k)=4$ in which case the $\frac{\pi}{6}$ from $c(\Delta)$ is added to $c(\hat{\Delta}_{k+1})$ as shown in Figure 5.8(ii), where $k=3$.  If $d(\hat{\Delta}_{k+1})>5$ then the process terminates (and note that $l(\hat{\Delta}_{k+1})=d^{-1} b^{-1} c^{-1} w$); otherwise
$l(\hat{\Delta}_{k+1})=d^{-1} b^{-1} c^{-1} ad^{-1}$ and the $\frac{\pi}{6}$ from $c(\Delta)$ is added to $c(\hat{\Delta}_{k+2})$ where $\hat{\Delta}_{k+2}$ is the region shown in Figure 5.8(ii) with $k=3$.  Observe that
$l(\hat{\Delta}_{k+2})=ba^{-1} c^{-1} w$ so $d(\hat{\Delta}_{k+2}) >5$ and the process terminates.  This completes the distribution of curvature that occurs.  It follows that if $\hat{\Delta}$ receives positive curvature across an edge $e_i$ say then $\hat{\Delta}$ does not receive any curvature across the adjacent edges $e_{i-1}$, $e_{i+1}$ except when $\hat{\Delta}$ is given by $\hat{\Delta}_{k+1}=\hat{\Delta}_4$ in Figure 5.8(ii).  Therefore if $l(\hat{\Delta})$ does not involve $(cbd)^{\pm 1}$ then Lemma 3.5(vi) applies and $c^{\ast}(\hat{\Delta}) \leq 0$ for $d(\hat{\Delta}) \geq 7$; and if $d(\hat{\Delta})=6$ then checking for $(bd)^{\pm 1}$, $(ca)^{\pm 1}$ and $(cba^{-1})^{\pm 1}$ shows that $\hat{\Delta}$ receives positive curvature across at most two edges and $c^{\ast}(\hat{\Delta}) \leq 0$.  Finally if $l(\hat{\Delta})=cbdw$ then we see from Figure 5.8(ii) that the maximum amount $\hat{\Delta}$ receives is on average $\frac{\pi}{3}$ across $\frac{2}{3}$ of its edges and so if
$d(\hat{\Delta}) \geq 8$ then $c^{\ast}(\hat{\Delta}) \leq 0$ by Lemma 3.5(vii).  Checking shows that if
$6 \leq d(\hat{\Delta}) \leq 7$ then $l(\hat{\Delta}) \in \{ cbdb^{-1} ab^{-1}, cbda^{-1} bdb^{-1}, cbda^{-1} c^{-1} ba^{-1}, cbda^{-1} cba^{-1} \}$ and so if $d(\hat{\Delta})=6,7$ then $\hat{\Delta}$ receives curvature across at most $2,3$ edges (respectively) and $c^{\ast}(\hat{\Delta}) \leq 0$.

In conclusion $\mathcal{P}$ fails to be aspherical in this case when $H$ is non-cyclic, $bda^{-1} c^{-1}=1$ and $|b| \in \{ 2,3,4,5 \}$; or when $H$ is cyclic and $bda^{-1} c^{\pm 1}=1$.

\textbf{(B4)} $\boldsymbol{|c|=2}$, $\bs{|d|>3,~ a^{-1}b\ne1,~ c^{\pm1}ab^{-1}=1,~ d^{\pm1}b^{-1}a\ne1}$.\newline
If $|d|<\infty$ then there is a sphere by Lemma 3.1(a)(i) and if $|d|=\infty$ then $\mathcal{P}$ is aspherical by Lemma 3.4(ii).

\textbf{(B5)} $\boldsymbol{|c|=2}$, $\bs{|d|>3}$, $\bs{a^{-1}b\ne1}$, $\bs{c^{\pm1}ab^{-1}\ne1}$, $\bs{db^{-1}a=1}$, $\bs{d^{-1}b^{-1}a\ne1}$.\newline
If $d(\Delta)=2$ then $\Delta$ is given by Figure 5.1(ii) in which case add $\frac{1}{2}c(\Delta)=\frac{\pi}{2}$ to each of $c(\Delta_1)$ and $c(\Delta_2)$ as shown; and if $d(\Delta)=3$ then $\Delta$ is given by Figure 4.10(iii) in which case add $c(\Delta)=\frac{\pi}{2}$ to $c(\Delta_5)$ as shown.  Assume that $H$ is non-cyclic.  If say $d(\Delta_1)=4$ in Figure 5.1(ii) then $l(\Delta_1) \in \{ bd^2 a^{-1}, bda^{-1} c^{\pm 1}, bdb^{-1} c^{\pm 1} \}$ which contradicts
$|d| > 3$ or $H$ non-cyclic; and if $d(\Delta_5)=4$ in Figure 4.10(iii) then $l(\Delta_5) \in \{ c^{-1} ad^{\pm 1} a^{-1}, c^{-1} ad^{\pm 1} b^{-1} \}$ which contradicts $|d| > 3$ or $H$ non-cyclic.  Observe that in Figure 5.1(ii) $\Delta_1$ say does not receive positive curvature from $\Delta_3$ or $\Delta_4$; and in Figure 4.10(iii) $\Delta_5$ does not receive any from $\Delta_4$ or $\Delta_6$.  Therefore $c^{\ast} (\hat{\Delta}) \leq 0$ for $d(\hat{\Delta}) \geq 7$ by Lemma 3.5(vi) and so it remains to consider $5 \leq d(\hat{\Delta}) \leq 6$.  But checking for subwords
$(bd)^{\pm 1}, (c^{-1} a)^{\pm 1}$ shows that if $d(\hat{\Delta})=5,6$ then $\hat{\Delta}$ receives positive curvature across $1,2$ edges respectively except when $l(\hat{\Delta})=d^2 a^{-1} cb$.  But $H$ would then be cyclic, so $c^{\ast}(\hat{\Delta}) \leq 0$.

Now assume that $H$ is cyclic.  If $bda^{-1}c=1$ then the conditions are $T$-equivalent to those of Lemma 3.1(a)(ii) and there is a sphere, so assume otherwise.  Then $d(\Delta_1) > 4$ in Figure 5.1(ii) and $d(\Delta_5) > 4$ in Figure 4.10(iii).  The argument above now applies except when $l(\Delta_1) = bd^2 a^{-1} c$ in which case add $\frac{1}{2} c(\Delta)=\frac{\pi}{2}$ to $c(\Delta_3)$ via $\Delta_1$ as shown in Figure 5.8(iii).  We see from Figures 4.10(iii), 5.1(ii) and 5.8(iii) that if $\hat{\Delta}$ receives positive curvature then it does so across at most half of its edges and so $c^{\ast}(\hat{\Delta}) \leq 0$ for $d(\hat{\Delta}) \geq 7$ by Lemma 3.5(vi), so let $4 \leq d(\hat{\Delta}) \leq 6$.  Then either $l(\hat{\Delta}) = bd^2 a^{-1} c^{-1}$ and $c^{\ast}(\hat{\Delta}) \leq c(\Delta) + \frac{\pi}{2} = 0$; or
$d(\hat{\Delta})=6$, $\hat{\Delta}$ receives across at most two edges (checking for $(bd)^{\pm 1}, (c^{-1}a)^{\pm 1}, (ad)^{\pm 1}$) and again $c^{\ast}(\hat{\Delta}) \leq 0$.

In conclusion $\mathcal{P}$ is aspherical if and only if $bda^{-1} c \neq 1$.

\textbf{(B6)} $\boldsymbol{|c|=2}$, $\bs{|d|=3}$, $\bs{a^{-1}b\ne1}$, $\bs{c^{\pm1}ab^{-1}\ne1}$, $\bs{db^{-1}a=1}$, $\bs{d^{-1}b^{-1}a\ne1}$.\newline
If $d(\Delta)=2$ then $\Delta$ is given by Figure 5.1(ii). If $d(\Delta)=3$ then $\Delta$ is given by
Figures 4.3(i) and 4.10(iii).
In Figure 5.1(ii) if $d(\Delta_{1})>4$ and $d(\Delta_{2})>4$ then add $\frac{1}{2}c(\Delta) =\frac{\pi}{2}$ to $c(\Delta_{1})$ and $c(\Delta_{2})$. If $d(\Delta_{1})>4$ and $d(\Delta_{2})=4$ then $l(\Delta_{2})\in\{bd^{2}a^{-1},bda^{-1}c^{\pm1},bdb^{-1}c^{\pm1}\}$ and so $l(\Delta_{2})=bd^{2}a^{-1}$ otherwise there is a contradiction to $|d|=3$. Therefore $\Delta_{2}$ is given by Figure 5.1(iii) forcing $l(\Delta_{3})=ca\omega$ and so $d(\Delta_{3})>4$, otherwise there is a contradiction to $|c|=2$, so add $\frac{\pi}{2}$ to $c(\Delta_{1})$ and to  $c(\Delta_{3})$ via $\Delta_2$. If $d(\Delta_{1})=4$ then add $\frac{1}{2}c(\Delta)=\frac{\pi}{2}$ to $c(\Delta_4)$ via $\Delta_1$. Observe in Figure 4.10(iii) that if $d(\Delta_{5})=4 $ then $l(\Delta_{5})\in\{c^{-1}ad^{\pm1}a^{-1},c^{-1}ad^{\pm1}b^{-1}\}$ which contradicts $|c|=2$ so add $c(\Delta)=\frac{\pi}{2}$ to $c(\Delta_{5})$ as shown.
In Figure 4.3(i) if $d(\Delta_{i})=4$ where $i \in \{ 1,3,5 \}$ then $l(\Delta_{i})\in\{cad^{\pm1}a^{-1},cad^{\pm1}b^{-1}\}$ which contradicts $|c|=2$ so add $c(\Delta)=\frac{\pi}{6}$ to $c(\Delta_{i})$ as shown.
We see from Figures 4.3(i), 4.10(iii) and 5.1(ii)--(iii) that if $\hat{\Delta}$ receives positive curvature then it does so across the edges $bd$, $ca$ or $c^{-1}a$. It follows that the only word of length 3 that contains no gaps is $a^{-1}c^{\pm1}a$. Suppose that $l(\hat{\Delta})=\omega_{1} a^{-1}c^{\pm1}a\omega_{2}\omega$ where $\omega_{1}$ and $\omega_{2}$ have length 2 and $\omega$ has length at least 0. Then there is a gap preceding $a^{-1}$ in $\omega_{1}a^{-1}$ and after $a$ in $a\omega_{2}$. Moreover, if $\omega_{1}$ does not contain a gap then $\omega_{1}=bd$ and if $\omega_{2}$ does not contain a gap then $\omega_{2}=d^{-1}b^{-1}$. It follows that if $l(\omega)>0$ then $l(\hat{\Delta})$ contain at least 4 gaps and $c^{*}(\hat{\Delta})\le0$. If $l(\omega)=0$ and $l(\hat{\Delta})=\omega_{1}a^{-1}c^{\pm1}ad^{-1}b^{-1}$ then $l(\hat{\Delta})\in\{ c^{\pm1}ba^{-1}c^{\pm1}ad^{-1}b^{-1}, ad^{\pm1}a^{-1}c^{\pm1}ad^{-1}b^{-1}\}$ and there are 4 gaps; if $l(\omega)=0$ and $l(\hat{\Delta})=bda^{-1}c^{\pm1}a\omega_{2}$ then $l(\hat{\Delta})\in\{bda^{-1}c^{\pm1}ad^{\pm1}a^{-1},bda^{-1}c^{\pm1}ab^{-1}c^{\pm1}\}$ and again there are 4 gaps. If $l(\hat{\Delta})=a^{-1}c^{\pm1}a\omega$ and $d(\hat{\Delta})\le6$ then either we obtain a contradiction to $|c|=2$ or $l(\hat{\Delta})=a^{-1}c^{\pm1}ab^{-1}c^{\pm1}b$ and there are 4 gaps. Now suppose that $l(\hat{\Delta})$ does not contain the subword $a^{-1}c^{\pm1}a$. Then $\hat{\Delta}$ receives positive curvature across at most half of its edges and so $d(\hat{\Delta})\ge7$ implies $c^{*}(\hat{\Delta})\le0$ by Lemma 3.5(vi). So let $d(\hat{\Delta})\le6$ and 

\begin{figure}
\begin{center}
\psfig{file=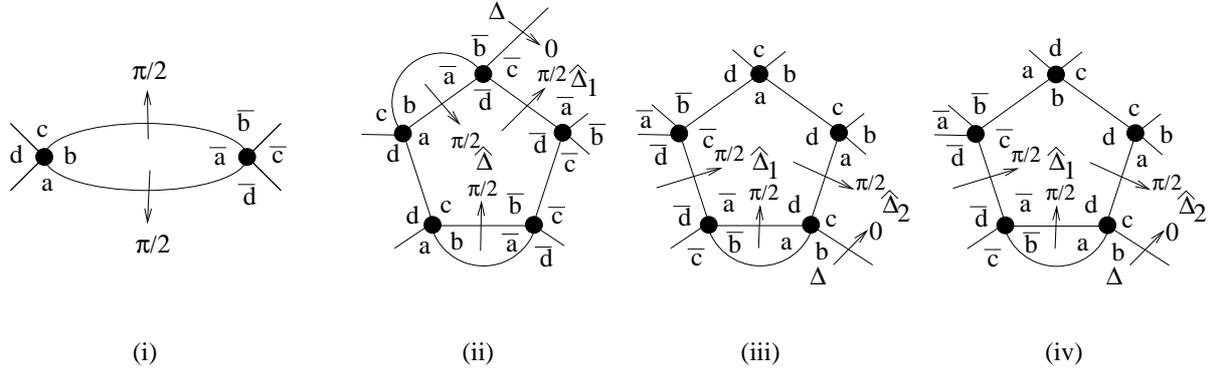}
\end{center}
\caption{curvature distribution for Cases (B7,8,9,11,12)}
\end{figure}

\noindent $l(\hat{\Delta})\in\{ca\omega,bd\omega,c^{-1}a\omega\}$. If $d(\hat{\Delta})<6$ then checking shows that there is a contradiction to a (\textbf{B6}) assumption or $l(\hat{\Delta})=bda^{-1} ba^{-1}$ and there are four 
gaps; and if 
$d(\hat{\Delta})=6$ then checking for $c^{-1}a\omega$ shows that either there is a contradiction to a (\textbf{B6}) assumption or $\hat{\Delta}$ receives positive curvature across at most two edges and so $c^{*}(\hat{\Delta})\le0$. Therefore $\mathcal{P}$ is aspherical.

\textbf{(B7)} $\boldsymbol{|c|>3}$, $\bs{|d|>3}$, $\bs{a^{-1}b=1}$, $\bs{c^{\pm1}ab^{-1}\ne1}$, $\bs{d^{\pm1}b^{-1}a\ne1}$.\newline
If $d(\Delta)=2$ then $\Delta$ is given by Figure 5.9(i).
Here $l(\Delta_{1})=b^{-1}c\omega$ and $l(\Delta_{2})=ad^{-1}\omega$.
First assume that $H$ is non-cyclic.
This implies that $d(\Delta_{1})>4$, $d(\Delta_{2})>4$, otherwise there is a contradiction to $H$ non-cyclic, $|c|>3$, or $|d|>3$. So add $\frac{1}{2} c(\Delta)=\frac{\pi}{2}$ to $c(\Delta_{1})$ and $c(\Delta_{2})$.
Observe that each of $\Delta_{1}$ and $\Delta_{2}$ does not receive positive curvature from $\Delta_{j}$ for $j\in\{3,4\}$. It follows that if $\hat{\Delta}$ receives positive curvature then it does so across at most half of its edges and so $d(\hat{\Delta})\ge7$ implies that $c^{*}(\hat{\Delta})\le0$ by Lemma 3.5(vi).
Checking shows that if $d(\hat{\Delta})=5$ then $l(\hat{\Delta})$ contradicts $|c|>3$, $|d|>3$ or $H$ non-cyclic. Also if $d(\hat{\Delta})=6$ then $\hat{\Delta}$ receives positive curvature across at most two edges and so $c^{*}(\hat{\Delta})\le0$.

Now assume that $H$ is cyclic.
If $c=d^{\pm 1}$ there is a sphere by Lemma 3.1(a)(iii) and Lemma 3.2(iii) so assume from now on that $c \neq d^{\pm 1}$.  If $c \neq d^2$ and $d \neq c^2$ then $\mathcal{P}$ is aspherical by the above argument.  If $c=d^2$ and $d=c^2$ then $c^3=1$, a contradiction, so assume that $c=d^2$ and $d \neq c^2$.  Then
$c^{\ast}(\hat{\Delta}) \leq 0$ as above except when $l(\hat{\Delta})=b^{-1} cad^{-2}$ and $\hat{\Delta}$ is given by Figure 5.9(ii).  In this case add $c^{\ast}(\hat{\Delta}) = \frac{\pi}{2}$ to $c(\hat{\Delta}_1)$ as shown noting that $\hat{\Delta}_1$ does not receive any curvature from $\Delta$.  Therefore if a region receives positive curvature across an edge it is across a $b^{-1} c$ or $ad^{-1}$ or $a^{-1} c^{-1}$ edge; and if across consecutive edges then the sublabel is $(da^{-1}c^{-1})^{\pm 1}$.
It follows that if $l(\hat{\Delta})$ does not involve $(da^{-1} c^{-1})^{\pm 1}$ and $d(\hat{\Delta}) \geq 7$ then $c^{\ast}(\hat{\Delta}) \leq 0$.  Observe that each occurrence of $da^{-1} c^{-1}$ contributes two gaps.  It follows that if $l(\hat{\Delta})$ contains at least two occurrences of $(da^{-1} c^{-1})^{\pm 1}$ or if $l(\hat{\Delta})$ contains exactly one occurrence and $d(\hat{\Delta}) \geq 8$ then there are four gaps and $c^{\ast}(\hat{\Delta}) \leq 0$.  But if $l(\hat{\Delta})=da^{-1} c^{-1} w$ and $d(\hat{\Delta})=7$ then checking possible $l(\hat{\Delta})$ shows that again there are four gaps and $c^{\ast}(\hat{\Delta}) \leq 0$.
If $5 \leq d(\hat{\Delta}_1) \leq 6$ then either there is a contradiction to one of our assumptions or $c^{\ast}(\hat{\Delta}_1) \leq 0$ except when
$l(\hat{\Delta}_1) \in \{ a^{-1} c^{-2} ad^{-1}, a^{-1} c^{-2} bd^{-1}, a^{-1} c^{-1} ad^2, a^{-1} c^{-1} bd^2, a^{-1} c^{-3} bd \}$. If $d=c^{-2}$ or $d=c^3$ then $|c|=5$ and there is a sphere by Lemma 3.2(vii).  This leaves $l(\hat{\Delta}_1) \in \{ a^{-1} c^{-1} ad^2, a^{-1} c^{-1} bd^2 \}$ and these are given by Figure 5.9(iii), (iv).  Add $c^{\ast} (\hat{\Delta}_1) = \frac{\pi}{2}$ to $c(\hat{\Delta}_2)$ as shown.  Repeat the argument for $\hat{\Delta}_2$ noting that as before $\hat{\Delta}_2$ does not receive positive curvature from the corresponding region $\Delta$.  This procedure will terminate at a region $\hat{\Delta}_k$, say, such that $\hat{\Delta}_i \neq \hat{\Delta}_j$ for $i \neq j$ where $1 \leq i,j \leq k$ and $c^{\ast} (\hat{\Delta}_k) \leq 0$.  Repeat this for each copy of the region $\hat{\Delta}_1$ to conclude that $\mathcal{P}$ is aspherical.  If $d=c^2$ and $c \neq d^2$ then the argument is the same by symmetry.

In conclusion $\mathcal{P}$ is aspherical in this case except when either $c=d^{\pm 1}$ or $c^5=1$ and $c=d^2$ or $d^5=1$ and $d=c^2$.

\textbf{(B8)} $\boldsymbol{|c|=2}$, $\bs{|d|>3}$, $\bs{a^{-1}b=1}$, $\bs{c^{\pm1}ab^{-1}\ne1}$, $\bs{d^{\pm1}b^{-1}a\ne1}$.\\
If $d(\Delta)=2$ then $\Delta$ is given by Figures 5.1(ii) and 5.9(i).
In Figure 5.9(i) $l(\Delta_{1})=b^{-1}c\omega$ and $l(\Delta_{2})=ad^{-1}\omega$. This implies that $d(\Delta_{1})>4$,$d(\Delta_{2})>4$, otherwise
there is a contradiction to $|d| > 3$ so add $\frac{1}{2} c(\Delta)=\frac{\pi}{2}$ to $c(\Delta_{1})$ and $c(\Delta_{2})$.
In Figure 5.1(ii) $l(\Delta_{1})=l(\Delta_{2})=bd\omega$. This similarly implies that $d(\Delta_{1})>4$ and $d(\Delta_{2})>4$, so add $\frac{1}{2} c(\Delta)=\frac{\pi}{2}$ to $c(\Delta_{1})$ and $c(\Delta_{2})$.
Observe that in Figure 5.9(i) $\Delta_{1}$ does not receive positive curvature from $\Delta_{4}$; and $\Delta_{2}$ does not receive positive curvature from $\Delta_{4}$.
Observe also that if $\hat{\Delta}$ receives positive curvature then it does so across the edges $b^{-1}c$, $ad^{-1}$ or $bd$. Thus there is always a gap immediately preceding $c^{-1}$ and $a$; and there is a gap after $c$ and $a^{-1}$. This implies that if there are at least four occurrences of $c^{\pm1}$ then $l(\hat{\Delta})$ contains at least four gaps and so $c^{*}(\hat{\Delta})\le0$ by Lemma 3.6. We will proceed according to the number of occurrences of $c^{\pm1}$ in $l(\hat{\Delta})$.
If there are no occurrences of $c^{\pm1}$ then either $l(\hat{\Delta})=(ab^{-1})^{k}$ where $|k|\ge 2$ and there are four gaps, or $l(\hat{\Delta})=d(ab^{-1})^{k_{1}} \ldots d(ab^{-1})^{k_{m}}$ where $k_{i}\in \mathbb{Z}$ $(1 \le i \le m)$. But since there is always at least one gap between any two occurrences of $d^{\pm1}$ it follows that again there are four gaps or $|d|\le 3$, a contradiction. So $c^{*}(\hat{\Delta})\le 0$ in this case.

Assume first that $H$ is non-cyclic.
If there is exactly one occurrence of $c^{\pm1}$ in $l(\hat{\Delta})$ then $H$ is cyclic so suppose
that there are either two or three occurrences of $c^{\pm1}$. Then either the label contains at least four gaps or it contradicts one of the (\textbf{B8}) assumptions or one of the following cases $\hat{\Delta}_i$ ($1 \leq i \leq 9$) occurs:
\begin{enumerate}
\item[(1)]
$cad^{-1}b^{-1}cad^{-1}b^{-1}$;
\item[(2)]
$cad^{-1}b^{-1}cbd^{-1}b^{-1}$;
\item[(3)]
$cad^{-1}b^{-1}c^{-1}ad^{-1}b^{-1}$;
\item[(4)]
$cad^{-1}b^{-1}c^{-1}bd^{-1}b^{-1}$;
\item[(5)]
$cad^{-1}b^{-1}cbda^{-1}$;
\item[(6)]
$cad^{-1}b^{-1}cbdb^{-1}$;
\item[(7)]
$cad^{-1}b^{-1}c^{-1}bda^{-1}$;

\begin{figure}
\begin{center}
\psfig{file=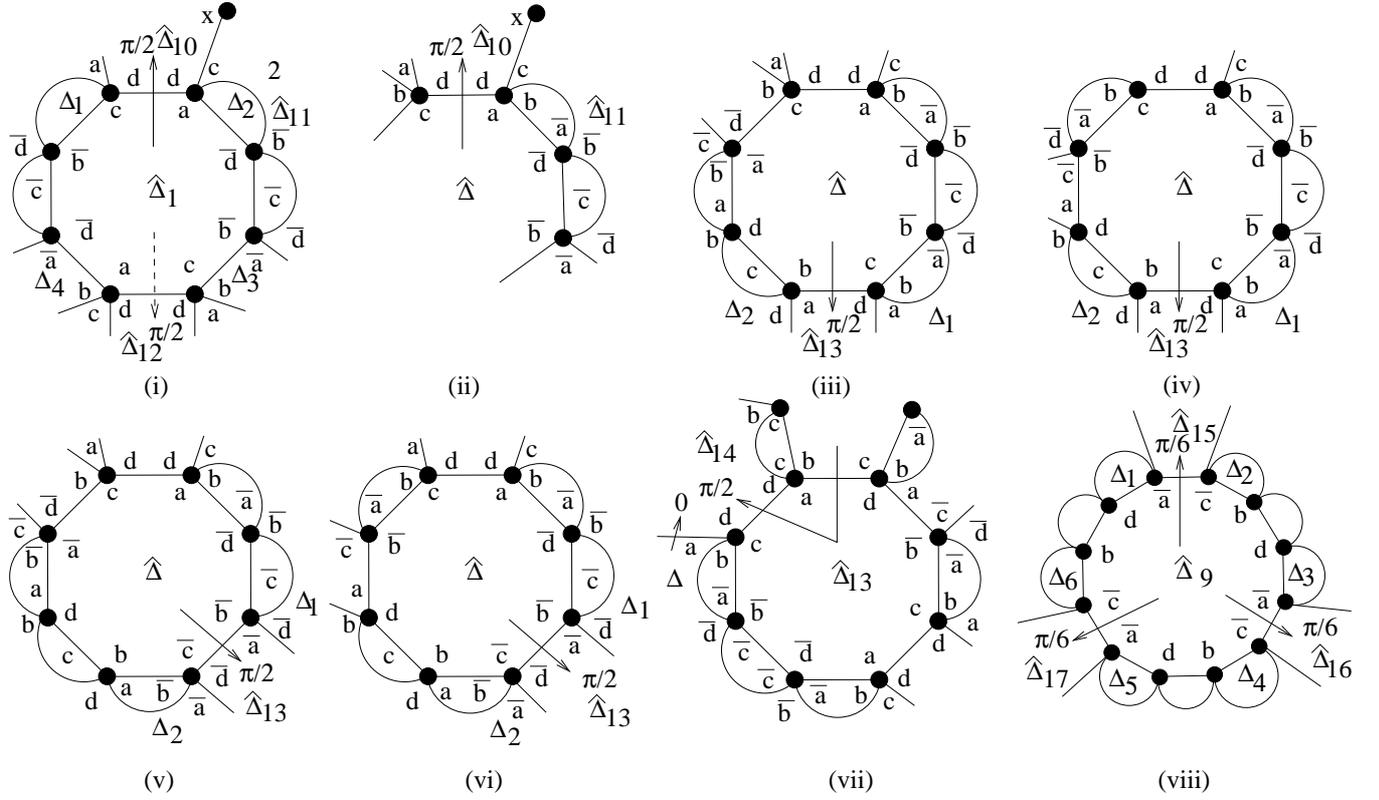}
\end{center}
\caption{curvature distribution for Case (B8)}
\end{figure}

\item[(8)]
$cad^{-1}b^{-1}c^{-1}bdb^{-1}$;
\item[(9)]
$(bda^{-1}c^{-1})^{3}$.
\end{enumerate}

If any of (1)--(4) occurs with any of (5)--(8) or with (9) then $|d|=2$, a contradiction.  Also if any of (5)--(8) occurs with (9) then $c=d^3$ and $H$ is cyclic, so assume otherwise.

Consider (1)--(4).  These yield the relator $(cd)^2$ and it follows that $cd^k = d^{-k} c$ for $k \in \mathbb{Z}$.  Moreover if $|d| < \infty$ then there is a sphere by Lemma 3.1(c)(ii) so it can be assumed that $|d| = \infty$.  In case (1) $\hat{\Delta}_1$ is given by Figure 5.10(i) where, given that
$c^{\ast}(\hat{\Delta}_1) > 0$, it can be assumed that $d(\Delta_1)=d(\Delta_2)=2$ and at least one of $d(\Delta_3)$, $d(\Delta_4)$ equals 2.
Add $\frac{\pi}{2}$ from $c(\hat{\Delta}_1)$ to $c(\hat{\Delta}_{10})$ as shown in Figure 5.10(i); and if $d(\Delta_3)=d(\Delta_4)=2$ add a further
$\frac{\pi}{2}$ of $c(\hat{\Delta}_1)$ to $c(\hat{\Delta}_{12})$ as shown.
In cases (2)--(4) $c^{\ast}(\hat{\Delta}) \leq \frac{\pi}{2}$ where
$\hat{\Delta} \in \{ \hat{\Delta}_2, \hat{\Delta}_3, \hat{\Delta}_4 \}$ and $\frac{\pi}{2}$ is added from $c(\hat{\Delta})$ to $c(\hat{\Delta}_{10})$ as shown in Figure 5.10(ii).  Observe that $x \neq b$ in Figure 5.10(i), (ii) for otherwise $c^2$ would be a proper sublabel, and so $x \in \{ a,d \}$.  If $x=a$ then the sublabel $ad$ yields a gap so let $x=d$.  Then either $dd$ yields a gap or $\hat{\Delta}_{11} \in \{ \hat{\Delta}_i \colon 1 \leq i \leq 4 \}$ and
$\frac{\pi}{2}$ is added to $c(\hat{\Delta}_{10})$ from $c(\hat{\Delta}_{11})$.  Continuing this way, since $|d| = \infty$, eventually we get a sublabel $ad$ or $dd$ which contributes a gap.  Consider $l(\hat{\Delta}_{10})$.  If it contains an odd number of occurrences of $c$ then $cd^k = d^{-k}c$ implies that
$c \in \langle d \rangle$ and $H$ is cyclic.  This leaves the case when there are exactly two occurrences of $c$ and $cd^{\alpha_1} cd^{\alpha_2}=1$ for
$\alpha_1,\alpha_2 \in \mathbb{Z} \backslash \{ 0 \}$.  If $| \alpha_1 |, | \alpha_2 | > 1$ then there are four gaps and
$c^{\ast}(\hat{\Delta}_{10}) \leq 0$; and if $| \alpha_1 | > 1$, $| \alpha_2 | = 1$ this implies $|d| < \infty$, a contradiction.

Consider (5)--(8).  These yield the relator $cdcd^{-1}$ and $H$ is Abelian.  Observe that $|d|=4$ yields (\textbf{E}), so assume otherwise.
In each case add $c^{\ast}(\hat{\Delta})=\frac{\pi}{2}$ to $c(\hat{\Delta}_{13})$ as shown in Figure 5.10(iii)--(vi).  Observe that $\hat{\Delta}_{13}$ receives no curvature from $\Delta_1$ or $\Delta_2$; that $l(\hat{\Delta}_{13})=adw$ implies $d(\hat{\Delta}_{13}) > 4$ otherwise there is a contradiction to
$|d| > 3$; and there is still a gap between each pair of occurrences of $d$.  If $l(\hat{\Delta}_{13})$ contains an odd number of occurrences of $c$ then $H$ is cyclic so it can be assumed that $l(\hat{\Delta}_{13})$ yields the relator $cd^{\beta_1} cd^{\beta_2}$.  If $| \beta_1 | > 1$ and $| \beta_2 | > 1$ then there are four gaps and if $(\beta_1,\beta_2) \in \{ (2,1),(2,-1),(1,1)\}$ then $|d| \leq 3$, so this leaves the case $\beta_1=1$, $\beta_2=-1$.  Again there are four gaps except when $l(\hat{\Delta}_{13})=adb^{-1} cad^{-1} b^{-1} c$ and this is shown in Figure 5.10(vii): add $c^{\ast}(\hat{\Delta}_{13}) = \frac{\pi}{2}$ to
$c(\hat{\Delta}_{14})$ and observe that $\hat{\Delta}_{14}$ does not receive positive curvature from $\Delta$.  Consider
$l(\hat{\Delta}_{14})=bddw$.  If there are at least four occurrences of $c$ then $c^{\ast}(\hat{\Delta}_4) \leq 0$; and if there is an odd number of occurrences then $H$ is cyclic.  Suppose firstly that there are no occurrences of $c$ in $l(\hat{\Delta}_{14})$.  Since $|d| \geq 5$, if there is one occurrence of $b$ then
$l(\hat{\Delta}_4)=a^{-1} bd^k$ ($k \geq 5$) and there are four gaps; and since each $(a^{-1}b)^{\pm 1}$ yields a gap and each $(bd^l)^{\pm 1}$
($l \geq 2$) yields a gap it follows that if there are at least two occurrences of $b$ then again $c^{\ast}(\hat{\Delta}_{14}) \leq 0$.  Suppose finally that there are two occurrences of $c$ and so $cd^{\beta_1} cd^{\beta_2}=1$ where $\beta_1 \geq 2$ and $|\beta_2| \geq 0$.  If $| \beta_2 | > 1$ then there are four gaps; and if $|\beta_2|=1$ then $\beta_1 \geq 4$, otherwise there is a contradiction to $|d| > 4$, and again there are four gaps, so $c^{\ast}(\hat{\Delta}_{14}) \leq 0$.

Finally consider case (9). In this case $\hat{\Delta}_{9}$ is given by Figure 5.10(viii).
Suppose that $c^{*}(\hat{\Delta}_{9})> 0$. Then it can be assumed that $d(\Delta_{i})=2$ for $1 \leq i \leq 6$ and $c^{*}(\hat{\Delta})=\frac{\pi}{2}$ so add $\frac{1}{3} c^{*}(\hat{\Delta})=\frac{\pi}{6}$ to $c(\hat{\Delta}_{l})$ for $l\in \{15,16,17\}$.
In this case if $|d|\in\{4,5\}$ then we obtain a sphere by Lemma 3.1(c)(iii). Now if $|d|\ge 6$ then as shown in Figure 5.10(viii) $l(\Delta_{l})=d^{-2}\omega$ and $d^{-2}$ will contribute two-thirds of a gap.  If there are now at least two occurrences of $c$ then either $|d|<6$ or $H$ is cyclic, a contradiction, or there are four gaps; if there is exactly one occurrence of $c$ then this contradicts $H$ non-cyclic; and if there are no occurrences of $c$ then $l(\Delta_{l})=d^{k_{1}}(b^{-1}a)^{m_{1}} \ldots d^{k_{n}}(b^{-1}a)^{m_{n}}$ where $m_{i}\in \mathbb{Z}$,$k_{i}\ge1$. Since $k_{1}+ \ldots +k_{n} \ge6$ it follows that there are at least four gaps and $c^{*}(\hat{\Delta}_{l})\le 0$.

Now let $H$ be cyclic.
If $c=d^2$ or $c=d^3$ then there is a sphere by $T$-equivalence and Lemma 3.2(vii), (viii); and $c=d^4$ is (\textbf{E4}), so assume otherwise.  In particular, $|d| > 4$.  We follow the same argument as above and so if $l(\hat{\Delta})$ contains no occurrences or at least four occurrences of $c$ then, as before,
$c^{\ast} (\hat{\Delta}) \leq 0$; and if $l(\hat{\Delta})$ contains an odd number of occurrences of $c$ then $c=d^k$ for some $k \geq 4$ which implies there are at least four gaps and $c^{\ast}(\hat{\Delta}) \leq 0$.  Suppose then that $l(\hat{\Delta})$ involves $c$ exactly twice.  Subcases (1)--(4) imply $d^2=1$ and (9) implies $c=d^3$, a contradiction.  This leaves subcases (5)--(8).

Add $c^{\ast}(\hat{\Delta}) = \frac{\pi}{2}$ to $c(\hat{\Delta}_{13})$ as in Figure 5.10(iii)-(vi).
Since there is still a gap between each pair of occurrences of $d$ it follows from the above paragraph and the previous argument that $c^{\ast}(\hat{\Delta}_{13}) \leq 0$ except when $l(\hat{\Delta}_{13}) = adb^{-1} cad^{-1} b^{-1} c$.  Again add $c^{\ast} (\hat{\Delta}_{13}) = \frac{\pi}{2}$ to
$c(\hat{\Delta}_{14})$ as shown in Figure 5.10(vii).  If $l(\hat{\Delta}_{14})=bddw$ involves at least three occurrences of $c$ then, since $\hat{\Delta}_4$ does not receive positive curvature from $\Delta$ in Figure 5.10(vii), there are at least four gaps and $c^{\ast}(\hat{\Delta}_4) \leq 0$.  Otherwise checking the possible labels for $l(\hat{\Delta}_4)=bddw$ shows that there are four gaps or a contradiction to $|d| > 4$ or $c \notin \{ d^3,d^4 \}$.

In conclusion $\mathcal{P}$ is aspherical except when $|cd|=2$, $|d| < \infty$ or when $| cd | = 3$, $|d| \in \{ 4,5 \}$ or when $H$ is cyclic and $c=d^2$ or $d^3$.

\textbf{(B9)} $\boldsymbol{|c|=3}$, $\bs{|d|>3}$, $\bs{a^{-1}b=1}$, $\bs{c^{\pm1}ab^{-1}\ne1}$, $\bs{d^{\pm1}b^{-1}a\ne1}$.\newline
If $d(\Delta)=2$ then $\Delta$ is given by Figure 5.9(i). If $d(\Delta)=3$ then $\Delta$ is given by Figure 4.1(ii). In Figure 5.9(i) $l(\Delta_{1})=b^{-1}c\omega$ and $l(\Delta_{2})=ad^{-1}\omega$ which implies that $d(\Delta_{1})>4$ and $d(\Delta_{2})>4$ otherwise there is a contradiction to $|d|>3$ or $|c|=3$, so add $\frac{1}{2} c(\Delta)=\frac{\pi}{2}$ to $c(\Delta_{1})$ and $c(\Delta_{2})$.
In Figure 4.1(ii) $l(\Delta_{1})=l(\Delta_{3})=l(\Delta_{5})=bd\omega$ which implies that $d(\Delta_{i})>4$ where $i\in\{1,3,5\}$ otherwise there is a contradiction to $|d|>3$ so add $\frac{1}{3} c(\Delta)=\frac{\pi}{6}$ to $c(\Delta_{i})$.
Observe that if $\hat{\Delta}$ receives positive curvature then it receives $\frac{\pi}{2}$ across the edges $b^{-1}c$, $ad^{-1}$ or receives $\frac{\pi}{6}$ across the edge $bd$. Thus there is always a gap immediately preceding $c^{-1}$ and $a$; and there is a gap after $c$ and $a^{-1}$. Also there is a two thirds of a gap across the edge $bd$.
This implies that if there are at least four occurrences of $c^{\pm1}$ then $l(\hat{\Delta})$ contains at least four gaps and so $c^{*}(\hat{\Delta})\le0$ by Lemma 3.6. We will proceed according to the number of occurrences of $c^{\pm1}$ in $l(\hat{\Delta})$. If there are no occurrences of $c^{\pm1}$ then either $l(\hat{\Delta})=(ab^{-1})^{k}$ where $|k|\ge 3$ and there are four gaps, or $l(\hat{\Delta})=d(ab^{-1})^{k_{1}} \ldots d(ab^{-1})^{k_{m}}$ where $k_{i}\in \mathbb{Z}$ $(1 \le i \le m)$. But since there is always at least one gap between any two occurrences of $d^{\pm1}$ it follows that again there are four gaps or $|d|\le 3$, a contradiction. So $c^{*}(\hat{\Delta})\le 0$ in this case.
On the other hand if there are between one and three occurrences of $c^{\pm 1}$ inclusive and $c^{\ast}(\hat{\Delta}) > 0$ then either there is a contradiction to one of the (\textbf{B9}) 

\begin{figure}
\begin{center}
\psfig{file=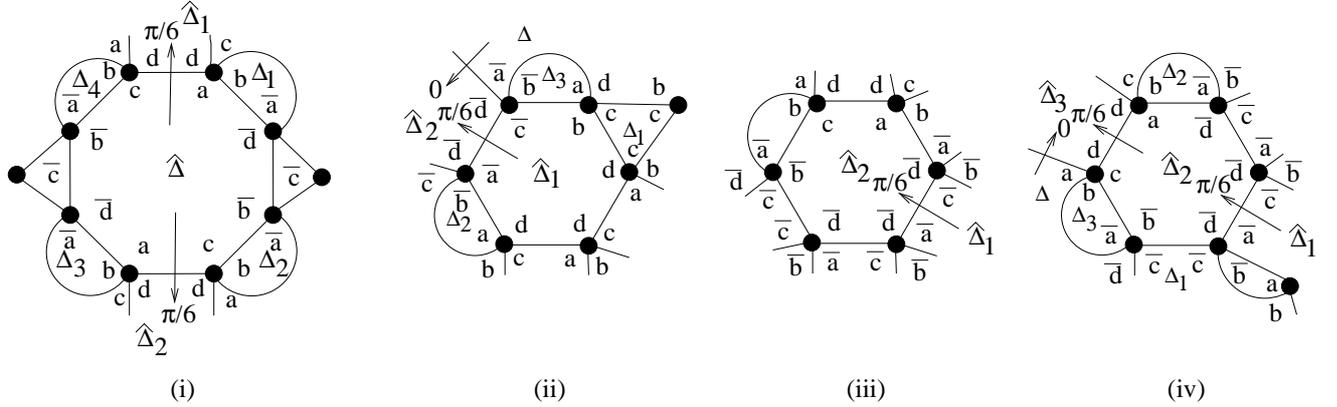}
\end{center}
\caption{curvature distribution for Case (B9)}
\end{figure}

\noindent assumptions or $c \in \{ d^{\pm 2}, d^3 \}$ or
$l(\hat{\Delta})=cad^{-1}b^{-1}cad^{-1}b^{-1}$.
If $c=d^{\pm 2}$ then (\textbf{E2}) or (\textbf{E3}) occurs, so assume otherwise.

First suppose that $l(\hat{\Delta})=cad^{-1}b^{-1}cad^{-1}b^{-1}$ (in particular, $H$ is non-cyclic) as in Figure 5.11(i).
If $c^{*}(\hat{\Delta})> 0$ then it can be assumed that $d(\Delta_{i})=2$ for $1 \leq i \leq 4$ and $c^{*}(\hat{\Delta})=\frac{\pi}{3}$, so add $\frac{1}{2} c^{*}(\hat{\Delta})=\frac{\pi}{6}$ to $c(\hat{\Delta}_{j})$ for $j\in \{1,2\}$. Therefore $d^{\pm2}$ contributes a two-thirds gap.
If $|d|\in\{4,5\}$ then we obtain a sphere by Lemma 3.1(c)(iv).
Assume now that $|d|\ge 6$ and without any loss of generality assume that $\hat{\Delta}_{l}\ne \hat{\Delta}$ of Figure 5.11(i). If $l(\hat{\Delta}_{l})$ has at least two occurrences of $c$ then either there is a contradiction to $H$ non-cyclic or one of the (\textbf{B9}) assumptions or there are four gaps and $c^{*}(\hat{\Delta}_{l})\le0$; if there is exactly one occurrence of $c$ then this contradicts $H$ non-cyclic; and if there are no occurrences of $c$ then $l(\hat{\Delta}_{l})=d^{k_{1}}(b^{-1}a)^{m_{1}}...d^{k_{n}}(b^{-1}a)^{m_{n}}$ where $m_{i}\in \mathbb{Z}$, $k_{i}\ge1$. Since $k_{1}+ \ldots +k_{n} \ge6$ it follows that there are at least four gaps and $c^{*}(\hat{\Delta}_{l})\le 0$.

Now let $c=d^3$.  It follows from the above argument and checking possible vertex labels
that if $\hat{\Delta}_1$ receives positive curvature then
$c^{\ast}(\hat{\Delta}_1) \leq 0$ except when $l(\hat{\Delta}_1) = c^{-1} bd^3 a^{-1}$ and $c^{\ast}(\hat{\Delta}) = \frac{\pi}{6}$.
Add $c^{\ast} (\hat{\Delta}_1)$ to $c(\hat{\Delta}_2)$ as shown in Figure 5.11(ii) (where we note that $\hat{\Delta}_2$ receives no curvature from $\Delta$).
Checking $\{ d^{-1},a,b\} d^{-2} w$ with the understanding that there is a gap preceding the $d^{-2}$ and each $d^{-2}$ contributes a two-thirds gap it follows that either there is a contradiction or
$c^{\ast}(\hat{\Delta}_2) \leq 0$ except
when $l(\hat{\Delta}_2) = l(\hat{\Delta}_1)^{-1} = cad^{-3} b^{-1}$.
If $\hat{\Delta}_2$ is given by Figure 5.11(iii) then $c^{\ast} (\hat{\Delta}_2) \leq c(\hat{\Delta}_2) + \frac{\pi}{2} + 2( \frac{\pi}{6} ) <0$, so let
$\hat{\Delta}_2$ be given by Figure 5.11(iv).  Observe that the neighbour $\Delta_1$ has label $c^{-2} bw$ and no longer has degree 3 and that if either
$d(\Delta_2) > 2$ or $d(\Delta_3)>2$ then $c^{\ast}(\hat{\Delta}) < 0$ and so Figure 5.11(iv) shows the only possibility for $c^{\ast}(\hat{\Delta}_2)>0$.
Now repeat the argument for $\hat{\Delta}_2$ and $\hat{\Delta}_3$.  The conditions on the neighbouring regions $\Delta_1$, $\Delta_2$ and $\Delta_3$ ensure that this procedure will terminate at a region $\hat{\Delta}_k$, say, such that $\hat{\Delta}_i \neq \hat{\Delta}_j$ for $i \neq j$ where $1 \leq i,j \leq k$ and $c^{\ast} (\hat{\Delta}_k ) \leq 0$.  Repeat this for each copy of the region $\hat{\Delta}_1$ to conclude $\mathcal{P}$ is aspherical.

In conclusion $\mathcal{P}$ is aspherical except when $H$ is non-cyclic, $| c^{-1} d | = 2$ and $|d| \in \{ 4,5 \}$.

\begin{figure}
\begin{center}
\psfig{file=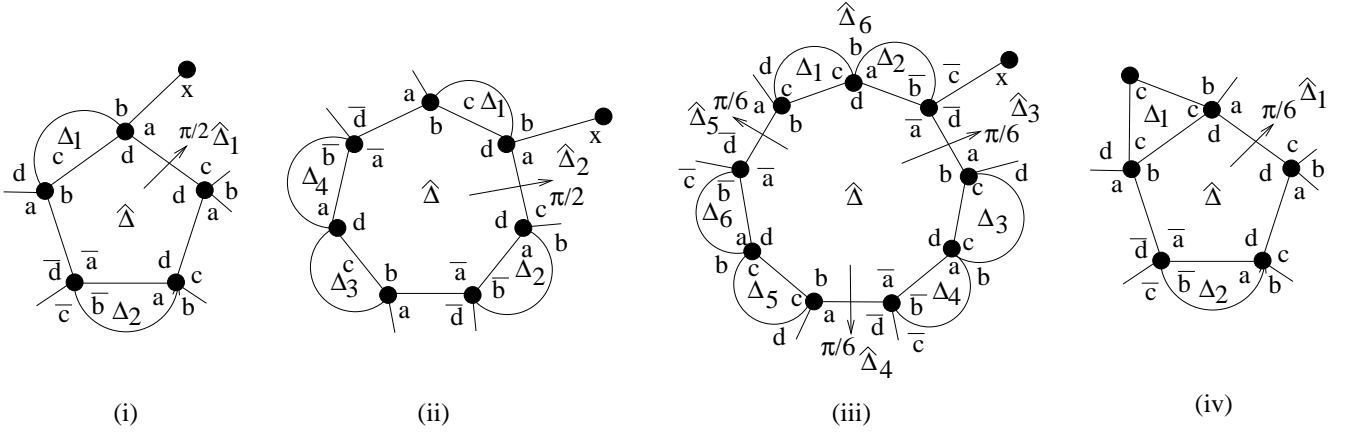}
\end{center}
\caption{curvature distribution for Cases (B11) and (B12)}
\end{figure}

\textbf{(B10)} $\boldsymbol{|c|=2}$, $\bs{|d|=2,~a^{-1}b=1,~c^{\pm1}ab^{-1}\ne1,~ d^{\pm1}b^{-1}a\ne1}$.\newline
If $|cd| < \infty$ then we obtain a sphere by Lemma 3.1(c)(i), so assume otherwise (in particular, $H$ is non-cyclic).
Observe that
if a label involves $c$ and is not $c^{\pm 2}$ then it yields the relation $(cd)^k=1$ and so it follows that if $c(\Delta) >0$ then $\Delta$ is given by Figure 5.1(ii) in which $l(\Delta_1)=l(\Delta_2)=bdw$.
Add $\frac{1}{2} c(\Delta)=\frac{\pi}{2}$ to $c(\Delta_1)$ and $c(\Delta_2)$ as shown.
This completes the distribution of curvature.
Consider $\Delta_1$ of Figure 5.1(ii).
Since positive curvature only crosses the edge $(bd)^{\pm 1}$ it follows that if $l(\Delta_1)$ involves at least four occurrences of $d$ then
$c^{\ast}(\Delta_1) \leq 0$.
This leaves the case when $l(\Delta_1)$ involves exactly two occurrences of $d$.
But $d^2$ cannot be a sublabel since there would then be a region with label $caw$, and so $l(\Delta_1)=bd \{ a^{-1} b, b^{-1} a \} da^{-1}$ has four gaps
and $c^{\ast}(\Delta_1) \leq 0$.
In conclusion $\mathcal{P}$ is aspherical in this case if and only if $|cd|=\infty$.

\textbf{(B11)} $\boldsymbol{|c|=2}$, $\bs{|d|=3,~ a^{-1}b=1,~ c^{\pm1}ab^{-1}\ne1,~ d^{\pm1}b^{-1}a\ne1}$.\newline
If $H$ is cyclic then there is a sphere by Lemma 3.2(vi), so assume otherwise.
If $d(\Delta)=2$ then
$\Delta$ is given by Figures 5.1(ii) and 5.9(i); and if $d(\Delta)=3$ then $\Delta$ is given by Figure 4.3(i).
In Figure 5.9(i) $l(\Delta_{1})=b^{-1}c\omega$ and $l(\Delta_{2})=ad^{-1}\omega$. This implies that $d(\Delta_{1})>4$ and $d(\Delta_{2})>4$ otherwise there is a contradiction to $|d|=3$ so add $\frac{1}{2} c(\Delta)=\frac{\pi}{2}$ to $c(\Delta_{1})$ and $c(\Delta_{2})$.
In Figure 5.1(ii) $l(\Delta_{1})=l(\Delta_{2})=bd\omega$.  This implies that $d(\Delta_{1})>4$ and $d(\Delta_{2})>4$ otherwise there is a contradiction to $|d|=3$ so add $\frac{1}{2} c(\Delta)=\frac{\pi}{2}$ to $c(\Delta_{1})$ and $c(\Delta_{2})$.
In Figure 4.3(i) $l(\Delta_{1})=l(\Delta_{3})= l(\Delta_{5})=ca\omega$. This implies that $d(\Delta_{1})>4$, $d(\Delta_{3})>4$ and $d(\Delta_{5})>4$ otherwise there is a contradiction to $|d|=3$ so add $\frac{1}{3} c(\Delta)=\frac{\pi}{6}$ to $c(\Delta_{1})$, $c(\Delta_{3})$ and $c(\Delta_{5})$.

Observe that if $\hat{\Delta}$ receives positive curvature then it does so across the edges $b^{-1}c$, $ad^{-1}$, $bd$, $ca$ and there is a two-thirds gap in $ca$.
If there is only one occurrence of $c$ or $d$ then this contradicts $H$ non-cyclic.
Also if $|cd|\in\{2,3,4,5\}$ then we obtain a sphere by Lemma 3.1(c)(ii),(v).
So assume that $|cd|>5$.
If there are between two and five occurrences of $c$ inclusive then the resulting relator contradicts $H$ non-cyclic or $| cd | > 5$ except for
$(cd)^2 (cd^{-1})^2$ and $(cdcd^{-1})^2$.  But each of these yields four gaps.
Since at least six occurrences of $c$ yields four gaps, this
leaves the case when $c$ does not occur and this implies that there is at least one gap between any two occurrences of $d$.  It follows that there can only be exactly three occurrences of $d$.  Checking shows there are four gaps
except when the following subcases occur:
\begin{enumerate}
\item[(1)]
$l(\hat{\Delta})=bddda^{-1}$;
\item[(2)]
$l(\hat{\Delta})=bd^{2}a^{-1}bda^{-1}$;
\item[(3)]
$l(\hat{\Delta})=bda^{-1}bda^{-1}bda^{-1}$.
\end{enumerate}
Consider (1), this case is given by Figure 5.12(i).
If $d(\Delta_{2})>2$ then $c^{*}(\hat{\Delta})\le 0$. Assume otherwise and add $c^{*}(\hat{\Delta})=\frac{\pi}{2}$ to $c(\hat{\Delta}_{1})$ as shown. Observe that $l(\hat{\Delta}_{1})=cax\omega$ and $x\in\{b^{-1},d,d^{-1}\}$. If $x\in\{b^{-1},d\}$ then $\hat{\Delta}_{1}$ cannot receive positive curvature across the edges $ab^{-1}$ or $ad$. Suppose that $x=d^{-1}$.  Then $\hat{\Delta}_{1}$ cannot receive positive curvature across the edge $ad^{-1}$ otherwise $l(\Delta_{1})\ne c^2$, a contradiction. So there will be a gap between $c$ and the next occurrence of $d$ and the previous argument does not change.  It follows that $c^{*}(\hat{\Delta}_{1})\le 0$.

Consider (2), this case is given by Figure 5.12(ii).
If $d(\Delta_{2})>2$ or $d(\Delta_{4})>2 $ then $c^{*}(\hat{\Delta})\le 0$.
Assume otherwise and add $c^{*}(\hat{\Delta})=\frac{\pi}{2}$ to $c(\hat{\Delta}_{2})$ as shown.
Observe that $l(\hat{\Delta}_{2})=cax\omega$ and
the same argument as for (1) can be applied to deduce that $c^{*}(\hat{\Delta}_{2})\le 0$.

Consider (3), this case is given by Figure 5.12(iii).
If $d(\Delta_{2})>2$ or $d(\Delta_{4})>2$ or $ d(\Delta_{6})>2$ then $c^{*}(\hat{\Delta})\le 0$.
Assume otherwise and add $\frac{1}{3} c^{*}(\hat{\Delta})=\frac{\pi}{6}$ to $c(\hat{\Delta}_{i})$ for $i\in\{3,4,5\}$ as shown.
Consider $l(\hat{\Delta}_{3})=ad^{-1}x\omega$ where $x\in\{b^{-1},a^{-1},d^{-1}\}$.
If $x=b^{-1}$ then $l(\hat{\Delta}_{6})=c^{-2}bd\omega$, a contradiction.
If $x=a^{-1}$ then exactly as above $l(\hat{\Delta}_{3})$ either contains four gaps
or there is a contradiction to either $H$ non-cyclic, $|cd| > 5$ or one of the (\textbf{B11}) assumptions
and so $c^{*}(\hat{\Delta}_{3})\le 0$.
Finally let $x=d^{-1}$. If $l(\hat{\Delta}_{3})\notin \{ad^{-3}b^{-1},ad^{-2}b^{-1}ad^{-1}b^{-1}\}$ then again as above
$c^{*}(\hat{\Delta}_{3})\le 0$ so assume otherwise.
But if $l(\hat{\Delta}_3)=ad^{-3} b^{-1}$ then we are back in subcase (1) and if $l(\hat{\Delta}_3)=ad^{-2} b^{-1} ad^{-1} b^{-1}$ then we are back in subcase (2), the only difference being that $\hat{\Delta}$ in Figures 5.12(i), (ii) this time receives $\pi /6$ from $\Delta_2$ rather than $\pi /2$.

In conclusion $\mathcal{P}$ is aspherical except when $H$ is non-cyclic and
$| cd | \in \{ 2,3,4,5 \}$ or when $H$ is cyclic.

\textbf{(B12)} $\bs{|c|=3}$, $\bs{|d|=3,~ a^{-1}b=1,~ c^{\pm1}ab^{-1}\ne1,~ d^{\pm1}b^{-1}a\ne1}$.\newline
If $H$ is cyclic then $c=d^{\pm 1}$ and so there is a sphere by Lemma 3.1(a)(iii) or Lemma 3.2(iii), so assume otherwise.
If $d(\Delta)=2$ then $\Delta$ is given by Figure 5.9(i); and if $d(\Delta)=3$ then $\Delta$ is given by Figures 4.1(ii) and 4.3(i).
In Figure 5.9(i) $l(\Delta_{1})=b^{-1}c\omega$ and $l(\Delta_{2})=ad^{-1}\omega$. This implies that $d(\Delta_{1})>4$ $d(\Delta_{2})>4$, otherwise there is a contradiction to $H$ non-cyclic or $|c|=3$ or $|d|=3$ so add $\frac{1}{2} c(\Delta)=\frac{\pi}{2}$ to $c(\Delta_{1})$ and $c(\Delta_{2})$.
In Figure 4.1(ii) $l(\Delta_{1})=l(\Delta_{3})=l(\Delta_{5})=bd\omega$.  This implies that $d(\Delta_{1})>4$, $d(\Delta_{3})>4$ and $d(\Delta_{5})>4$ otherwise there is a contradiction to $H$ non-cyclic or $|c|=3$ or $|d|=3$ so add $\frac{1}{3} c(\Delta)=\frac{\pi}{6}$ to $c(\Delta_{1})$, $c(\Delta_{3})$ and $c(\Delta_{5})$.
In Figure 4.3(i) $l(\Delta_{1})=l(\Delta_{3})= l(\Delta_{5})=ca\omega$. This implies that $d(\Delta_{1})>4$, $d(\Delta_{3})>4$ and $d(\Delta_{5})>4$ otherwise there is a contradiction to $H$ non-cyclic or $|c|=3$ or $|d|=3$ so add $\frac{1}{3} c(\Delta)=\frac{\pi}{6}$ to $c(\Delta_{1})$, $c(\Delta_{3})$ and $c(\Delta_{5})$.

Observe that if $\hat{\Delta}$ receives positive curvature then it does so across the edges $b^{-1}c$, $ad^{-1}$, $bd$, $ca$ and there is a two-thirds gap in $bd$ and $ca$. If there is only one occurrence of $c$ or $d$ in $l(\hat{\Delta})$ then this contradicts $H$ non-cyclic so assume otherwise. Also if $|c^{-1}d|=2$ then we obtain a sphere as by Lemma 3.1(c)(iv).
So assume that $|c^{-1}d|>2$.
Checking now shows that there is at least one gap between any two occurrences of $c$ or between any two occurrences of $d$; and there is at least a two-thirds gap between any two occurrences of $c$ and $d$.  It follows that if there are at least four occurrences of $c$ or of $d$ then
$c^{\ast} (\hat{\Delta}) \leq 0$.  Suppose there are three occurrences of $c$.  If there are three occurrences of $d$ then there are at least
$6. \frac{2}{3} = 4$ gaps; and if there are two occurrences of $d$ then there are four gaps or $l(\hat{\Delta})$ yields the relator
$c^2 d^{-1} cd^{-1}$ which forces $H$ to be cyclic.  Suppose there are two occurrences of $c$.  If there are three occurrences of $d$ then there are four gaps or $l(\hat{\Delta})$ yields the relator $c^{-1} dc^{-1} d^2$ which forces $H$ to be cyclic; and if there are two occurrences of $d$ then $l(\hat{\Delta})$ forces $| c^{-1} d | = 2$ or $H$ cyclic or is one of
$c\{a,b\} d \{a^{-1},b\} c \{a,b\} d \{a^{-1},b^{-1}\}$, $c\{a,b\} d \{a^{-1},b^{-1}\} c^{-1} \{a,b\} d^{-1} \{a^{-1},b^{-1}\}$ and there are four gaps.
This leaves the case when $c$ does not occur.
Checking $l(\hat{\Delta})$ shows that there is a contradiction to $|d|=3$ or there are at least four gaps except when $l(\hat{\Delta})=bddda^{-1}$ as shown in Figure 5.12(iv).
In this case add $c^{*}(\hat{\Delta})=\frac{\pi}{6}$ to $c(\hat{\Delta}_{1})$ as shown. Since $\frac{\pi}{6}$ is distributed across $ca$ it follows that $c^{*}(\hat{\Delta}_{1})\le 0$ as above.

In conclusion $\mathcal{P}$ is aspherical except when $H$ is cyclic or when $H$ is non-cyclic and $| cd | = 2$.

It follows from (\textbf{B1})--(\textbf{B12}) that either $\mathcal{P}$ is aspherical or modulo $T$-equivalence one of the conditions of Theorem 1.1 (i)--(iii) or Theorem 1.2 (i), (ii), (iv)--(x) is satisfied and so Theorems 1.1 and 1.2 are proved for Case B.




\begin{center}
Abd Ghafur Bin Ahmad\\
School of Mathematical Sciences\\
Faculty of Sciences and Technology\\
Universiti Kebangsaan Malaysia\\
43600 UKM, Bangi\\
Selangor Darul Ehsan, Malaysia\\
ghafur@ukm.my

\bigskip

Muna A Al-Mulla\footnote{Research supported by Ministry of Education, Kingdom of Saudi Arabia}\\
Department of Mathematics and Statistics\\
King Faisal University\\
PO Box 400\\
Alahssa 31982, Kingdom of Saudi Arabia\\
malmulla@fku.edu.sa

\bigskip

Martin Edjvet\\
School of Mathematical Sciences\\
The University of Nottingham\\
University Park, Nottingham NG7 2RD, UK\\
martin.edjvet@nottingham.ac.uk
\end{center}


\begin{thebibliography}{99}
\bibitem{}
Y G Baik, W A Bogley and S J Pride, On the asphericity of length four relative group presentations, Int.\ J.\ Algebra Comput.\ \textbf{7} (1997), 277--312.
\bibitem{}
W A Bogley and S J Pride, Aspherical relative presentations, Proc.\ Edinburgh Math.\ Soc.\ \textbf{35} (1992), 1--39.
\bibitem{}
D E Cohen, Combinatorial group theory: a topological approach, LMS Student Texts \textbf{14} (CUP) (1989).
\bibitem{}
D J Collins and J Huebschmann, Spherical diagrams and identities among relations, Math.\ Ann.\ \textbf{261} (1982), 155--183.
\bibitem{}
P J Davidson, On the asphericity of a family of relative group presentations, Int.\ J.\ Algebra Comput.\ \textbf{19} (2009), 159--185.
\bibitem{}
M Edjvet, Equations over groups and a theorem of Higman, Neumann and Neumann, Proc.\ London Math.\ Soc.\ \textbf{62} (1991), 563--589.
\bibitem{}
M Edjvet, On the asphericity of one-relator relative presentations, Proc.\ Roy.\ Soc.\ Edinburgh Sec.\ \textbf{124} (1994), 713--728.
\bibitem{}
M Edjvet and J Howie, The solution of length four equations over groups, Trans.\ Amer. Math.\ Soc.\ \textbf{326} (1991), 345--369.
\bibitem{}
The GAP Group: GAP -- Groups, Algorithms and Programming; Version 4.6.4, 2013.  http://www.gap-system.org.
\bibitem{}
N D Gilbert and J Howie, LOG groups and cyclically presented groups, J.\ Algebra \textbf{174} (1995), 118--131.
\bibitem{}
J Howie and V Metaftsis, On the asphericity of length five relative group presentations, Proc.\ London Math.\ Soc.\ \textbf{82} (2001), 173--194.
\bibitem{}
V Metaftsis, On the asphericity of relative group presentations of arbitrary length, Int.\ J.\ Algebra Comput.\ \textbf{13} (2003), 323--339.
\end{thebibliography}
\end{document}